\documentclass[onefignum,onetabnum]{siamart171218}

\usepackage{bm}
\usepackage{lipsum}
\usepackage{amsfonts,amssymb}
\usepackage{graphicx,xcolor}
\usepackage{epstopdf}
\usepackage[version=4]{mhchem}
\newcommand{\mathcolortext}[2]{\!\colorbox{#1}{$\scriptstyle #2$}\!}
\usepackage{myalgpseudocode}
\DeclareMathOperator{\arcosh}{arcosh}
\renewcommand{\algorithmicrequire}{\textbf{Input:}}
\renewcommand{\algorithmicensure}{\textbf{Output:}}
\newcommand{\vertiii}[1]{{\left\vert\kern-0.25ex\left\vert\kern-0.25ex\left\vert #1 
    \right\vert\kern-0.25ex\right\vert\kern-0.25ex\right\vert}}
\ifpdf
  \DeclareGraphicsExtensions{.eps,.pdf,.png,.jpg}
\else
  \DeclareGraphicsExtensions{.eps}
\fi

\newcommand{\Expect}{\operatorname{\mathbb{E}}}

\newsiamremark{remark}{Remark}
\newsiamremark{hypothesis}{Hypothesis}
\crefname{hypothesis}{Hypothesis}{Hypotheses}
\newsiamthm{claim}{Claim}

\headers{Randomized algorithms for low-rank matrix approximation}{J. A. Tropp and R. J. Webber}

\title{Randomized algorithms for low-rank matrix approximation: Design, analysis, and applications%
\thanks{
\textbf{Funding:} JAT and RJW acknowledge partial support from the Office of Naval Research through BRC Award N00014-18-1-2363, from the National Science Foundation through FRG Award 1952777, and from Caltech through the Carver Mead New Adventures Fund.}}

\author{Joel A. Tropp\thanks{Department of Computing and Mathematical Sciences, California Institute of Technology, Pasadena, CA 91125
  (\email{jtropp@cms.caltech.edu}, \email{rwebber@caltech.edu}).}
\and Robert J. Webber\footnotemark[2]}

\usepackage{amsopn}
\DeclareMathOperator{\diag}{diag}

\makeatletter
\newcommand*{\addFileDependency}[1]{
  \typeout{(#1)}
  \@addtofilelist{#1}
  \IfFileExists{#1}{}{\typeout{No file #1.}}
}
\makeatother

\ifpdf
\hypersetup{
  pdftitle={Randomized algorithms for low-rank matrix approximation},
  pdfauthor={J. A. Tropp and R. J. Webber}
}
\fi

\begin{document}

\maketitle


\begin{abstract}
This survey explores modern approaches for computing low-rank approximations of high-dimensional matrices by means of the randomized SVD, randomized subspace iteration, and randomized block Krylov iteration.
The paper compares the procedures via theoretical analyses and numerical studies to highlight how the best choice of algorithm depends on spectral properties of the matrix and the computational resources available.

Despite superior performance for many problems, randomized block Krylov iteration has not been widely adopted in computational science.  This paper strengthens the case for the method in three ways.  First, it presents new pseudocode that can significantly reduce computational costs.  Second, it provides a new analysis that yields simple, precise, and informative error bounds.  Last, it showcases applications to challenging scientific problems, including principal component analysis for genetic data and spectral clustering for molecular dynamics data.
\end{abstract}



\section{Motivation} \label{sec:introduction}

A core problem in numerical linear algebra is to produce
a low-rank, factorized approximation of a high-dimensional
input matrix $\bm{A} \in \mathbb{R}^{L \times N}$:
\[
\begin{array}{cccc}
\bm{A} & \approx & \bm{B} & \bm{C} \\
L \times N & & L \times k & k \times N
\end{array}
\]
We think of the inner dimension $k$ as much smaller than the outer dimensions $L$ and $N$.
In this case, the factorized approximation $\bm{BC}$ has fewer degrees of freedom than the
input matrix $\bm{A}$, so the factorization is easier to store and manipulate.
Furthermore, to attain a small error,
the approximation must expose structure in the input matrix.
Particular examples of ``structure'' include the range of the input matrix
(rank-revealing QR factorization), principal components
(truncated singular value decomposition), or salient
columns (CUR or interpolative decomposition).

Low-rank approximation serves as a fundamental tool for computational science.
Application areas include
fluid dynamics~\cite{berkooz1993proper,schmid2022dynamic},
uncertainty quantification~\cite{bui2012extreme,constantine2015active},
genetics~\cite{galinsky2016fast},
climate science~\cite{wallace1992singular},
geophysics~\cite{yang2021low},
astronomical imaging~\cite{levis2021inference},
and beyond.  Indeed, low-rank approximation
is helpful whenever we need to extract dominant
patterns in data, remove unwanted noise,
or perform compression.

\begin{figure}[t]
    \centering
    \includegraphics[width = .5\textwidth]{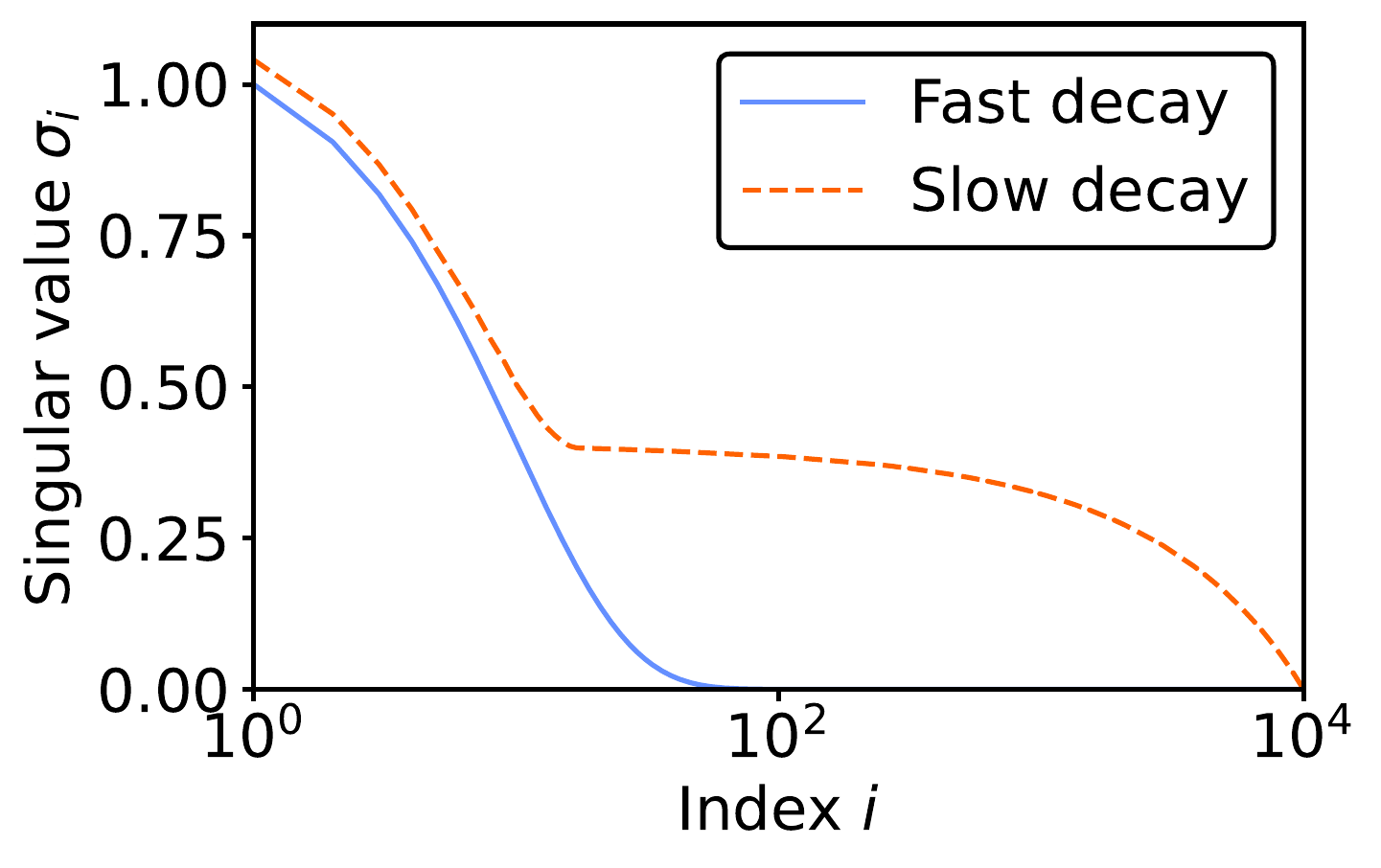}
    \caption{\textbf{(Singular value decay profiles).}  Fast versus slow singular value decay; for details of the matrices, see \cref{sec:illustrative}.}
    \label{fig:singular_decay}
\end{figure}

For high-dimensional matrices, the computational difficulty
of low-rank approximation depends on the singular value spectrum of the input matrix;
see~\cref{fig:singular_decay} for a schematic.  While we can apply
simple algorithms to approximate a matrix with a rapidly
decaying spectrum, we need more powerful algorithms
to approximate a matrix with a slowly decaying spectrum.
This challenge arises in many modern applications.  For example,
consider this quotation from the genetics literature~\cite{chen2013improved}:
\vspace{.2cm}
\begin{quote}
``[I]n large datasets, eigenvalues may be highly significant (reflecting real population structure in the data) but only slightly larger than background noise eigenvalues...''
\end{quote}
\vspace{.2cm}
In this setting, we require scalable techniques
that can filter out the noise components, while accurately
approximating the signal components.  So what algorithms
should we use?

This survey explores a family of powerful methods
for low-rank approximation: the randomized singular value
decomposition (\texttt{RSVD}), randomized subspace iteration
(\texttt{RSI}), and randomized block Krylov iteration (\texttt{RBKI}).
We will introduce these techniques in \cref{sec:hmt-framework,sec:basic-algorithms}.
They are related because each one collects information about the input matrix
using matrix--vector products with random vectors.
We will develop a systematic comparison of the three approaches
by examining their mathematical structure,
providing new implementations,
refining theoretical analyses,
and presenting numerical studies.

As compared with classic methods, the randomized
algorithms are more scalable, reliable, and robust.
Our investigation shows that the simplest method, \texttt{RSVD},
returns accurate approximations for matrices with rapid singular value decay,
but we must use the more sophisticated \texttt{RBKI} method for challenging problems.
In case storage is the limiting factor, \texttt{RSI} may offer a reasonable compromise.
These conclusions will be familiar to experts.

Nevertheless, the \texttt{RBKI} method still has not reached a
wide audience in computational science.
(For example, see the recent genetics review~\cite{tsuyuzaki2020benchmarking}.)
Therefore, the second goal of this survey is to establish firm foundations for the \texttt{RBKI}
method to speed its adoption.  We develop a new formulation of the algorithm
that is significantly faster than previous implementations.  We prove
new theoretical guarantees that demonstrate exactly how \texttt{RBKI}
improves over \texttt{RSVD} and \texttt{RSI}.  Finally,
we show applications to benchmark problems in computational science
that confirm the benefits of \texttt{RBKI} in these settings.

\begin{remark}[Finding structure with randomness] 
The paper~\cite{halko2011finding} of
Halko et al.~made a comprehensive case for
\texttt{RSVD} and \texttt{RSI}, which led to
broader usage of these methods.
The goal of this survey is to perform
the same mission for \texttt{RBKI}.
\end{remark}

\section{Randomized matrix approximation: Overview}

This section offers an introduction to the randomized
algorithms that we study in this paper.  It provides
a first look at the numerical behavior of these methods,
and it presents simple theoretical bounds that allow us
to compare their performance.

\subsection{Framework} \label{sec:hmt-framework}

In 2009, Halko et al.~\cite{halko2011finding} developed a framework
for designing randomized algorithms for low-rank matrix approximation.  In this
section, we outline the key concepts from their approach.

Consider the problem of approximating a high-dimensional matrix $\bm{A} \in \mathbb{R}^{L \times N}$.
The matrix can be an array of data stored on a computer,
or it can be an abstract linear operator defined by its action on test vectors.
Regardless, we only access $\bm{A}$ by performing multiplications $\bm{x} \mapsto \bm{A} \bm{x}$
and $\bm{y} \mapsto \bm{A}^\ast \bm{y}$ with input vectors $\bm{x} \in \mathbb{R}^N$ and $\bm{y} \in \mathbb{R}^L$.  As usual, $\bm{A}^{\ast}$ is the (conjugate) transpose.

The first idea underlying low-rank approximation is that we can find important
directions in the range of $\bm{A}$ by applying the matrix to a \emph{random} vector.
Suppose that $\bm{A}$ has the singular value decomposition (SVD) 
\begin{equation*}
    \bm{A} = \sum\nolimits_{i=1}^N \sigma_i \bm{u}_i \bm{v}_i^\ast
    \quad\text{with $\sigma_1 \geq \sigma_2 \geq \dots \geq \sigma_N \geq 0$.}
\end{equation*}
The ``leading'' left singular vectors $\bm{u}_i$,
associated with the largest singular values
$\sigma_i$, are significant directions in the range.
Draw a random vector $\bm{\omega} \in \mathbb{R}^N$ from the standard normal distribution $\mathcal{N}(\bm{0}, \mathbf{I}_N)$.
By rotational invariance, we can express $\bm{\omega}$ in the basis of right singular vectors:
\[
\bm{\omega} = \sum\nolimits_{i=1}^N Z_i \,\bm{v}_i
\quad\text{where $Z_i \overset{\mathrm{iid}}{\sim} \mathcal{N}(0,1)$.}
\]
With probability one, all of the random coefficients $Z_i$
are nonzero, and they are on the same scale because
each coefficient $Z_i$ has mean zero and variance one.
The matrix--vector product takes the form
\[
\bm{A} \bm{\omega} = \sum\nolimits_{i=1}^N Z_i \sigma_i \, \bm{u}_i.
\]
We see that the image $\bm{A}\bm{\omega}$ is strongly correlated with
the leading left singular vectors, while
it is weakly correlated with the trailing left singular vectors.

Now, observe that \emph{repeated} multiplication with the input
matrix $\bm{A}$ and its transpose $\bm{A}^\ast$ amplify the
coefficients associated with large singular values, while
attenuating coefficients with small singular values.
More precisely,
\[
(\bm{AA}^\ast)^{q-1} \bm{A} \bm{\omega} = \sum\nolimits_{i=1}^N Z_i \sigma_i^{2q-1} \, \bm{u}_i.
\]
As we increase the depth $q = 1, 2, 3, \dots$, the small singular values
are suppressed exponentially fast.  Therefore, the output vector is more
and more likely to be aligned with the leading left singular vectors.  In fact, we can achieve this goal even when the depth is a moderate constant, say, $q \leq 5$.  There is no need to take a limit.

Next, we double down on this insight by multiplying the input matrix
with a \emph{block} of $k \gg 1$ vectors to find many important
directions in the range.  To express this operation,
we introduce a standard normal matrix
$\bm{\Omega} = \begin{bmatrix} \bm{\omega}_1 & \dots & \bm{\omega}_k \end{bmatrix} \in \mathbb{R}^{N \times k}$
with independent columns.  By sequential multiplication, we obtain the
matrices
\begin{equation*}
    \bm{A} \bm{\Omega}, \quad (\bm{A} \bm{A}^\ast) \bm{A} \bm{\Omega}, \quad (\bm{A} \bm{A}^\ast)^2 \bm{A} \bm{\Omega}, \quad \ldots, \quad (\bm{A} \bm{A}^\ast)^{q-1} \bm{A} \bm{\Omega}.
\end{equation*}
As the powers increase, the ranges of these matrices
align better with the $k$ leading left singular vectors
$\bm{u}_1, \dots, \bm{u}_k$ of the input matrix.
We will quantify the extent of the overlap
in \cref{sec:theory1,sec:theory2}.

Finally, suppose we have constructed a matrix $\bm{M} \in \mathbb{R}^{L \times k}$
whose range aligns with the leading left singular vectors of the input matrix $\bm{A}$.
We can build the desired low-rank approximation of $\bm{A}$ by
\emph{compressing} $\bm{A}$ into the range of $\bm{M}$.
To that end, we orthogonalize the columns of $\bm{M}$ to obtain
a matrix $\bm{X} = \texttt{orth}(\bm{M})$,
and we construct the orthogonal projection $\bm{\Pi}_{\bm{M}} = \bm{XX}^*$
onto the range of $\bm{M}$.  Then we form the approximation
\begin{equation}
\label{eq:last_step}
    \hat{\bm{A}} = \bm{\Pi}_{\bm{M}} \bm{A} 
    = \bm{X} \bm{X}^\ast \bm{A}
    = \bm{X} \bm{Y}^\ast \qquad \text{where } \bm{Y} = \bm{A}^\ast \bm{X}.
\end{equation}
Since $\bm{X}$ has at most $k$ columns, so does the matrix $\bm{Y}$.
We have produced a factorized rank-$k$ approximation
$\hat{\bm{A}} = \bm{XY}^\ast$ of the input matrix $\bm{A}$.

At this point, the approximation $\hat{\bm{A}}$ can be manipulated
into alternative forms by standard transformations.
For example, we can easily produce an SVD,
a QR factorization, or a CUR decomposition.
See~\cite[Sec.~5]{halko2011finding} for details.

\subsection{Matrix approximations and matrix computations} \label{sec:basic-algorithms}

Combining these basic ideas, we obtain several fundamental approaches to low-rank matrix approximation.

In the first approach, called ``randomized singular value decomposition'' or \texttt{RSVD}, we generate the low-rank approximation
\begin{equation*}
\tag{\texttt{RSVD}}
    \hat{\bm{A}} = \bm{\Pi}_{\bm{A} \bm{\Omega}} \bm{A}.
\end{equation*}
This procedure requires a first multiplication with $\bm{A}$ and a second multiplication with $\bm{A}^\ast$,
while storage is limited to two blocks of $k$ vectors.
See \cref{alg:rsvd} for pseudocode.

In the second approach, called ``randomized subspace iteration'' or \texttt{RSI}, we generate
\begin{equation*}
\tag{\texttt{RSI}}
    \hat{\bm{A}} = \bm{\Pi}_{(\bm{A} \bm{A}^\ast)^{q-1} \bm{A} \bm{\Omega}} \bm{A}.
\end{equation*}
This procedure requires $2q$ multiplications, alternating between $\bm{A}$ and $\bm{A}^\ast$,
while storage is limited to two blocks of $k$ vectors.
We think about the depth parameter $q$ as a \emph{fixed} number, rather than treating
the algorithm as a limiting process.  See~\cref{alg:rsi} for pseudocode.

The third approach, called ``randomized block Krylov iteration'' or \texttt{RBKI},
forms the approximation
\begin{equation*}
\tag{\texttt{RBKI}}
    \hat{\bm{A}} = \bm{\Pi}_{[\begin{smallmatrix}
    \bm{A} \bm{\Omega} & \cdots &  (\bm{A} \bm{A}^{\ast})^{q-1} \bm{A} \bm{\Omega}
    \end{smallmatrix}]} \bm{A}.
\end{equation*}
In this case, we project onto the range of the block matrix $\begin{bmatrix}
\bm{A} \bm{\Omega} & \cdots &  (\bm{A} \bm{A}^{\ast})^{q-1} \bm{A} \bm{\Omega}
\end{bmatrix}$.
\Cref{sec:rbki} shows how to carry out this construction using $2q$ multiplications.
The storage cost is
higher than the other algorithms, typically $2q$ blocks of $k$ vectors. 
As before, the depth parameter $q$ is viewed as a fixed number.
See \cref{alg:rbki} for pseudocode, which reduces the number of matrix--vector products by 33\% compared to previous \texttt{RBKI} implementations.

Our primary objective in this work is to address the following question: 
\vspace{.2cm}
\begin{quote}
    What is the most appropriate algorithm for low-rank approximation of a given matrix: \texttt{RSVD}, \texttt{RSI}, or \texttt{RBKI}?
\end{quote}
\vspace{.2cm}
The answer depends on the singular value spectrum of the matrix, the required accuracy, and limitations
on computational resources (such as storage).

We focus on the high-dimensional setting in which $L$ or $N$ is large (say, $\geq 10^5$).
In this case, the predominant operating cost for each approximation method is the repeated application
of $\bm{A}$ and $\bm{A}^\ast$ to a block of $k$ vectors.
Therefore, we measure the operating cost in terms of the number $m$ of matrix--matrix multiplications, together with the block size $k$.
The number of matrix--vector products is $km$ for all the algorithms we consider.

Our goal is to reduce the parameters $m$ and $k$ as much as possible, while ensuring a robust and accurate low-rank approximation.  We accept certain tradeoffs between $m$ and $k$, but we focus mainly on decreasing the number $m$ of multiplications.
Reducing $m$ gives a linear improvement in the runtime because the repeated matrix multiplications are inherently serial.

\begin{remark}[Matrix multiplication]
On modern computers, matrix multiplications with a large block size $k$
are optimized to take advantage of architectural features,
including
multithreading \cite{smith2014anatomy},
caching \cite{goto2008anatomy,low2016analytical},
parallelization \cite{bui2012extreme}, and
single-instruction multiple data processing \cite{VSM11}.
By designing multiplication-rich algorithms,
we harness the full power of modern computing
platforms.  This is one of the main advantages
of block Krylov methods over Krylov methods
with a single starting vector ($k = 1$). 
See \cref{sec:block_size} for runtime comparisons.
\end{remark}

\subsection{Illustrative comparison} \label{sec:illustrative}

In this section, we present an illustrative comparison of \texttt{RSVD}, \texttt{RSI}, and \texttt{RBKI}.
\Cref{sec:experiments} presents a more comprehensive series of experiments, including real-world applications.

We consider approximating two matrices with dimensions $L = N = 10^4$:
\begin{equation*}
    \bm{A} = \textup{diag}(1, \mathrm{e}^{-.1}, \mathrm{e}^{-.2}, \ldots, \mathrm{e}^{-999.9}), \qquad \bm{B} = \bm{A} + \bm{Z}, \quad \text{where $Z_{ij} \overset{\mathrm{iid}}{\sim} \mathcal{N}(0, .002^2)$.}
\end{equation*}
The first matrix $\bm{A}$ has
exponentially decaying singular values; see the ``fast decay'' profile in \cref{fig:singular_decay}.
The second matrix $\bm{B}$ is a noisy version of $\bm{A}$ that has been perturbed
by adding an independent Gaussian random variable to each entry.
For a particular realization of the noise, the upper left submatrices are
\begin{equation*}
\bm{A} \!=\! 
\left[\!\begin{smallmatrix}
~~\mathcolortext{blue!20}{1.000} & ~~0.000 & ~~0.000 & ~~0.000 & \\
~~0.000 & ~~\mathcolortext{blue!20}{0.905} & ~~0.000 & ~~0.000 & \\
~~0.000 & ~~0.000 & ~~\mathcolortext{blue!20}{0.819} & ~~0.000 & \\
~~0.000 & ~~0.000 & ~~0.000 & ~~\mathcolortext{blue!20}{0.741} & \\
& & & & \ddots
\end{smallmatrix}\right]\!,
\quad
\bm{B} \!=\! 
\left[\!\begin{smallmatrix}
~~\mathcolortext{red!20}{1.001} & -0.002 & ~~0.002 & ~~0.002 & \\
~~0.000 & ~~\mathcolortext{red!20}{0.907} & -0.003 & -0.001 & \\
-0.003 & ~~0.002 & ~~\mathcolortext{red!20}{0.819} & -0.001 & \\
-0.002 & ~~0.000 & -0.003 & ~~\mathcolortext{red!20}{0.739} & \\
& & & & \ddots
\end{smallmatrix}\right]\!.
\end{equation*}
The matrix entries are similar, up to variations on the scale $\pm 0.002$.
However, the Gaussian noise transforms the structure of the singular values.
The matrix $\bm{B}$ has a long tail of singular values with magnitudes $0.0$ to $0.4$;
see the ``slow decay'' profile in \cref{fig:singular_decay}.

First, we exhibit the output of \texttt{RSVD} (\cref{alg:rsvd}) with a block size $k = 100$ when applied
to the original matrix $\bm{A}$:
\begin{equation*}
\tag{\texttt{RSVD}}
\hat{\bm{A}}_{k=100} =\! 
\left[\!\begin{smallmatrix}
~~\mathcolortext{blue!20}{1.000} & ~~0.000 & ~~0.000 & ~~0.000 & \\
~~0.000 & ~~\mathcolortext{blue!20}{0.905} & ~~0.000 & ~~0.000 & \\
~~0.000 & ~~0.000 & ~~\mathcolortext{blue!20}{0.819} & ~~0.000 & \\
~~0.000 & ~~0.000 & ~~0.000 & ~~\mathcolortext{blue!20}{0.741} & \\
& & & & \ddots
\end{smallmatrix}\right]\!.
\end{equation*}
The upper left submatrix coincides with $\bm{A}$ up to three decimal places.
\texttt{RSVD} performs ideally for this example, and there is no reason
to increase the block size or use a more elaborate approximation.

In contrast, \texttt{RSVD} with a block size $k = 100$
produces an appalling approximation of the noisy matrix $\bm{B}$,
with major discrepancies in the first four rows and columns:
\begin{equation*}
\tag{\texttt{RSVD}}
\hat{\bm{B}}_{k=100} =\! 
\left[\!\begin{smallmatrix}
~~\mathcolortext{red!20}{0.270} & -0.018 & ~~0.024 & -0.006 & \\
-0.019 & ~~\mathcolortext{red!20}{0.159} & ~~0.021 & -0.007 & \\
~~0.029 & ~~0.024 & ~~\mathcolortext{red!20}{0.106} & -0.015 & \\
-0.008 & -0.008 & -0.018 & ~~\mathcolortext{red!20}{0.086} & \\
& & & & \ddots
\end{smallmatrix}\right]\!.
\end{equation*}
We can improve the approximation quality by using a larger block size $k = 1{,}000$, but the errors remains large:
\begin{equation*}
\tag{\texttt{RSVD}}
\hat{\bm{B}}_{k = 1000} =\! 
\left[\!\begin{smallmatrix}
~~\mathcolortext{red!20}{0.779} & -0.001 & ~~0.004 & ~~0.006 & \\
~~0.002 & ~~\mathcolortext{red!20}{0.655} & ~~0.000 & -0.004 & \\
~~0.003 & ~~0.005 & ~~\mathcolortext{red!20}{0.569} & -0.006 & \\
~~0.004 & -0.005 & -0.009 & ~~\mathcolortext{red!20}{0.487} & \\
& & & & \ddots
\end{smallmatrix}\right]\!.
\end{equation*}
The off-diagonal entries are too large in magnitude, and the diagonal entries are too small.

As an alternative to \texttt{RSVD}, we can approximate the noisy matrix $\bm{B}$ using \texttt{RSI} (\cref{alg:rsi2}) or \texttt{RBKI} (\cref{alg:rbki2}) with a block size $k = 100$.  After $m = 5$ multiplications, \texttt{RSI} improves on \texttt{RSVD},
but it still produces an underestimate of the diagonal:
\begin{equation*}
\tag{\texttt{RSI}}
\hat{\bm{B}}_{\textup{RSI}} =\! 
\left[\!\begin{smallmatrix}
~~\mathcolortext{red!20}{0.997} & -0.002 & ~~0.001 & ~~0.001 & \\
~~0.001 & ~~\mathcolortext{red!20}{0.899} & -0.002 & ~~0.000 & \\
-0.003 & ~~0.002 & ~~\mathcolortext{red!20}{0.806} & -0.001 & \\
-0.003 & ~~0.000 & -0.002 & ~~\mathcolortext{red!20}{0.723} & \\
& & & & \ddots
\end{smallmatrix}\right].
\end{equation*}
In contrast, after $m = 5$ multiplications,
\texttt{RBKI} yields the superior approximation
\begin{equation*}
\tag{\texttt{RBKI}}
\hat{\bm{B}}_{\textup{RBKI}} =\! 
\left[\!\begin{smallmatrix}
~~\mathcolortext{red!20}{0.999} & -0.002 & ~~0.001 & ~~0.002 & \\
~~0.000 & ~~\mathcolortext{red!20}{0.905} & -0.003 & -0.001 & \\
-0.003 & ~~0.002 & ~~\mathcolortext{red!20}{0.816} & -0.001 & \\
-0.002 & ~~0.000 & -0.003 & ~~\mathcolortext{red!20}{0.736} & \\
& & & & \ddots
\end{smallmatrix}\right]\!.
\end{equation*}
In fact, the approximation $\hat{\bm{B}}_{\textup{RBKI}}$ agrees
with the best rank-$100$ approximation of $\bm{B}$ up to
three decimal places.  The remaining error is a consequence
of the low-rank approximation, rather than a deficiency in
the \texttt{RBKI} algorithm.

The takeaway from these simple experiments is that it can be difficult to approximate a matrix such as $\bm{B}$ that has been corrupted by noise.
For approximation problems with noisy matrices, \texttt{RBKI} typically offers the best combination of speed and accuracy, and it outperforms the other randomized low-rank approximation algorithms.

\subsection{Theoretical analysis}

In this section, we introduce our main error bounds for randomized low-rank matrix approximation, which help to explain the outcomes of the experiments in \cref{sec:illustrative}. Detailed derivations and additional bounds can be found in \cref{sec:theory1,sec:theory2}.

Return to the setting where $\bm{A} \in \mathbb{R}^{L \times N}$ is an arbitrary matrix.
We measure the quality of a computed approximation by comparison with an optimal rank-$r$ approximation.  For a parameter $r \geq 1$, the minimal rank-$r$ approximation error in the spectral norm is given by
the $(r + 1)$st singular value:
\begin{equation*}
\min \{ \lVert \bm{B} - \bm{A} \rVert : \operatorname{rank}(\bm{B}) = r \}
    = \sigma_{r+1}(\bm{A}).
\end{equation*}
Here, $\lVert \cdot \rVert$ is the spectral norm, and
$\sigma_{j}(\cdot)$ returns the $j$th largest singular value. 
Equality holds for any $r$-truncated SVD of the matrix $\bm{A}$.

A randomized low-rank approximation algorithm produces a random approximation $\hat{\bm{A}}$
of the input matrix $\bm{A}$.  
We evaluate the computed approximation
$\hat{\bm{A}}$ using the mean-square relative error:
\begin{equation*}
\Expect \Biggl[ \frac{\lVert  \bm{A} - \hat{\bm{A}} 
\rVert^2}{\sigma_{r+1}(\bm{A})^2} \Biggr].
\end{equation*}
The expectation $\Expect$ averages over the randomness in the algorithm, which comes from the choice of the random initialization matrix $\bm{\Omega}$. When we compare an approximation with rank $k$ to a best approximation with rank $r < k$, the relative error could be smaller than one.

Here is the core technical question:
\vspace{.2cm}
\begin{quote}
For a given target rank $r$, what block size $k$ and number $m$ of matrix multiplications suffice to make the mean-square relative error small?
\end{quote}
\vspace{.2cm}
Our theoretical analysis elucidates how the parameters $k, m$ and the
choice of approximation subspace affect the quality of the approximation.
Let us emphasize that
the user of the algorithm specifies $k$ and $m$, while the comparison rank $r$ only appears in the theoretical analysis.

As a first result, we present an error bound for \texttt{RSVD} with $k$
Gaussian vectors.  For each block size $k \geq r + 2$,
we have the estimate
\begin{equation*}
    \tag{\texttt{RSVD}}
    \frac{\Expect \lVert \bm{A}  - \hat{\bm{A}}\rVert^2}{\sigma_{r + 1}(\bm{A})^2} \leq 1
    + \frac{r}{k - r - 1} \sum_{i > r} \frac{\sigma_i(\bm{A})^2}{\sigma_{r + 1}(\bm{A})^2}.
\end{equation*}
As the block size $k$ increases, the polynomial factor $r/(k-r-1)$ decreases.
Typical choices for the block size are $k = r + 2$ and $k = 2r + 1$.
Once $k$ and $r$ are fixed, the bound depends only on the tail singular values
$\sigma_{r + 1}(\bm{A}), \sigma_{r + 2}(\bm{A}), \ldots$.
When these tail singular values decay quickly, the error is proportional to the optimal rank-$r$ error. 
For example, when the singular values of $\bm{A}$ decay exponentially fast, the error satisfies $\Expect \lVert \bm{A} - \hat{\bm{A}} \rVert^2 \leq (1 + \varepsilon) \cdot \sigma_{r+1}(\bm{A})^2$ for $k \geq r + \mathrm{const} \cdot \log (1 + r/\varepsilon)$.

The next result shows that randomized subspace iteration, \texttt{RSI},
is a better alternative for matrices with slowly decaying singular values.
Consider the approximation $\hat{\bm{A}}$ produced by the \texttt{RSI}
algorithm with $k$ Gaussian vectors and $m$ matrix
multiplications.  For each block size $k \geq r + 2$,
\begin{equation*}
    \tag{\texttt{RSI}}
    \log\biggl(\frac{\Expect \lVert \bm{A} - \hat{\bm{A}} \rVert^2}{\sigma_{r+1}(\bm{A})^2}\biggr)
    \leq \frac{1}{m - 1} \log\biggl(1 + \frac{r}{k - r - 1}
    \sum_{i > r} \frac{\sigma_i(\bm{A})^{2}}
    {\sigma_{r+1}(\bm{A})^{2}}
    \biggr).
\end{equation*}
The \emph{logarithm} of the mean-square relative error decreases in proportion
to the number $m$ of multiplications. 
The logarithm ensures that the error cannot depend too strongly on the singular value decay.
For example, the error always satisfies $\Expect \lVert \bm{A} - \hat{\bm{A}} \rVert^2 \leq \rm{e} \cdot \sigma_{r+1}(\bm{A})^2$ for $k = 2r + 1$ and
$m \geq \log \bigl(1 + \min\{L, N\} \bigr) + 1$.

Now, consider the approximation $\hat{\bm{A}}$ produced by the
\texttt{RBKI} algorithm with $k$ Gaussian initialization vectors
and $m$ matrix multiplications.  For each block size $k \geq r + 2$,
\begin{equation*}
\tag{\texttt{RBKI}}
    \log\biggl(\frac{\Expect \lVert \bm{A} - \hat{\bm{A}} \rVert^2}{\sigma_{r+1}(\bm{A})^2}\biggr)
    \leq \frac{1}{4(m - 2)^2} \biggl[\log\biggl(4 + \frac{4r}{k - r - 1}
    \sum_{i > r} \frac{\sigma_i(\bm{A})^2}{\sigma_{r+1}(\bm{A})^2}
    \biggr)\biggr]^2.
\end{equation*}
This bound shows that \texttt{RBKI} reduces the logarithm
of the relative error proportionally to the \emph{square} $m^2$
of the number $m$ of multiplications, which suggests that
\texttt{RBKI} can provide a far more accurate approximation than \texttt{RSI}.
For example, we might require $m = 100$ multiplications to obtain an accurate approximation using $\texttt{RSI}$ but only $m = 10$ multiplications using $\texttt{RBKI}$.

As stated, these error bounds are new, but they build on prior work.
The results for \texttt{RSVD} and \texttt{RSI} are patterned on
arguments from Halko et al.~\cite[Thm.~10.6, Cor.~10.10]{halko2011finding}.
Musco and Musco~\cite[Thm.~1]{musco2015randomized} obtained
a qualitative result for \texttt{RBKI} that identifies
the $\mathcal{O}(m^{-2})$ scaling, but our detailed
analysis of \texttt{RBKI} does not have a precedent
in the literature.
For mathematical derivations
and a more comprehensive discussion,
see \cref{sec:theory1,sec:theory2}.

\subsection{Positive-semidefinite matrices}

A special priority of this survey is to develop methods for approximating \emph{positive-semidefinite} (psd) matrices.
Of course, we can approximate a psd matrix using the general-purpose algorithms \texttt{RSVD}, \texttt{RSI}, and \texttt{RBKI},
but this strategy is wasteful.
A more accurate alternative that also preserves the psd property is to employ algorithms based on Nystr\"om approximation (discussed in \cref{sec:randomized_nystrom}).

We discuss three methods: \texttt{NysSVD}, \texttt{NysSI}, and \texttt{NysBKI},
which are the Nystr\"om-based extensions of \texttt{RSVD}, \texttt{RSI}, or \texttt{RBKI}.
We address the following question:
\vspace{.2cm}
\begin{quote}
    What is the most appropriate algorithm for low-rank approximation of a given psd matrix: \texttt{NysSVD}, \texttt{NysSI}, or \texttt{NysBKI}?
\end{quote}
\vspace{.2cm}
Our results show that \texttt{NysSVD} is the most efficient algorithm for rapidly decaying eigenvalues, while \texttt{NysBKI} is the most efficient algorithm for slowly decaying eigenvalues.
Additionally, our theory and experiments suggest that \texttt{NysBKI} improves on \texttt{RBKI} by a factor of $\sqrt{2}$ matrix--matrix multiplications, demonstrating the potential for speedups in the psd setting.

The \texttt{NysRBKI} algorithm is new, as is most of the theory for the Nystr{\"o}m methods.
For a real-world application, see the spectral clustering problem in \cref{sec:clustering}.

\subsection{Plan for the paper}
The rest of the paper is organized as follows.
\Cref{sec:history} presents historical background,
\cref{sec:design} defines goals of low-rank matrix approximation,
\cref{sec:pseudocode} presents pseudocode,
\cref{sec:parameter_choices} discusses parameter choices,
\cref{sec:experiments} showcases numerical experiments, \cref{sec:theory1,sec:theory2} derive error bounds, and \cref{sec:conclusion} offers conclusions.

\subsection{Notation}

For simplicity, we work with real-valued matrices only.  There is an analogous theory
for complex-valued matrices, although the constants may differ.  We use the term
\emph{Gaussian} random matrix (or vector) as a shorthand for a matrix with independent
Gaussian entries with mean zero and variance one.

We adopt standard notation for various matrix operations.
The transpose of a matrix $\bm{A} \in \mathbb{R}^{L \times N}$ is represented by $\bm{A}^\ast$.
The Moore--Penrose pseudoinverse of $\bm{A}$ is denoted by $\bm{A}^{\dagger}$, and the orthogonal projection onto the range space of $\bm{A}$ is written as $\bm{\Pi}_{\bm{A}} = \bm{A} (\bm{A}^\ast \bm{A})^{\dagger} \bm{A}^\ast$.
The eigenvalues of a positive-semidefinite matrix $\bm{A}$ are ordered from largest to smallest as $\lambda_1(\bm{A}) \geq \lambda_2(\bm{A}) \geq \cdots$, as are the singular values $\sigma_1(\bm{A}) \geq \sigma_2(\bm{A}) \geq \dots$ of a general matrix $\bm{A}$.
The symbol $\preccurlyeq$ denotes the positive-semidefinite (psd) order on symmetric matrices: $\bm{A} \preccurlyeq \bm{B}$ if and only if $\bm{B} - \bm{A}$ is psd.

We measure the size of a matrix $\bm{M} \in \mathbb{R}^{L \times N}$ using the Schatten $p$-norm \cite[Sec.~4.2]{bhatia2013matrix}, defined by
\begin{equation*}
\lVert \bm{M} \rVert_p = \Bigl(\sum\nolimits_{i=1}^N \sigma_i(\bm{M})^p\Bigr)^{1 / p}
\end{equation*}
for $1 \leq p < \infty$ and $\lVert \bm{M} \rVert_{\infty} = \sigma_1(\bm{M})$.
The most commonly used Schatten $p$-norms are the
trace norm
$\lVert \cdot \rVert_1 = \lVert \cdot \rVert_\ast$, the
Frobenius norm
$\lVert \cdot \rVert_2 = \lVert \cdot \rVert_{\rm F}$, and the spectral norm $\lVert \cdot \rVert_{\infty} = \lVert \cdot \rVert$, each with its distinctive notation.
Schatten $p$-norms are unitarily invariant, invariant under (conjugate) transposition, and satisfy the matrix inequality $\lVert \bm{A} \bm{M} \bm{B} \rVert_p \leq \lVert \bm{A} \rVert \lVert \bm{M} \rVert_p \lVert \bm{B} \rVert$.

In addition to the standard conventions, we introduce the following special definitions.
A matrix is said to be rank-$r$ when its rank does not exceed $r$.
We use $\lfloor\bm{A}\rfloor_r$ to denote any optimal rank-$r$ approximation of the matrix $\bm{A} \in \mathbb{R}^{L \times N}$, derived from an $r$-truncated singular value decomposition:
\[
\bm{A} = \sum\nolimits_{i=1}^N \sigma_i \bm{u}_i \bm{v}_i^\ast
\quad\text{yields}\quad
\lfloor \bm{A} \rfloor_r = \sum\nolimits_{i=1}^r \sigma_i \bm{u}_i \bm{v}_i^\ast.
\]
This low-rank approximation may not be unique, so we employ the notation only in contexts where it leads to an unambiguous statement.
We represent the Nystr\"om approximation of a psd matrix $\bm{A} \in \mathbb{R}^{N \times N}$ with respect to a test matrix $\bm{X} \in \mathbb{R}^{N \times k}$ as $\bm{A}\langle \bm{X} \rangle = \bm{A} \bm{X} (\bm{X}^{\ast} \bm{A} \bm{X})^{\dagger} \bm{X}^{\ast} \bm{A}$.

\section{History} \label{sec:history}

In this section, we provide a concise history of randomized low-rank matrix approximation.

\subsection{Numerical methods for spectral computation}

Computational tools for low-rank approximation have their
roots in numerical methods for spectral computation.  We focus on methods that extract a few eigenvalues and eigenvectors, rather than returning a full spectral decomposition.

In the early 20th century, the power method and variants, such as inverse iteration, emerged
as popular techniques for computing a few eigenvalues and eigenvectors of
a square matrix (see \cite[Sec.~2.1.3]{stewart2001matrix}).
Around 1950, Lanczos~\cite{lanczos1950iteration} introduced
an improvement of the power method for Hermitian matrices
that fully exploits Krylov information (matrix powers applied to vectors).

Both the power method and the Lanczos algorithm employ just a single starting vector (i.e., block size $k = 1$),
making them inadequate for resolving eigenvalues with multiplicity higher than one.
To address this shortcoming, between the 1950s and 1970s, computational mathematicians
developed versions of the power method and the Lanczos algorithm that use multiple starting vectors,
called ``subspace iteration'' and ``block Krylov iteration'' (see \cite[Secs.~5.3.4, 6.1.4]{stewart2001matrix}).

These classic iterative methods have been extensively discussed in the literature;
see Golub and van der Vorst~\cite{golub2000eigenvalue} for an historical overview.
For modern accounts of the power method and the Lanczos algorithm, refer to the textbooks of
Partlett \cite{parlett1998symmetric} and Saad \cite{saad2011numerical}.
For analyses of subspace iteration and block Krylov iteration, see Stewart \cite{stewart1976simultaneous} and Li and Zhang \cite{li2015convergence}.

All the traditional iterative methods for calculating the dominant eigenvalues and eigenvectors of a Hermitian matrix can be
adapted to calculate the dominant singular values and singular vectors of a rectangular matrix $\bm{A} \in \mathbb{R}^{L \times N}$ by applying them to the Jordan--Wielandt matrix $\bigl[\begin{smallmatrix} \bm{0} & \bm{A} \\ \bm{A}^\ast & \bm{0} \end{smallmatrix}\bigr]$; see \cite{GoKa65,GLO81}.
When initialized with $\bm{A} \bm{\Omega}$ where $\bm{\Omega}$ is a random matrix, this approach matches the approximate singular value decomposition $\hat{\bm{A}} = \bm{U} \bm{\Sigma} \bm{V}^\ast$ from \texttt{RSI} or \texttt{RBKI}.
However, the iterative approach to singular value decomposition was historically used to calculate a few (say, 1--4) leading singular values and singular vectors \cite{GLO81}.
The potential for converting these schemes into algorithms for low-rank approximation apparently went unrecognized until the 21st century.

Early iterative algorithms for eigenvalue or singular vector computations used relatively
few starting vectors (say, 3 or 4), which were often generated randomly, e.g., by random Gaussians.
The random initialization was viewed as undesirable but necessary to ensure nonzero overlap
with the leading singular vectors \cite{OSV79,CMSW79}.
To obtain accurate results,
the algorithms employed a large number of matrix multiplications (10s or 100s).
Numerical analysts emphasized the benefit of these repeated multiplications, since they viewed these computations as limiting processes, rather than finite algorithms \cite[Sec.~13.2]{parlett1998symmetric}.

\subsection{The benefits of random sampling}

The idea to calculate a singular value decomposition using just a few matrix multiplications was not seriously considered in the literature on classical iterative algorithm.
Computational scientists needed a new probabilistic perspective to bring this computational strategy to light.

Dixon (1983) \cite{dixon1983estimating} and Kuczy\'nski and Wo\'zniakowski (1992) \cite{kuczynski1992estimating} introduced a fresh analysis by applying probability theory to analyze the power method
and the Lanczos algorithm.
When the starting vector is Gaussian, these authors proved that the iterative algorithms can approximate the largest eigenvalue of a psd matrix after a \emph{fixed} number of steps that depends logarithmically on the matrix dimension.  They also demonstrated that the algorithms succeed even when the maximum eigenvalue is not separated from the remaining eigenvalues.

In the late 1990s and early 2000s, theoretical computer scientists began to analyze a different type of randomized low-rank approximation, through an approach called column sampling \cite{frieze1998fast,DFK+99,DV06,DMM06}.
The simplest column sampling approximation takes the form $\hat{\bm{A}} = \bm{\Pi}_{\bm{C}} \bm{A}$, where $\bm{C}$ is a random subset of the columns of $\bm{A}$ \cite{DMM06}, chosen from an appropriate probability distribution.
The probabilities can be defined proportionally to the square column norms \cite{frieze1998fast,DFK+99}, or using a more complicated distribution based on adaptive sampling \cite{DV06} or subspace sampling \cite{DMM06}.
Additionally as a post-processing procedure, many column sampling algorithms apply an $r$-truncated singular value decomposition to produce the rank-$r$ approximation $\hat{\bm{A}} = \lfloor \bm{\Pi}_{\bm{C}} \bm{A} \rfloor_r$.
Various analyses \cite{drineas2006fast2,rudelson2007sampling,DV06,BMD09} demonstrate that column sampling can lead to a more accurate low-rank approximation than deterministic algorithms based on rank-revealing QR \cite{Gol65,GuEi96}.
For a more detailed history of column-sampling approaches, see the survey \cite[Sec.~2.1.2]{halko2011finding}.

Around the same time, researchers began investigating randomized low-rank approximations similar to \texttt{RSVD} and \texttt{RSI}.
Papadimitriou et al.~(1998) \cite{PTRV98,papadimitriou2000latent}
introduced a low-rank approximation $\hat{\bm{A}} = \bm{\Pi}_{\lfloor \bm{A} \bm{\Omega}\rfloor_r} \bm{A}$, where $\bm{\Omega} \in \mathbb{R}^{N \times k}$ is a random orthogonal matrix and $\lfloor \bm{A} \bm{\Omega} \rfloor_r$ is a $r$-truncated singular value decomposition of $\bm{A} \bm{\Omega}$.
Sarl{\'os} (2006) \cite{Sar06} proposed a similar approximation $\hat{\bm{A}} = \lfloor \bm{\Pi}_{\bm{A} \bm{\Omega}} \bm{A} \rfloor_r$, where $\bm{\Omega} \in \mathbb{R}^{N \times k}$ is a subsampled trigonometric transform.
Rokhlin et al. (2010) \cite{rokhlin2010randomized} acknowledged the difficulty of approximating a matrix with slowly decaying singular values, and they proposed a more accurate low-rank approximation $\bm{\Pi}_{\lfloor (\bm{A} \bm{A}^\ast)^{q-1} \bm{A} \bm{\Omega}\rfloor_r} \bm{A}$, where $\bm{\Omega} \in \mathbb{R}^{N \times k}$ is a Gaussian matrix and $q \geq 1$ is a depth parameter; see also the 1997 paper of Roweis \cite{Row97}.
These early works argue that a small number of multiplications with a large block of random vectors (typically, $10 \leq k \leq 1{,}000$) is sufficient to approximate a matrix with rapidly decaying singular values,
and they provide preliminary theory
to support their argument.

In 2008 and 2009, Halko, Martinsson, and Tropp \cite{halko2011finding} developed a framework for designing
randomized low-rank approximation algorithms, which leads to the standard versions of \texttt{RSVD} and \texttt{RSI} that treat the the rank-$r$ truncation as an optional last step.
Halko et al. also obtained the first
theoretical analysis for these algorithms that is precise enough to predict their empirical behavior.  They showed that \texttt{RSI} yields provable guarantees after a fixed number of matrix multiplications,
even when the input matrix does not have gaps between the singular values.
Following the publication of \cite{halko2011finding}, researchers from various fields have applied randomized low-rank approximation to large-scale matrix computations in spectral clustering \cite{boutsidis2015spectral}, genetics \cite{abraham2014fast,galinsky2016fast}, uncertainty quantification \cite{bui2012extreme,isaac2015scalable}, and image processing \cite{cheng2016fast,kumar2019target,yang2021low}, among other areas.

\subsection{Proposals for more efficient methods} \label{sec:proposals}

In the next stage of inquiry, researchers proposed ways to improve the efficiency of randomized low-rank approximation, going beyond the basic \texttt{RSVD} and \texttt{RSI}.
First, they introduced new approaches for the low-rank approximation of psd matrices based on \emph{Nystr\"om approximation} \cite{halko2011finding,li2017algorithm,tropp2017fixed}.
For a mathematical description of Nystr\"om approximation, see \cref{sec:randomized_nystrom} or \cite[Sec.~14]{martinsson2020randomized}.
Nystr\"om approximation can be applied with any choice of test vectors, including random Gaussian vectors \cite{tropp2017fixed} or random coordinate vectors \cite{williams2001using}.
The approximation is always psd, and it always yields more accurate results than the simpler approximation based on orthogonal projection (see \cref{lem:nystrom}).

In another development, researchers improved the theoretical understanding of structured random matrices $\bm{\Omega} \in \mathbb{R}^{N \times k}$.
Random matrices constructed using sparse sign vectors \cite{BDN15,cohen2016nearly} or subsampled trigonometric transforms \cite{tropp2011improved,Foho13,PiWa15} can be used for initializing randomized low-rank approximation, and they ensure a high quality of approximation.
Due to fast matrix multiplications, these structured initialization matrices significantly improve runtimes for algorithms that require just a single matrix--matrix multiplication, such as \texttt{NysSVD} (\cref{alg:nyssvd}) or one-pass \texttt{RSVD} algorithms \cite{woolfe2008fast,ClWo09,tropp2017practical,TYUC19}.
However, when low-rank matrix approximation requires repeatedly multiplying the iterates with $\bm{A}$ and $\bm{A}^\ast$, as in \texttt{RSVD} or \texttt{RBKI}, the computational benefits of a structured random initialization are limited due to loss of structure after the first multiplication.

Last, one of the most significant advances in low-rank matrix approximation has been the introduction \cite{rokhlin2010randomized,halko2011algorithm} and analysis \cite{musco2015randomized,wang2015improved,drineas2018structural,tropp2018analysis,tropp2022randomized,BCW22} of randomized block Krylov methods.
In 2015, Musco and Musco \cite{musco2015randomized} proved that \texttt{RBKI} can be substantially more accurate than \texttt{RSI}, 
and their insight has spurred a variety of follow-up works \cite{wang2015improved,drineas2018structural,tropp2018analysis,tropp2022randomized,BCW22,meyer2023unreasonable,bakshi2023krylov}.
In spite of this progress, the literature still lacks a crisp treatment
of \texttt{RBKI} that parallels the earlier work on \texttt{RSVD} and \texttt{RBKI}.
As a consequence, while \texttt{RBKI} is well-known among experts
in randomized numerical linear algebra, it remains under-utilized in computational science
(for instance, see the recent survey article \cite{tsuyuzaki2020benchmarking} in genetics).
Therefore, our objective in this work is to increase the use of \texttt{RBKI}
by providing efficient pseudocode, detailed theory, and numerical experiments
that highlight the benefits of this algorithm.

\section{Goals and consequences} \label{sec:design}

In this section, we identify the goals and consequences of randomized low-rank matrix approximation.
We answer the questions, ``What is the purpose of randomized low-rank approximation?''
and ``What would an ideal low-rank approximation algorithm accomplish?''

\subsection{Goals} \label{sec:goals}

Ideally, a low-rank approximation algorithm should satisfy three main criteria: (1) speed, (2) accuracy, and (3) utility for downstream matrix computations.

The first design criterion is speed.
Randomized low-rank approximation is intended to quickly extract structure
from a high-dimensional matrix $\bm{A} \in \mathbb{R}^{L \times N}$.
To that end, we develop multiplication-rich algorithms, which harness the efficiency of matrix multiplications on modern computing platforms. 
We perform most of the computations in a sequence of $m$ matrix multiplications between $\bm{A}$ or $\bm{A}^\ast$ and a block of $k$ vectors.
Even with dense matrix--matrix multiplication, this approach costs only $\mathcal{O}(kmLN)$ arithmetic operations, so it can be much faster than a traditional (full) singular value decomposition or eigendecomposition, which expends $\mathcal{O}(LN \min\{L, N\})$ operations.
See \cref{sec:experiments} for examples of computational speedups spanning 1 to 3 orders of magnitude.

The second design criterion is accuracy.  Our goal is to produce an approximation $\hat{\bm{A}}$ that
is competitive with the best rank-$r$ approximation in spectral norm:
\begin{equation}
\label{eq:to_achieve}
    \lVert \bm{A} - \hat{\bm{A}} \rVert \approx \lVert \bm{A} - \lfloor \bm{A} \rfloor_r \rVert.
\end{equation}
We allow the rank of the approximation $\hat{\bm{A}}$ to be slightly larger than the comparison rank $r$.
Our experiments (\cref{sec:experiments}) and theory (\cref{sec:theory1,sec:theory2}) show
how randomized low-rank approximation achieves the goal  \cref{eq:to_achieve}.

The third design criterion is utility for downstream matrix calculations.
We focus on algorithms that return a factorized approximation, which can be stored and manipulated efficiently.
In the general setting, we seek a
singular value decomposition
$\hat{\bm{A}} = \bm{U} \bm{\Sigma} \bm{V}^{\ast}$.
In the psd setting, we instead seek an eigenvalue decomposition $\hat{\bm{A}} = \bm{U} \bm{\Lambda} \bm{U}^{\ast}$.
The output of these algorithms identifies the singular vectors or eigenvectors of the approximation $\hat{\bm{A}}$, which serve as proxies for the singular vectors or eigenvectors of the input matrix $\bm{A}$.
Additionally, the factorized output leads to accelerated routines for computing matrix--vector products or solving regularized linear systems.

\subsection{Consequences}
\label{sec:consequences}

As we have noted, our goal is to produce an approximation
that is close to the input matrix in spectral norm, say,
$\lVert \bm{A} - \hat{\bm{A}} \rVert \leq \varepsilon$.
This type of bound allows us to substitute the approximation
in place of the input matrix in many different contexts.

\vspace{0.2cm}

\begin{enumerate}
\item
\textbf{Matrix--vector multiplications}.
We can use the low-rank approximation $\hat{\bm{A}}$ for matrix--vector products,
and the error is controlled as
\begin{equation*}
    \lVert \bm{A} \bm{v} -
    \hat{\bm{A}} \bm{v} \rVert
    \leq \lVert \bm{A} - \hat{\bm{A}} \rVert \lVert \bm{v} \rVert.
\end{equation*}
\item
\textbf{Regularized linear systems}.
If $\bm{A}$ is psd and we generate a psd Nystr\"om approximation $\hat{\bm{A}}$, we can use the approximation to solve the regularized linear system $(\bm{A} + \mu \mathbf{I}) \bm{x} = \bm{y}$ with $\mu > 0$. By the resolvent identity, the error is controlled by
\begin{align*}
    \lVert (\bm{A} + \mu \mathbf{I})^{-1} - (\hat{\bm{A}} + \mu \mathbf{I})^{-1}\rVert
    &= \lVert (\bm{A} + \mu \mathbf{I})^{-1}(\bm{A} - \hat{\bm{A}})(\hat{\bm{A}} + \mu \mathbf{I})^{-1}\rVert \\
    &\leq \frac{\lVert \bm{A} - \hat{\bm{A}} \rVert}{\mu^2}.
\end{align*}
Hence, the linear solve is guaranteed to be accurate as long as $\lVert \bm{A} - \hat{\bm{A}} \rVert \ll \mu^2$.
\item
\textbf{Singular values}.
We can use the low-rank approximation $\hat{\bm{A}}$ to compute singular values.
Weyl's inequality \cite[Thm.~4.3.1]{horn2012matrix} guarantees that the error is controlled by
\begin{equation*}
    |\sigma_i(\bm{A}) - \sigma_i(\hat{\bm{A}})| \leq \lVert \bm{A} - \hat{\bm{A}} \rVert
\end{equation*}
for each $1 \leq i \leq \min\{L, N\}$.
\item
\textbf{Singular vectors}.
The low-rank approximation $\hat{\bm{A}}$ can be used to compute singular vectors.
Wedin's theorem \cite{Wedin_1972,stewart1991perturbation} (also see \cite[Thm.~4]{orourke2018random}) shows that
\begin{equation}
\label{eq:singular_vector}
    \langle \bm{u}_i(\bm{A}), \bm{u}_i(\hat{\bm{A}}) \rangle^2
    + \langle \bm{v}_i(\bm{A}), \bm{v}_i(\hat{\bm{A}})\rangle^2
    \geq 2 - \frac{8 \lVert \bm{A} - \hat{\bm{A}} \rVert^2}{\delta_i^2(\bm{A})}.
\end{equation}
for each $1 \leq i < \min\{L, N\}$.
Here, $\bm{u}_i(\bm{M})$ denotes the $i$th left singular vector, $\bm{v}_i(\bm{M})$ denotes the $i$th right singular vector, and
$\delta_i(\bm{M}) = \min_{j \neq i} |\sigma_j(\bm{M}) - \sigma_i(\bm{M})|$
denotes the $i$th singular value gap.
As a consequence of \cref{eq:singular_vector}, the $i$th singular vectors of $\bm{A}$ and $\hat{\bm{A}}$ are closely aligned as long as $\lVert \bm{A} - \hat{\bm{A}} \rVert \ll \delta_i(\bm{A})$.
\end{enumerate}
\vspace{0.2cm}
In summary, randomized low-rank approximation can be a helpful tool for performing many fast, approximate matrix computations.

\section{Pseudocode} \label{sec:pseudocode}

In this section, we present pseudocode for all of the algorithms we study. 
Efficient versions of \texttt{RSVD}, \texttt{RSI}, \texttt{NysSVD} and \texttt{NysSI}
have already appeared in \cite{rokhlin2010randomized,halko2011finding,gu2015subspace,tropp2017fixed},
but we have developed a faster implementation of \texttt{RBKI}.
Our \texttt{RBKI} pseudocode stores and reuses matrix multiplications,
which reduces the number of matrix--vector products by roughly $33\%$.  Additionally, to the best of our knowledge,
the \texttt{NysBKI} algorithm is new,
and it provides an additional factor-of-$\sqrt{2}$ savings in the number of matrix--vector products needed to achieve a fixed accuracy.

The section is structured as follows:
\cref{sec:randomized_svd,sec:rsi,sec:rbki} provide pseudocode for \texttt{RSVD}, \texttt{RSI}, and \texttt{RBKI};
\cref{sec:randomized_nystrom} discusses Nystr\"om approximation algorithms.
and \cref{sec:further} discusses modifications to the pseudocode that reduce or eliminate orthogonalization steps.

\subsection{\texttt{RSVD}}
\label{sec:randomized_svd}

The simplest approach for randomized low-rank matrix approximation is the \emph{randomized singular value decomposition}, which was introduced in \cite{papadimitriou2000latent,Sar06} and refined in \cite{halko2011finding}.
\texttt{RSVD} has been applied to thousands of problems in
spectral clustering \cite{boutsidis2015spectral}, biology \cite{hie2019efficient,tsuyuzaki2020benchmarking}, uncertainty quantification \cite{bui2012extreme,isaac2015scalable}, and image processing \cite{cheng2016fast,kumar2019target}.
The algorithm is now implemented in scikit-learn using the command ``randomized\_svd'' with parameter ``n\_iter'' set to 0 \cite{scikit-learn} and in
Matlab using the
command ``svdsketch''
with parameter ``NumPowerIterations'' set to 0 \cite{matlab2020version}.

\begin{algorithm}[t]
\caption{\texttt{RSVD} \label{alg:rsvd}}
\begin{algorithmic}[1]
\Require Matrix $\bm{A} \in \mathbb{R}^{L \times N}$; block size $k$
\Ensure Orthogonal $\bm{U} \in \mathbb{R}^{L \times k}$,
orthogonal $\bm{V} \in \mathbb{R}^{N \times k}$, and diagonal $\bm{\Sigma} \in \mathbb{R}^{k \times k}$ such that $\bm{A} \approx \bm{U} \bm{\Sigma} \bm{V}^{\ast}$
\State Generate a random matrix $\bm{\Omega} \in \mathbb{R}^{N \times k}$
\State $\bm{X} = \bm{A} \bm{\Omega}$
\State $[\bm{X}, \sim] = \textup{qr\_econ}
(\bm{X})$ \Comment{economy-sized QR or stabilized QR \cref{eq:stabilized_qr}}
\State $\bm{Y} = \bm{A}^{\ast} \bm{X}$
\State $[\hat{\bm{U}}, \bm{\Sigma}, \bm{V}] = 
\textup{svd\_econ}(\bm{Y}^\ast)$
\Comment{economy-sized SVD factorization}
\State $\bm{U} = \bm{X} \hat{\bm{U}}$
\end{algorithmic}
\end{algorithm}

\texttt{RSVD} generates a low-rank approximation
\begin{equation*}
\tag{\texttt{RSVD}}
\hat{\bm{A}} = \bm{\Pi}_{\bm{A} \bm{\Omega}} \bm{A},
\end{equation*}
where $\bm{\Omega} \in \mathbb{R}^{N \times k}$ is a random matrix, typically Gaussian, and $k$ is a block size parameter chosen by the user, typically $10 \leq k \leq 1{,}000$.
\Cref{alg:rsvd} provides \texttt{RSVD} pseudocode, which is optimized for efficiency in three ways:
\vspace{0.2cm}
\begin{enumerate}
\item
The low-rank approximation
\begin{equation*}
    \bm{\Pi}_{\bm{A} \bm{\Omega}} \bm{A} = \bm{A} \bm{\Omega} (\bm{A} \bm{\Omega})^{\dagger} \bm{A}
\end{equation*}
is generated without forming the pseudoinverse $(\bm{A} \bm{\Omega})^{\dagger}$.
As a cheaper and stabler approach, the columns of $\bm{A} \bm{\Omega}$ are orthogonalized to produce a matrix $\bm{X} = \texttt{orth}(\bm{A} \bm{\Omega})$ and the low-rank approximation is generated using
\begin{equation*}
    \bm{\Pi}_{\bm{A} \bm{\Omega}} \bm{A} = \bm{X} (\bm{X}^{\ast} \bm{A}).
\end{equation*}
\item The orthogonalization step $\bm{X} = \texttt{orth}(\bm{A} \bm{\Omega})$ can be performed 
by taking the $\bm{Q}$ matrix from a standard economy-sized QR factorization.
However, this leads to numerical issues when the requested rank of $\hat{\bm{A}}$ exceeds the true rank of $\bm{A}$.
For such problems, we advocate a stabler approach in which we first calculate the SVD
\begin{equation*}
    \bm{A} \bm{\Omega} = \bm{U} \bm{\Sigma} \bm{V}^\ast.
\end{equation*}
Then, we identify all the indices $i = s_1, \ldots, s_r$ for which the diagonal entry $\sigma_{ii}$ exceeds a threshold proportional to the machine precision $\varepsilon_{\textup{mach}}$
and return the approximate factorization
\begin{equation}
\label{eq:stabilized_qr}
    \bm{A} \bm{\Omega} \approx \bm{Q} \bm{R}, \text{ where }
    \begin{cases}
    \bm{Q} = \begin{bmatrix} \bm{u}_{s_1} & \cdots & \bm{u}_{s_r} \end{bmatrix}, \\[4pt]
    \bm{R} = \begin{bmatrix} \sigma_{s_1} \bm{v}_{s_1} & \cdots & \sigma_{s_r} \bm{v}_{s_r} \end{bmatrix}^\ast.
    \end{cases}
\end{equation}
This is not a traditional QR factorization since $\bm{R}$ is not upper triangular, but it is a stabilized QR factorization since it reliably exposes the range of $\bm{A} \bm{\Omega}$ \cite[Sec.~5.5.8]{golub2013matrix}.
We can take $\texttt{orth}(\bm{A} \bm{\Omega})$ to be the $\bm{Q}$ factor from the stabilized QR factorization, as usual.
\item
The SVD of $\bm{X} (\bm{X}^\ast \bm{A})$ is evaluated by first computing the SVD of the wide matrix 
\begin{equation*}
\bm{X}^\ast \bm{A} = \hat{\bm{U}} \bm{\Sigma} \bm{V}^\ast
\end{equation*}
and then applying the rotation
\begin{equation*}
	\bm{X} (\bm{X}^\ast \bm{A})
	= (\bm{X} \hat{\bm{U}}) \bm{\Sigma} \bm{V}^{\ast}.
\end{equation*}
\end{enumerate}
\vspace{0.2cm}
These implementation techniques, which are mostly standard~\cite[Sec.~1.5]{halko2011finding},
lead to the fast and stable pseudocode that is presented in \cref{alg:rsvd}.

Next, we comment on the main user choice when implementing \texttt{RSVD}: choosing the block size $k$.
As a general rule, increasing the block size $k$ increases the cost of \texttt{RSVD}, but it also improves the approximation quality.
We provide a quick proof of this fact below:

\begin{lemma}[Bigger projection, smaller error] \label{lem:bigger_projection}
Consider a matrix $\bm{A} \in \mathbb{R}^{L \times N}$ 
and orthogonal projections $\bm{P}_1, \bm{P}_2 \in \mathbb{R}^{L \times L}$ such that 
$\textup{range}(\bm{P}_1) \subseteq \textup{range}(\bm{P}_2)$. Then,
\begin{equation}
\label{eq:bigger_smaller}
\lVert \bm{A} - \bm{P}_2 \bm{A} \rVert_p \leq \lVert \bm{A} - \bm{P}_1 \bm{A} \rVert_p
\end{equation}
for any Schatten $p$-norm with $1 \leq p \leq \infty$.
The same inequality holds for any unitarily invariant norm.
\end{lemma}
\begin{proof}
Because
$\textup{range}(\bm{P}_1) \subseteq \textup{range}(\bm{P}_2)$,
we have $\bm{P}_1 \preccurlyeq \bm{P}_2$, and
the complementary projectors satisfy 
$\mathbf{I} - \bm{P}_2 \preccurlyeq \mathbf{I} - \bm{P}_1$.
Since conjugation preserves the psd order,
\begin{equation*}
\bm{A}^{\ast} (\mathbf{I} - \bm{P}_2) \bm{A}
\preccurlyeq
\bm{A}^{\ast} (\mathbf{I} - \bm{P}_1) \bm{A}.
\end{equation*}
Consequently, for each $1 \leq i \leq \min\{L, N\}$,
\begin{align*}
    \sigma_i(\bm{A} - \bm{P}_2 \bm{A})^2
    &= \lambda_i(\bm{A}^{\ast} (\mathbf{I} - \bm{P}_2) \bm{A}) \\
    &\leq \lambda_i(\bm{A}^{\ast} (\mathbf{I} - \bm{P}_1) \bm{A})
    = \sigma_i(\bm{A} - \bm{P}_1 \bm{A})^2.
\end{align*}
The inequality is Weyl's monotonicity theorem \cite[Thm.~4.3.1]{horn2012matrix}.
This statement implies the $p$-norm inequality~\cref{eq:bigger_smaller} and the result for unitarily invariant norms.
\end{proof}

Since increasing the block size $k$ improves the approximation quality, the user of \texttt{RSVD} can always raise $k$ to address more difficult problems.
As an adaptive strategy, the user can increase $k$ until the Frobenius--norm error hits a target precision, i.e.,
$\lVert \bm{A} - \hat{\bm{A}} \rVert_{\rm F}^2 = \lVert \bm{A} \rVert_{\rm F}^2 - \lVert \hat{\bm{A}} \rVert_{\rm F}^2 < \varepsilon$.
This adaptive version of \texttt{RSVD} is currently implemented in the ``svdsketch'' function for Matlab \cite{matlab2020version}.
For more discussion of the block size and adaptive stopping rules, see \cref{sec:parameter_choices}.

\subsection{\texttt{RSI}} \label{sec:rsi}

The early developers of \texttt{RSVD} acknowledged that the algorithm leads to large random errors when applied to a matrix with slowly decaying singular values.
When targeting a matrix with slow singular value decay, they proposed an alternative strategy called \emph{randomized subspace iteration} \cite{rokhlin2010randomized,halko2011finding}.
\texttt{RSI} is based on the randomized low-rank approximation
\begin{equation*}
\tag{\texttt{RSI}}
\hat{\bm{A}} = \bm{\Pi}_{(\bm{A} \bm{A}^{\ast})^{q-1} \bm{A} \bm{\Omega}} \bm{A},
\end{equation*}
where $\bm{\Omega} \in \mathbb{R}^{N \times k}$ is a random matrix, typically Gaussian, and $q$ is a depth parameter, typically $2 \leq q \leq 5$.
\texttt{RSI} is implemented in scikit-learn using the command ``randomized\_svd'' \cite{scikit-learn} and in
Matlab using the
command ``svdsketch'' \cite{matlab2020version}.

\texttt{RSI} requires $q$ matrix multiplications with $\bm{A}$ and
$q$ multiplications with $\bm{A}^{\ast}$.
The algorithm with $q = 1$ is identical to \texttt{RSVD}.
Increasing the depth $q$ systematically improves the approximation quality,
but the cost grows linearly with $q$.

\begin{algorithm}[t]
\caption{\texttt{RSI}, simple version \label{alg:rsi}}
\begin{algorithmic}[1]
\Require Matrix $\bm{A} \in \mathbb{R}^{L \times N}$; block size $k$; iteration count $q$
\Ensure Orthogonal $\bm{U} \in \mathbb{R}^{L \times k}$,
orthogonal $\bm{V} \in \mathbb{R}^{N \times k}$, and diagonal $\bm{\Sigma} \in \mathbb{R}^{k \times k}$ such that $\bm{A} \approx \bm{U} \bm{\Sigma} \bm{V}^{\ast}$
\State Generate a random matrix $\bm{Y} \in \mathbb{R}^{N \times k}$
\For{$i = 1, \ldots, q$}
\State $\bm{X} = \bm{A} \bm{Y}$
\State $[\bm{X}, \sim] = \textup{qr\_econ}(\bm{X})$
\Comment{Optional: stabilized QR \cref{eq:stabilized_qr}}
\State $\bm{Y} = \bm{A}^{\ast} \bm{X}$
\EndFor
\State $[\hat{\bm{U}}, \bm{\Sigma}, \bm{V}] = 
\textup{svd\_econ}(\bm{Y}^\ast)$
\State $\bm{U} = \bm{X} \hat{\bm{U}}$
\end{algorithmic}
\end{algorithm}

We provide simple, stable pseudocode for \texttt{RSI} in \cref{alg:rsi}.
However, note that the pseudocode has the slightly awkward requirement of performing an \emph{even} number of matrix multiplications ($q$ multiplications with $\bm{A}$ and $q$ multiplications with $\bm{A}^\ast$).
Therefore we ask, ``Is there any way to perform \texttt{RSI} with an \emph{odd} number of multiplications?''
Bjarkason \cite{bjarkason2019pass} answered this question in the affirmative (also see earlier work of \cite[Sec.~4.3]{rokhlin2010randomized}), by replacing the approximation $\hat{\bm{A}} = \bm{\Pi}_{(\bm{A} \bm{A}^{\ast})^{q-1} \bm{A} \bm{\Omega}} \bm{A}$
with the closely related approximation
\begin{equation}
    \label{eq:approximation_new}
\hat{\bm{A}} = \bm{A} \bm{\Pi}_{(\bm{A}^{\ast} \bm{A})^q \bm{\Omega}}.
\end{equation}
The approach requires $q + 1$ matrix-matrix multiplications with $\bm{A}$ and only $q$ matrix--matrix multiplications with $\bm{A}^{\ast}$.
When multiplications with $\bm{A}$ are much cheaper than multiplications with $\bm{A}^{\ast}$ (as in some PDE models \cite{bjarkason2019pass}),
we essentially get the final matrix multiplication for free, and we obtain an improved approximation.

\begin{algorithm}[t]
\caption{\texttt{RSI}, extended version \label{alg:rsi2}}
\begin{algorithmic}[1]
\Require Matrix $\bm{A} \in \mathbb{R}^{L \times N}$; block size $k$; stopping criterion
\Ensure Orthogonal $\bm{U} \in \mathbb{R}^{L \times k}$,
orthogonal $\bm{V} \in \mathbb{R}^{N \times k}$, and diagonal $\bm{\Sigma} \in \mathbb{R}^{k \times k}$ such that $\bm{A} \approx \bm{U} \bm{\Sigma} \bm{V}^{\ast}$
\State Generate a random matrix $\bm{Y} \in \mathbb{R}^{N \times k}$
\State $i = 1$
\Comment{$i$ counts the number of multiplications}
\State $\bm{X} = \bm{A} \bm{Y}$
\State $[\bm{X}, \sim] = \textup{qr\_econ}(\bm{X})$
\While{stopping criterion is not met}
\State $i = i + 1$
\If{$i$ is even}
    \State $\bm{Y} = \bm{A}^{\ast} \bm{X}$
    \State $[\bm{Y}, \bm{R}] = \textup{qr\_econ}(\bm{Y})$
    \Comment{Optional: stabilized QR \cref{eq:stabilized_qr}}
    \State $\bm{T} = \bm{R}^\ast$
    \Comment{$\bm{T}$ is lower triangular}
\Else
    \State $\bm{X} = \bm{A} \bm{Y}$
    \State $[\bm{X}, \bm{S}] = \textup{qr\_econ}(\bm{X})$
    \Comment{Optional: stabilized QR \cref{eq:stabilized_qr}}
    \State $\bm{T} = \bm{S}$
    \Comment{$\bm{T}$ is upper triangular}
\EndIf
\EndWhile
\State $[\hat{\bm{U}}, \bm{\Sigma}, \hat{\bm{V}}] = \textup{svd\_econ}(\bm{T})$
\State $\bm{U} = \bm{X} \hat{\bm{U}}$
\State $\bm{V} = \bm{Y} \hat{\bm{V}}$
\end{algorithmic}
\end{algorithm}

Bjarkason's insight leads to an extended version of \texttt{RSI} that uses either an odd or an even number of matrix multiplications.
We provide an implementation in \cref{alg:rsi2}, which incorporates the following features:
\begin{enumerate}
\item The algorithm generates orthonormal matrices $\bm{X} \in \mathbb{R}^{L \times k}$ 
and $\bm{Y} \in \mathbb{R}^{N \times k}$,
representing the range and co-range of the low-rank approximation.
The low-rank approximation takes the form $\hat{\bm{A}} = \bm{X} \bm{T} \bm{Y}^\ast$ for a core matrix $\bm{T} \in \mathbb{R}^{k \times k}$.
\item At odd iterations, the algorithm updates $\bm{X}$ and $\bm{T}$. At even iterations, the algorithm updates $\bm{Y}$ and $\bm{T}$.
\item The algorithm uses the fact that
\begin{equation*}
    \bm{T} = \hat{\bm{U}} \bm{\Sigma} \hat{\bm{V}}^\ast
    \quad \text{implies} \quad
    \hat{\bm{A}} = (\bm{X} \hat{\bm{U}}) \bm{\Sigma} (\bm{Y} \hat{\bm{V}})^\ast,
\end{equation*}
to efficiently compute an SVD for the low-rank approximation $\hat{\bm{A}}$.
\end{enumerate}
To optimize the number of multiplications, we can run \cref{alg:rsi2} with an adaptive stopping criterion.
For example, we can stop the algorithm as soon as the square Frobenius--norm error $\lVert \bm{A} - \hat{\bm{A}} \rVert_{\rm F}^2 = \lVert \bm{A} \rVert_{\rm F}^2 - \lVert \hat{\bm{A}} \rVert_{\rm F}^2$ falls below a target precision $\varepsilon^2 \cdot \lVert \bm{A} \rVert_{\rm F}^2$.
See \cref{sec:parameter_choices} for a discussion.

\subsection{\texttt{RBKI}} \label{sec:rbki}

In light of \cref{lem:bigger_projection}, it is always better to project onto a bigger subspace, and
\emph{randomized block Krylov iteration} projects onto the biggest subspace generated by the sequential
matrix products.  We define the Krylov subspace by
\begin{align*}
    \mathcal{K}_q(\bm{A} \bm{A}^{\ast}; \bm{A} \bm{\Omega})
    &= \textup{range} \begin{bmatrix}
    \bm{A} \bm{\Omega} & (\bm{A} \bm{A}^{\ast}) \bm{A} \bm{\Omega} & \cdots &  (\bm{A} \bm{A}^{\ast})^{q-1} \bm{A} \bm{\Omega}
    \end{bmatrix} \\
    &= \bigcup\nolimits_{\textup{deg}(\phi) \leq q - 1} \textup{range}[ \phi(\bm{AA}^\ast) \bm{A} \bm{\Omega} ].
\end{align*}
In this expression, $\phi$ varies over polynomials with degree not exceeding $q - 1$.
The \texttt{RBKI} method compresses the input matrix to the Krylov subspace:
\begin{equation*}
\tag{\texttt{RBKI}}
    \hat{\bm{A}} = \bm{\Pi}_{[\begin{smallmatrix}
    \bm{A} \bm{\Omega} & \cdots &  (\bm{A} \bm{A}^{\ast})^{q-1} \bm{A} \bm{\Omega}
    \end{smallmatrix}]} \bm{A}.
\end{equation*}
The Krylov subspace is much larger than the projection space used in \texttt{RSI},
which only contains the range of the final power $(\bm{AA}^\ast)^{q-1}\bm{A\Omega}$.
As a consequence, the \texttt{RBKI} approximation has the potential to be much more accurate.

\texttt{RBKI} was first introduced by Rokhlin, Szlam, and Tygert \cite{rokhlin2010randomized}.
The paper~\cite{halko2011algorithm} contains a numerical study in the context of principal
component analysis.  The first theoretical analysis is due to Musco and Musco \cite{musco2015randomized}.
Over the past decade, \texttt{RBKI} has been applied to problems in genetics \cite{galinsky2016fast}, meteorology \cite{bousserez2020enhanced},
and geophysical imaging \cite{yang2021low},
leading to reported higher accuracy than \texttt{RSI}. 
Nonetheless, \texttt{RBKI} has remained under-utilized in applications (as reflected in the computational genetics review \cite{tsuyuzaki2020benchmarking}).

\begin{algorithm}[t]
\caption{\texttt{RBKI}, simple version \label{alg:rbki}}
\begin{algorithmic}[1]
\Require Matrix $\bm{A} \in \mathbb{R}^{L \times N}$; block size $k$; iteration count $q$
\Ensure Orthogonal $\bm{U} \in \mathbb{R}^{L \times qk}$,
orthogonal $\bm{V} \in \mathbb{R}^{N \times qk}$, and diagonal $\bm{\Sigma} \in \mathbb{R}^{qk \times qk}$ such that $\bm{A} \approx \bm{U} \bm{\Sigma} \bm{V}^{\ast}$
\State Generate a random matrix $\bm{Y}_0 \in \mathbb{R}^{N \times K}$.
\For{$i = 1, \ldots, q$}
\State $\bm{X}_i = \bm{A} \bm{Y}_{i-1}$
\State $\bm{X}_i = \bm{X}_i - \sum_{j < i} \bm{X}_j (\bm{X}_j^{\ast} \bm{X}_i)$
\Comment{Orthogonalize w.r.t. past iterates}
\State $\bm{X}_i = \bm{X}_i - \sum_{j < i} \bm{X}_j (\bm{X}_j^{\ast} \bm{X}_i)$
\Comment{The repetition ensures stability}
\State $[\bm{X}_i, \sim] = \textup{qr\_econ}(\bm{X}_i)$
\Comment{Optional: stabilized QR \cref{eq:stabilized_qr}}
\State $\bm{Y}_i = \bm{A}^{\ast} \bm{X}_i$
\EndFor
\State $[\hat{\bm{U}}, \bm{\Sigma}, \bm{V}] = \textup{svd\_econ}\bigl(
\begin{bmatrix} \bm{Y}_1 & \cdots & \bm{Y}_q \end{bmatrix}^\ast\bigr)$
\State $\bm{U} = \begin{bmatrix} \bm{X}_1 & \cdots & \bm{X}_q \end{bmatrix} \hat{\bm{U}}$
\end{algorithmic}
\end{algorithm}

\begin{algorithm}[t]
\caption{\texttt{RBKI}, extended version \label{alg:rbki2}}
\begin{algorithmic}[1]
\Require Matrix $\bm{A} \in \mathbb{R}^{L \times N}$; block size $k$; stopping criterion
\Ensure Orthogonal $\bm{U}$,
orthogonal $\bm{V}$, and diagonal $\bm{\Sigma}$ such that $\bm{A} \approx \bm{U} \bm{\Sigma} \bm{V}^{\ast}$
\State Generate a random matrix $\bm{Y}_0 \in \mathbb{R}^{N \times k}$
\State $i = 1$
\Comment{$i$ counts the number of multiplications}
\State $\bm{X}_1 = \bm{A} \bm{Y}_0$
\State $[\bm{X}_1, \sim] = \textup{qr\_econ}(\bm{X}_1)$
\State $\bm{R} = [\,]; \bm{S} = [\,]$
\While{stopping criterion is not met}
\State $i = i + 1$
\If{$i$ is even}
    \State $\bm{Y}_i = \bm{A}^\ast \bm{X}_{i-1}$
    \State $\bm{R}_{\bullet i} = \begin{bmatrix} \bm{Y}_2 & \bm{Y}_4 & \cdots & \bm{Y}_{i-2} \end{bmatrix}^\ast \bm{Y}_i$
    \State $\bm{Y}_i = \bm{Y}_i - \sum_{\textup{even } j < i} \bm{Y}_j
    (\bm{Y}_j^\ast \bm{Y}_i)$
    \Comment{Orthog. w.r.t. even iterates (2$\times$)} 
    \State $[\bm{Y}_i, \bm{R}_{ii}] = \textup{qr\_econ}(\bm{Y}_i)$
    \Comment{Optional: stabilized QR \cref{eq:stabilized_qr}}
    \State $\bm{R} = \begin{bmatrix} \bm{R} & \bm{R}_{\bullet i} \\ \bm{0} & \bm{R}_{ii} \end{bmatrix}$
    \State $\bm{T} = \bm{R}^\ast$
\Else
    \State $\bm{X}_i = \bm{A} \bm{Y}_{i-1}$
    \State $\bm{S}_{\bullet i} = \begin{bmatrix} \bm{X}_1 & \bm{X}_3 & \cdots & \bm{X}_{i-2} \end{bmatrix}^\ast \bm{X}_i$
    \State $\bm{X}_i = \bm{X}_i - \sum_{\textup{odd } j < i} \bm{X}_j
    (\bm{X}_j^\ast \bm{X}_i)$
    \Comment{Orthog. w.r.t. odd iterates (2$\times$)} 
    \State $[\bm{X}_i, \bm{S}_{ii}] = \textup{qr\_econ}(\bm{X}_i)$
    \Comment{Optional: stabilized QR \cref{eq:stabilized_qr}}
    \State $\bm{S} = \begin{bmatrix} \bm{S} & \bm{S}_{\bullet i} \\ \bm{0} & \bm{S}_{ii} \end{bmatrix}$
    \State $\bm{T} = \bm{S}$
\EndIf
\EndWhile
\State $[\hat{\bm{U}}, \bm{\Sigma}, \hat{\bm{V}}] = \textup{svd\_econ} (\bm{T})$
\State $\bm{U} = \begin{bmatrix} \bm{X}_1 & \bm{X}_3 & \cdots \end{bmatrix} \hat{\bm{U}}$
\State $\bm{V} = \begin{bmatrix} \bm{Y}_2 & \bm{Y}_4 & \cdots \end{bmatrix} \hat{\bm{V}}$
\end{algorithmic}
\end{algorithm}

\Cref{alg:rbki} presents pseudocode for a simple version of \texttt{RBKI}, while
\Cref{alg:rbki2} displays pseudocode for an extended version of \texttt{RBKI} that uses an even or odd number of multiplications.  Note that the block Gram--Schmidt step (lines 4 and 5 in \cref{alg:rbki}) must be performed twice for numerical stability.
In subsequent algorithm displays, we indicate this repetition with the label (2$\times$) rather than writing out the
formula twice.

Our \texttt{RBKI} pseudocode is more efficient than existing \texttt{RBKI} implementations \cite{rokhlin2010randomized,halko2011algorithm,musco2015randomized,galinsky2016fast,yuan2018superlinear,drineas2018structural,bjarkason2019pass,yang2021low,martinsson2020randomized}.  Indeed, the existing procedures all require the following series of multiplications:
\begin{itemize}
    \item[(a)] $q$ times: multiply
    $\bm{A}$ with a matrix of size $N \times k$;
    \item[(b)] $q-1$ times: multiply $\bm{A}^{\ast}$ with a matrix of size $L \times k$;
    and 
    \item[(c)] $1$ time: multiply $\bm{A}^{\ast}$ with a block matrix of size $L \times kq$.
\end{itemize}
In contrast, \cref{alg:rbki} uses a reduced set of multiplications:
\begin{itemize}
\item[(a)] $q$ times: multiply
$\bm{A}$ with a matrix of size $N \times k$; and 
\item[(b)] $q$ times:
multiply $\bm{A}^{\ast}$ with a matrix of size $L \times k$.
\end{itemize}
Our new procedure removes step (c), which involves an expensive multiplication of $\bm{A}^{\ast}$
with a block matrix of size $L \times kq$.  We have substituted a single multiplication
with a smaller matrix of size $L \times k$.

If we measure computational cost in the simplest way, by counting matrix--vector products, our new pseudocode results in cost savings of roughly $33\%$.
The savings are even higher when matrix--vector multiplications are more expensive with $\bm{A}^\ast$ than with $\bm{A}$ (as in some PDE models \cite{bjarkason2019pass}).
Our new pseudocode does not change the asymptotic arithmetic cost of \texttt{RBKI}, which is still $\mathcal{O}(kmLN)$ for dense matrices, yet it makes a noticeable difference in applications.  There are also potential reductions in communication costs.

Our new pseudocode enables a clean comparison between \texttt{RSI} and \texttt{RBKI}
because the two algorithms now perform matrix multiplications
with precisely the same block size.
Although \texttt{RBKI} performs more
arithmetic than \texttt{RSI} outside the matrix products,
this arithmetic is rarely the predominant cost of the algorithm. 
We acknowledge that \texttt{RBKI} requires  $\mathcal{O}(km(L + N))$ storage
whereas \texttt{RSI} only requires $\mathcal{O}(k(L + N))$ storage.
However, in typical applications involving data matrices (see \cref{sec:experiments}), this is not an issue
because storing the low-rank approximation $\hat{\bm{A}}$ is much cheaper
than storing the original matrix $\bm{A}$.

As the main distinction,
\texttt{RBKI} makes better use of the matrix products than \texttt{RSI},
resulting in a far more accurate approximation.
Overall, we suspect that \texttt{RSI} is preferable to \texttt{RBKI}
only if $\bm{A}$ is an abstract linear operator and working storage is the
limiting factor for the computation.

\subsection{Positive-semidefinite matrices} \label{sec:randomized_nystrom}

Suppose that $\bm{A} \in \mathbb{R}^{N \times N}$ is psd and we have generated a matrix $\bm{M} \in \mathbb{R}^{N \times k}$
whose range is aligned with leading eigenvectors of $\bm{A}$.
Then, the simplest approximation for $\bm{A}$, which is used in \texttt{RSVD}, \texttt{RSI}, and \texttt{RBKI}, is based on compressing the range or co-range of $\bm{A}$ as follows:
\begin{equation*}
    \bm{\Pi}_{\bm{M}} \bm{A} = \bm{M} (\bm{M}^{\ast} \bm{M})^{\dagger} (\bm{AM})^{\ast}
    \qquad \text{or} \qquad
    \bm{A} \bm{\Pi}_{\bm{M}} = (\bm{A} \bm{M}) (\bm{M}^{\ast} \bm{M})^{\dagger} \bm{M}^{\ast}.
\end{equation*}
Since $\bm{A}$ is psd, we can employ the same data more effectively using
the \emph{Nystr\"om approximation} \cite[Sec.~14]{martinsson2020randomized}:
\begin{equation*}
    \bm{A}\langle\bm{M}\rangle := \bm{A}^{1 \slash 2} \bm{\Pi}_{\bm{A}^{1 \slash 2} \bm{M}} \bm{A}^{1 \slash 2}
    = (\bm{A} \bm{M}) (\bm{M}^{\ast} (\bm{A} \bm{M}))^{\dagger} (\bm{AM})^{\ast}.
\end{equation*}
The Nystr{\"o}m approximation always improves over the simpler approximations.

\begin{lemma}[Nystr{\"o}m helps] \label{lem:nystrom}
Consider a psd matrix $\bm{A} \in \mathbb{R}^{N \times N}$ and any matrix $\bm{M} \in \mathbb{R}^{N \times k}$.
Then
\begin{equation}
\label{eq:nystrom_error}
    \lVert \bm{A} - \bm{A}\langle\bm{M}\rangle \rVert_p
    \leq \lVert \bm{A} - \bm{\Pi}_{\bm{M}} \bm{A} \rVert_p = \lVert \bm{A} - \bm{A} \bm{\Pi}_{\bm{M}} \rVert_p
\end{equation}
for any Schatten $p$-norm with $1 \leq p \leq \infty$.  The same inequality holds for any unitarily invariant norm.
\end{lemma}
\begin{proof}
Since $\bm{\Pi}_{\bm{M}} \bm{A}$ and $\bm{A}\langle\bm{M}\rangle$ only depend on the range of $\bm{M}$, we can assume without loss of generality that $\bm{M}$ has orthonormal columns and change basis so that
$\bm{M} = \bigl[\begin{smallmatrix} \mathbf{I} \\ \bm{0} \end{smallmatrix}\bigr]$.
With this transformation, the matrices $\bm{A}$, $\bm{\Pi}_{\bm{M}} \bm{A}$, and $\bm{A}\langle\bm{M}\rangle$ can be written in block form:
\begin{equation*}
    \bm{A} = \begin{bmatrix}
    \bm{A}_{11} & \bm{A}_{12} \\
    \bm{A}_{21} & \bm{A}_{22}
    \end{bmatrix},
    \quad
    \bm{\Pi}_{\bm{M}} \bm{A} = \begin{bmatrix}
    \bm{A}_{11} & \bm{A}_{12} \\
    \bm{0} & \bm{0}
    \end{bmatrix},
    \quad \bm{A}\langle\bm{M}\rangle = \begin{bmatrix}
    \bm{A}_{11} & \bm{A}_{12} \\
    \bm{A}_{21} & \bm{A}_{21} \bm{A}_{11}^{\dagger} \bm{A}_{12}
    \end{bmatrix}.
\end{equation*}
Using the psd order relations $\bm{A}_{22} - \bm{A}_{21} \bm{A}_{11}^{\dagger} \bm{A}_{12} \preccurlyeq \bm{A}_{22}$ and $\bm{A}_{22}^2 \preccurlyeq \bm{A}_{22}^2 + \bm{A}_{21} \bm{A}_{12}$, we calculate
\begin{align*}
    \sigma_i(\bm{A} - \bm{A}\langle\bm{M}\rangle)^2
    &= \lambda_i\biggl(
    \begin{bmatrix}
    \bm{0} & \bm{0} \\
    \bm{0} & \bm{A}_{22} - \bm{A}_{21} \bm{A}_{11}^{\dagger} \bm{A}_{12}
    \end{bmatrix}
    \biggr)^2 \\
    &\leq \lambda_i\biggl(
    \begin{bmatrix}
    \bm{0} & \bm{0} \\
    \bm{0} & \bm{A}_{22}
    \end{bmatrix}
    \biggr)^2
    = \lambda_i\biggl(
    \begin{bmatrix}
    \bm{0} & \bm{0} \\
    \bm{0} & \bm{A}_{22}^2
    \end{bmatrix}
    \biggr) \\
    &\leq \lambda_i\biggl(
    \begin{bmatrix}
    \bm{0} & \bm{0} \\
    \bm{0} & \bm{A}_{22}^2 + \bm{A}_{21} \bm{A}_{12}
    \end{bmatrix}
    \biggr)
    = \lambda_i((\bm{A} - \bm{\Pi}_{\bm{M}}\bm{A}) (\bm{A} - \bm{\Pi}_{\bm{M}}\bm{A})^\ast) \\
    &= \sigma_i(\bm{A} - \bm{\Pi}_{\bm{M}}\bm{A})^2
\end{align*}
for each $1 \leq i \leq N$, which confirms \cref{eq:nystrom_error}.
\end{proof}

Next, we show how to incorporate Nystr\"om approximation in three different randomized low-rank approximation algorithms: \texttt{NysSVD} (\cref{sec:rnystrom}), \texttt{NysSI} (\cref{sec:nysrsi}), and \texttt{NysBKI} (\cref{sec:nysrbki}).

\subsubsection{\texttt{NysSVD}} \label{sec:rnystrom}

In \emph{randomized Nystr\"om approximation}, we generate a low-rank approximation:
\begin{equation*}
\tag{\texttt{NysSVD}}
\hat{\bm{A}} = \bm{A}\langle \bm{\Omega} \rangle,
\end{equation*}
where $\bm{\Omega} \in \mathbb{R}^{N \times k}$ is a random matrix.
In this work, $\bm{\Omega}$ is always a Gaussian matrix \cite{tropp2017fixed}, mimicking the strategy used in \texttt{RSVD}.  For parallelism, we call the resulting algorithm \texttt{NysSVD} even though it produces an eigenvalue decomposition.

\begin{algorithm}[t]
\caption{\texttt{NysSVD} \label{alg:nyssvd}}
\begin{algorithmic}[1]
\Require Psd matrix $\bm{A} \in \mathbb{R}^{N \times N}$; block size $k$; shift $\varepsilon > 0$
\Ensure Orthogonal $\bm{U} \in \mathbb{R}^{N \times k}$ and diagonal $\bm{\Lambda} \in \mathbb{R}^{k \times k}$ such that $\bm{A} \approx \bm{U} \bm{\Lambda} \bm{U}^{\ast}$
\State Generate a random matrix $\bm{\Omega} \in \mathbb{R}^{N \times k}$
\State $\bm{Y} = \bm{A} \bm{\Omega} + \varepsilon \bm{\Omega}$ \Comment{Apply $\varepsilon$ shift}
\State $\bm{C} = \textup{chol}(\bm{\Omega}^{\ast} \bm{Y})$
\State $\bm{Z} = \bm{Y} \bm{C}^{-1}$ \Comment{Triangular solve}
\State $[\bm{U}, \bm{\Sigma}, \sim] = \textup{svd\_econ}(\bm{Z})$
\State $\bm{\Lambda} = \max\left\{\bm{0}, \bm{\Sigma}^2 - \varepsilon \mathbf{I}\right\}$ \Comment{Remove $\varepsilon$ shift}
\end{algorithmic}
\end{algorithm}

\Cref{alg:nyssvd} presents \texttt{NysSVD} pseudocode, which is optimized for efficiency in the following ways:
\vspace{0.2cm}
\begin{enumerate}
\item The low-rank approximation
\begin{equation*}
    \bm{A}\langle \bm{\Omega} \rangle = \bm{A} \bm{\Omega} (\bm{\Omega}^\ast \bm{A} \bm{\Omega})^{\dagger} \bm{\Omega}^\ast \bm{A}^\ast,
\end{equation*}
is generated without ever forming the pseudoinverse $(\bm{\Omega}^\ast \bm{A} \bm{\Omega})^{\dagger}$ explicitly.
Instead, we compute an upper-triangular Cholesky factor $\bm{C} \in \mathbb{R}^{k \times k}$ so that $\bm{\Omega}^\ast (\bm{A} \bm{\Omega}) = \bm{C}^\ast \bm{C}$.  The low-rank approximation is generated as
\begin{equation*}
    \bm{A}\langle \bm{\Omega} \rangle = (\bm{A} \bm{\Omega} \bm{C}^{\dagger}) (\bm{A} \bm{\Omega} \bm{C}^{\dagger})^\ast.
\end{equation*}
\item To ensure numerical stability, the Nystr\"om approximation is applied to a shifted operator $\bm{A}_{\varepsilon} = \bm{A} + \varepsilon \mathbf{I}$, where the shift parameter $\varepsilon = \varepsilon_{\textup{mach}} \textup{tr}(\bm{A})$ depends on the machine precision.
To approximately counteract the shift, the eigenvalues of the approximation $\hat{\bm{A}}_{\varepsilon}$ are all reduced by $\varepsilon$.
\end{enumerate}
\vspace{0.2cm}
These improvements, previously recommended in \cite{li2017algorithm}, lead to the fast and robust \texttt{NysSVD} implementation in \cref{alg:nyssvd}.

\begin{remark}[Column Nystr{\"o}m approximation]
When evaluating each entry of the input matrix $\bm{A}$ is expensive,
we might prefer to construct a Nystr{\"o}m approximation
with the range determined from just a few columns of $\bm{A}$.  Equivalently,
the test matrix $\bm{\Omega}$ contains (random) coordinate vectors~\cite{williams2001using}.
To find an informative set of columns,
we need additional ideas; see the paper~\cite{chen2022randomly} and its background references.
\end{remark}

\subsubsection{\texttt{NysSI}} \label{sec:nysrsi}

In \emph{randomized subspace iteration with Nystr\"om approximation} \cite[Alg.~5.5]{halko2011finding}, we approximate a psd matrix $\bm{A} \in \mathbb{R}^{N \times N}$ as:
\begin{equation*}
\tag{\texttt{NysSI}}
    \hat{\bm{A}} = \bm{A}\langle \bm{M} \rangle, \quad \text{where} \quad \bm{M} = \bm{A}^{m-1} \bm{\Omega}.
\end{equation*}
\Cref{alg:nyssi} presents \texttt{NysSI} pseudocode, which can be implemented with any adaptive stopping criterion.
For example, we can stop the algorithm as soon as the trace--norm error 
$\lVert \bm{A} - \hat{\bm{A}} \rVert_{\ast} = \lVert \bm{A} \rVert_{\ast} - \lVert \hat{\bm{A}} \rVert_{\ast}$ falls below a target precision $\varepsilon \cdot \lVert \bm{A} \rVert_{\ast}$.
\texttt{NysSI} achieves a target precision more quickly than \texttt{RSI} by one-half of a matrix--matrix multiplication, as reflected in our experiments (\cref{fig:krylov_v_subspace,fig:eigenvectors}) and theory (\cref{sec:theory2}).

\begin{algorithm}[t]
\caption{\texttt{NysSI} \label{alg:nyssi}}
\begin{algorithmic}[1]
\Require Psd matrix $\bm{A} \in \mathbb{R}^{N \times N}$; block size $k$; stopping criterion; shift $\varepsilon > 0$
\Ensure Orthogonal $\bm{U} \in \mathbb{R}^{N \times k}$ and diagonal $\bm{\Lambda} \in \mathbb{R}^{k \times k}$ such that $\bm{A} \approx \bm{U} \bm{\Lambda} \bm{U}^{\ast}$
\State Generate a random matrix $\bm{Y} \in \mathbb{R}^{N \times k}$
\State $i = 0$
\While{stopping criterion is not met}
\State $i = i + 1$
\State $[\bm{X}, \sim] = \textup{qr\_econ}(\bm{Y})$
\Comment{Optional: stabilized QR \cref{eq:stabilized_qr}}
\State $\bm{Y} = \bm{A} \bm{X}$
\EndWhile
\State $\bm{Y} = \bm{Y} + \varepsilon \bm{X}$ \Comment{Apply $\varepsilon$ shift}
\State $\bm{C} = \textup{chol}(\bm{X}^{\ast} \bm{Y})$ \Comment{Cholesky decomposition}
\State $\bm{Z} = \bm{Y} \bm{C}^{-1}$ \Comment{Triangular solve}
\State $[\bm{U}, \bm{\Sigma}, \sim] = \textup{svd\_econ}(\bm{Z})$
\State $\bm{\Lambda} = \max\left\{\bm{0}, \bm{\Sigma}^2 - \varepsilon \mathbf{I}\right\}$ \Comment{Remove $\varepsilon$ shift}
\end{algorithmic}
\end{algorithm}

\subsubsection{\texttt{NysBKI}} \label{sec:nysrbki}

In \emph{randomized block Krylov iteration with Nystr\"om approximation}, we approximate a psd matrix $\bm{A} \in \mathbb{R}^{N \times N}$ as
\begin{equation*}
\tag{\texttt{NysBKI}}
    \hat{\bm{A}} = \bm{A}\langle \bm{M} \rangle, \quad \text{where} \quad \bm{M} = \begin{bmatrix}
    \bm{\Omega} & \bm{A} \bm{\Omega} & \bm{A}^2 \bm{\Omega} & \cdots &  \bm{A}^{m-1} \bm{\Omega}
    \end{bmatrix}.
\end{equation*}
\texttt{NysBKI} uses the output of \emph{every} matrix multiplication, not just every second matrix multiplication, to build an enriched approximation space.
As a result, \texttt{NysBKI}
is more efficient than \texttt{RBKI} by a factor of $\sqrt{2}$ matrix--matrix multiplications (\cref{fig:krylov_v_subspace,fig:eigenvectors,sec:theory2}).
To the best of our knowledge, \texttt{NysBKI} has not appeared in the matrix approximation literature before now.
\texttt{NysBKI} pseudocode appears in \cref{alg:nysbki}.

\begin{algorithm}[t]
\caption{\texttt{NysBKI} \label{alg:nysbki}}
\begin{algorithmic}[1]
\Require Psd matrix $\bm{A} \in \mathbb{R}^{N \times N}$; block size $k$; stopping criterion; shift $\varepsilon > 0$
\Ensure Orthogonal $\bm{U}$ and diagonal $\bm{\Lambda}$ such that $\bm{A} \approx \bm{U} \bm{\Lambda} \bm{U}^{\ast}$
\State Generate a random matrix $\bm{Y}_0 \in \mathbb{R}^{N \times k}$
\State $i = 0$
\While{stopping criterion is not met}
\State $\bm{X}_i = \bm{Y}_i$
\State $\bm{X}_i = \bm{X}_i - \sum_{j < i} \bm{X}_j (\bm{X}_j^{\ast} \bm{X}_i)$
\Comment{Orthog.~w.r.t.~past iterates (2$\times$)}
\State $[\bm{X}_i, \sim] = \textup{qr\_econ}(\bm{X}_i)$
\Comment{Optional: stabilized QR \cref{eq:stabilized_qr}}
\State $\bm{Y}_{i+1} = \bm{A} \bm{X}_i + \varepsilon \bm{X}_i$
\State $i = i + 1$
\EndWhile
\State $\bm{C} = \textup{chol}\bigl(\begin{bmatrix} \bm{X}_0 & \cdots & \bm{X}_{i-1} \end{bmatrix}^{\ast} \begin{bmatrix} \bm{Y}_1 & \cdots & \bm{Y}_i \end{bmatrix}\bigr)$
\State $\bm{Z} = \begin{bmatrix} \bm{Y}_1 & \cdots & \bm{Y}_i \end{bmatrix} \bm{C}^{-1}$
\State $[\bm{U}, \bm{\Sigma}, \sim] = \textup{svd\_econ}(\bm{Z})$
\State $\bm{\Lambda} = \max\{\bm{0}, \bm{\Sigma}^2 - \varepsilon \mathbf{I}\}$
\end{algorithmic}
\end{algorithm}

\subsection{Additional opportunities for speedups} \label{sec:further}

So far, we have presented simple, stable pseudocode that requires the minimal number of matrix multiplications with $\bm{A}$ and $\bm{A}^{\ast}$.
However, we have not necessarily optimized the operations involving smaller matrices, and the algorithms use a large number of QR and SVD factorizations.  Here we discuss some possibilities for computational speedups by removing or replacing these factorizations.  This section is intended for numerical linear algebra experts, and most readers can skip it with impunity.

As the first speedup opportunity, we can remove all the SVD steps in our pseudocode \cite{cheng2016fast}.
We can redesign \texttt{RSVD} to return a factorization $\hat{\bm{A}} = \bm{X} \bm{Y}^\ast$, without ever calculating the singular values and vectors.
We can make similar adjustments to the pseudocode for other randomized low-rank approximation algorithms.
After making these changes, the algorithms return the same low-rank approximation (in exact precision arithmetic), but the singular values and singular vectors of $\hat{\bm{A}}$ are no longer easily accessible as before.

As the second speedup opportunity, we can replace all the QR factorizations in the Nystr\"om-based algorithms
and all but the last QR factorization in \texttt{RSI} with cheaper matrix decompositions.
Indeed, the sole purpose of the QR factorizations is to
prevent the iterates from becoming ill-conditioned, but we can also avoid ill-conditioning
using an LU factorization \cite{li2017algorithm}
or a randomized QR decomposition \cite{balabanov2020randomized,balabanov2021randomized}.
These alternatives use roughly $50\%$ of the arithmetic of the standard pivoted QR factorization.

For \texttt{RBKI} and \texttt{NysBKI}, we can also replace the block orthogonalization steps (such
as lines 4--5 in \cref{alg:rbki}) by a shorter Lanczos-type recurrence \cite{nakatsukasa2021fast}:
\begin{equation*}
\bm{X}_i = \bm{X}_i - \sum\nolimits_{j = i-2}^{i-1} \bm{X}_j (\bm{X}_j^\ast \bm{X}_i).
\end{equation*}
This reduces the total orthogonalization cost from $\mathcal{O}(k^2 q^2(L + N))$ to $\mathcal{O}(k^2 q(L + N))$ arithmetic operations.
On the other hand, in finite precision arithmetic it may result in a basis $\bm{M} = \begin{bmatrix} \bm{X}_1 & \cdots & \bm{X}_q \end{bmatrix}$ that is not orthogonal.  

Last,
orthogonalization is used in \texttt{RSVD}, \texttt{RSI}, and \texttt{RBKI} to compute a low-rank approximation $\hat{\bm{A}} = \bm{\Pi}_{\bm{M}} \bm{A}$ from a nonorthogonal basis matrix $\bm{M}$.
However, if we are willing to accept some additional error, Nakatsukasa \cite[Eq.~3]{nakatsukasa2020fast} has proposed replacing $\bm{\Pi}_{\bm{M}} \bm{A} = \bm{M} \bm{M}^{\dagger} \bm{A}$
with a cheaper approximation
\begin{equation}
\label{eq:original}
    \hat{\bm{A}} = \bm{M} (\bm{S} \bm{M})^{\dagger} \bm{S} \bm{A},
\end{equation}
where $\bm{S}$ is a randomized embedding matrix that reduces the dimensionality of $\bm{M}$ and $\bm{A}$.
Nakatsukasa \cite{nakatsukasa2020fast} finds the additional error to be small in numerical tests.

\section{Parameter choices} \label{sec:parameter_choices}

In this section, we discuss strategies for choosing the block size $k$ and the number of matrix--matrix multiplications $m$ in low-rank approximation algorithms.
In \texttt{RSVD} and \texttt{NysSVD}, there is no choice but to increase $k$ to handle more difficult problems.
However, in the other algorithms, a major question concerns the relative advantages of using a large depth ($m \gg 1$) versus a large block size ($k \gg 1$).
Section \ref{sec:block_size} considers the tradeoff between $k$ and $m$,
while \cref{sec:truncate} discusses adaptive strategies for choosing $k$ and $m$ that ensure a high-quality approximation.

\subsection{Tradeoff between the block size and depth} \label{sec:block_size}

A recent analysis \cite{meyer2023unreasonable} of Meyer, Musco, and Musco gives evidence that $k = 1$ often leads to the highest accuracy in \texttt{RBKI}, assuming a fixed number of matrix--vector products $km$.
Setting $k = 1$ leads to the classical Lanczos iteration.
However, Meyer and coauthors also acknowledge two problems with setting $k = 1$.
First, the approximation quality can be poor if any singular values have multiplicity greater than one (or have exponentially small singular value gaps),
limiting the applicability to eigenvalue problems in physics and chemistry with high multiplicities (e.g., \cite[Sec.~14]{gilmore2008lie}).
Second, the Lanczos method is intrinsically serial, making it slow to run on modern computers.

\begin{figure}[t]
    \centering
    \includegraphics[width=\textwidth]{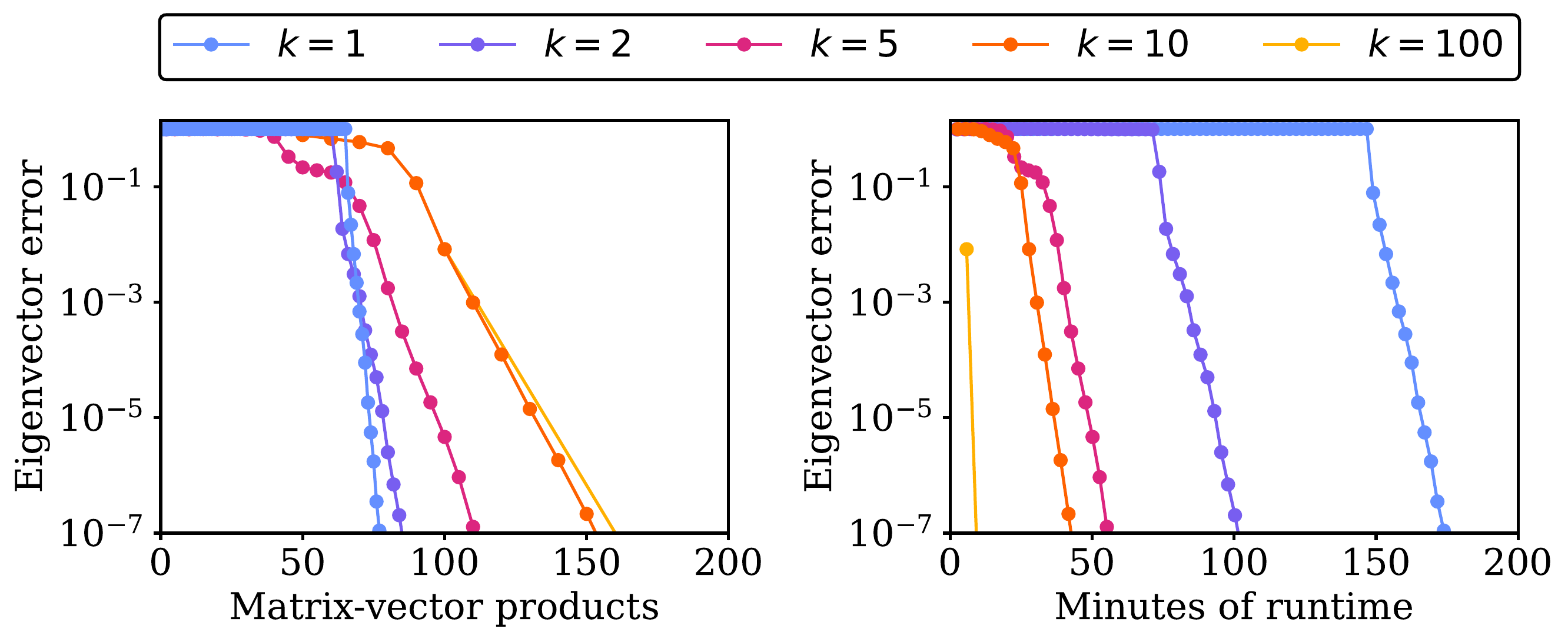}
    \caption{\textbf{Runtime comparisons.} Eigenvector approximation error of \texttt{NysBKI} with block size $k = 1$, $2$, $5$, $10$, or $100$ applied to the $250{,}000 \times 250{,}000$ kernel matrix, as constructed in \cref{sec:clustering}.
    The left panel measures the number of matrix--vector products while the right panel measures the runtime. 
    \label{fig:runtime}}
\end{figure}

We have performed a runtime analysis for \texttt{NysBKI} applied to the $250{,}000 \times 250{,}000$ kernel matrix that will be introduced later in \cref{sec:clustering}, leading to the results pictured in \cref{fig:runtime}.
The figure presents the eigenvector approximation error
\begin{equation*}
    \lVert \bm{\Pi}_3(\bm{A}) - \bm{\Pi}_3(\hat{\bm{A}}) \rVert,
\end{equation*}
where $\bm{\Pi}_3(\cdot)$ denotes the orthogonal projection onto the leadings three eigenvectors, and $\hat{\bm{A}}$ is the stochastic approximation of the kernel matrix $\bm{A}$ from a single run of \texttt{NysBKI}.
We vary the block size to be either $k = 1$, $2$, $5$, $10$, or $100$ (different color lines).
The results in the left panel of \cref{fig:runtime} support the conclusion of Meyer et al. \cite{meyer2023unreasonable} that
the smallest block size $k = 1$ leads to the highest accuracy, given a fixed number of matrix--vector products.
However, the number of matrix--vector multiplications is not a good indicator of the runtime cost, as shown in the right panel of \cref{fig:runtime}.
To reach a tolerance of $\varepsilon = 10^{-5}$, \texttt{NysBKI} with a block size $k = 1$ takes $167$ minutes on a laptop computer, whereas \texttt{NysBKI} with a block size $k = 100$ takes $11$ minutes, \textbf{making it $15 \times$ faster}.
We often anticipate an order-of-magnitude speedup by switching to a much larger block size.

In earlier work, Li and coauthors \cite{li2017algorithm} performed similar speed tests comparing a Matlab implementation of \texttt{RSI} against Lanczos implementations in ARPACK \cite{lehoucq1998arpack} and PROPACK \cite{larsen1998lanczos}.
Their results reinforce the practical advantages of using a large block size, both for reducing runtimes and for reducing errors:
\vspace{.2cm}
\begin{quote}
    On strictly serial processors with no complicated caching (such as the processors of
    many decades ago), the most careful implementations of Lanczos iterations... could likely attain performance nearing the randomized methods’... The randomized methods can
    attain much higher performance on parallel and distributed processors and generally
    are easier to use—setting their parameters properly is trivial (defaults perform well)... 
\end{quote}
\vspace{.2cm}
The conclusions of Li et al.~are based on hundreds of tests with dense matrices as large as $100{,}000 \times 100{,}000$ and sparse matrices as large as $3{,}000{,}000 \times 3{,}000{,}000$.

In light of the computational evidence, we advocate using a block size of at least $k = 10$ and potentially $k = 100$--$1{,}000$ for large-scale problems ($N \geq 10^5$).
After choosing the block size, we recommend running \texttt{RBKI} or \texttt{NysBKI} with an adaptive stopping rule to determine the minimal number of multiplications.
See \cref{sec:truncate} for a discussion of stopping rules that ensure the quality of the low-rank approximation.

\subsection{Quality assurance} \label{sec:truncate}

Here we evaluate two strategies for controlling the quality of a randomized low-rank approximation $\hat{\bm{A}}$.
The first strategy is based on measuring and controlling the global approximation error.
In this strategy, we increase the block size $k$ or depth parameter $m$ until achieving an error tolerance $\lVert \bm{A} - \hat{\bm{A}} \rVert_{\rm F} < \varepsilon \cdot \lVert \bm{A}  \rVert_{\rm F}$
or $\lVert \bm{A} - \hat{\bm{A}} \rVert_{\ast} < \varepsilon \cdot \lVert \bm{A} \rVert_{\ast}$.
These error bounds imply that $\hat{\bm{A}}$ can be reliably used in place of $\bm{A}$ in various matrix computations (\cref{sec:consequences}).

An \texttt{RSI} implementation that controls the global approximation error is provided in the ``svdsketch'' \cite{matlab2020version} function for Matlab, based on pseudocode from \cite{yu2018efficient}:
\begin{itemize}
    \item As the first step, this function evaluates the square Frobenius norm $\lVert \bm{A} \rVert_{\rm F}^2$ using a single pass through the entries of $\bm{A}$.
    \item At each subsequent step, the function adds new columns to the initialization matrix $\bm{\Omega}$, updates the low-rank approximation $\hat{\bm{A}}$, and updates the Frobenius norm $\lVert \hat{\bm{A}} \rVert_{\rm F}$.
    \item The procedure terminates as soon as
    \begin{equation*}
        \lVert \bm{A} - \hat{\bm{A}} \rVert_{\rm F}^2 = \lVert \bm{A} \rVert_{\rm F}^2 - \lVert \hat{\bm{A}} \rVert_{\rm F}^2 < \varepsilon^2 \cdot \lVert \bm{A} \rVert_{\rm F}^2.
    \end{equation*}
\end{itemize}
``Svdsketch'' always returns a low-rank approximation which accounts for a $1 - \varepsilon^2$ proportion of the square Frobenius norm of the target matrix.

Two improvements to ``svdsketch'' could make the procedure even more powerful.
First, ``svdsketch'' adaptively chooses the block size $k$, but our numerical experiments (\cref{sec:experiments}) suggest that adaptively choosing the number of multiplications $m$ would lead to even greater efficiency.
When computations are serial, adding to the depth more quickly reduces the errors than adding to the block size (\cref{fig:runtime} left panel).
Second, ``svdsketch'' can only be applied to matrices stored entry-wise on a computer, because of the way it calculates the Frobenius-norm error.
The recent paper \cite[Sec.~2]{epperly2023efficient} describes a different strategy for approximating the Frobenius norm error based on matrix multiplications, which remains valid even when $\bm{A}$ is an abstract linear operator.

Next we consider a different adaptive strategy, in which we increase $k$ or $m$ until each one of the leading singular vector triplets $(\hat{\bm{u}}_i, \hat{\sigma}_i, \hat{\bm{v}}_i)$ achieves the residual accuracy
\begin{equation}
\label{eq:residual}
    r_i = \Bigl(\lVert (\bm{A} - \hat{\bm{A}})^\ast \hat{\bm{u}}_i \rVert^2 
    + \lVert (\bm{A} - \hat{\bm{A}}) \hat{\bm{v}}_i \rVert^2\Bigr)^{1/2} < \varepsilon.
\end{equation}
As soon as $r_i < \varepsilon$, the $i$th singular vector triplet is stable in the sense of backward error \cite{sun93}.
Namely, there is a perturbation matrix $\bm{E} \in \mathbb{R}^{L \times N}$ with $\lVert \bm{E} \rVert_{\rm F} < \varepsilon$ such that $\bm{A} + \bm{E}$ has exactly the singular vector triplet $(\hat{\bm{u}}_i, \hat{\sigma}_i, \hat{\bm{v}}_i)$.
With assistance from Maksim Melnichenko and Riley Murray,
we have developed adaptive versions of \texttt{RBKI} and \texttt{NysBKI} implementing this strategy, with the pseudocode appearing in \cref{alg:rbki_quality,alg:nysbki_quality}.

\begin{figure}[t]
    \centering
    \includegraphics[width = .5\textwidth]{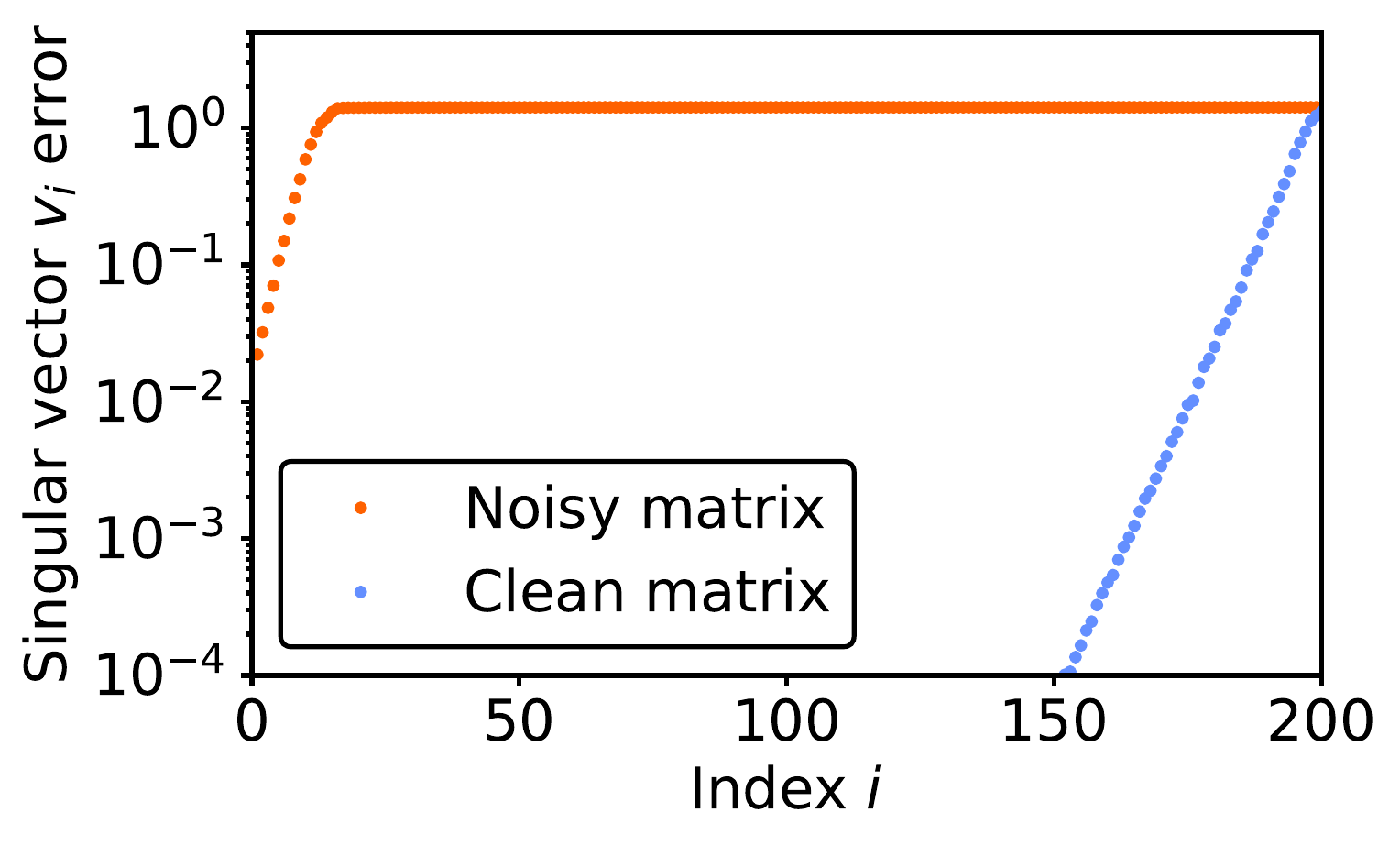}
    \caption{\textbf{(Accuracy of the singular vectors).} Singular vector approximation error for the clean matrix $\bm{A}$ (blue) and the noisy matrix $\bm{B}$ (orange), as constructed in \cref{sec:illustrative}. \label{fig:decaying_accuracy}}
\end{figure}

Both adaptive strategies presented here would lead to a more reliable approximation than the earlier approach
\cite{halko2011finding,martinsson2020randomized} of applying an $r$-truncated singular value decomposition to the low-rank approximation with a default value $r =  k-2$ or $r = \lfloor k / 2 \rfloor$.
A sufficiently small truncation parameter $r$ can in principle eliminate inaccurate singular vector estimates, but selecting $r$ is notoriously difficult in practice.
For example, \cref{fig:decaying_accuracy} evaluates the 
error
\begin{equation*}
    \bigl(\mathbb{E}
    \lVert \bm{\Pi}_{\bm{v}_i(\bm{A})} - \bm{\Pi}_{\bm{v}_i(\hat{\bm{A}})} \rVert^2 \bigr)^{1/2}
    \quad \text{or} \quad
    \bigl(\mathbb{E}\lVert \bm{\Pi}_{\bm{v}_i(\bm{B})} - \bm{\Pi}_{\bm{v}_i(\hat{\bm{B}})} \rVert^2\bigr)^{1/2},
\end{equation*}
for the estimated singular vectors of the clean matrix $\bm{A}$ and noisy matrix $\bm{B}$ described in \cref{sec:illustrative}.
Here, $\bm{\Pi}_{\bm{v}_i(\cdot)}$ denotes the orthogonal projection onto the $i$th right singular vector, and the expectation is evaluated over 100 runs of \texttt{RBKI} with $k = 100$ and $m = 5$.
From the results in \cref{fig:decaying_accuracy}, we would need to apply truncation with $r \leq 4$ for the noisy matrix but we could take $r$ as high as $186$ for the clean matrix, in order to achieve an error tolerance $\varepsilon = .1$.
Unfortunately, there is no default truncation parameter $r$ which performs well for both matrices.

\section{Applications} \label{sec:experiments}

In this section, we present experiments comparing different randomized low-rank approximation algorithms (\cref{sec:comparison}), and we apply the algorithms to scientific problems that require principal component analysis (\cref{sec:genetics}) or kernel spectral clustering (\cref{sec:clustering}).

\subsection{Comparison} \label{sec:comparison}

\begin{figure}[t]
    \centering
    \includegraphics[width = .5\textwidth]{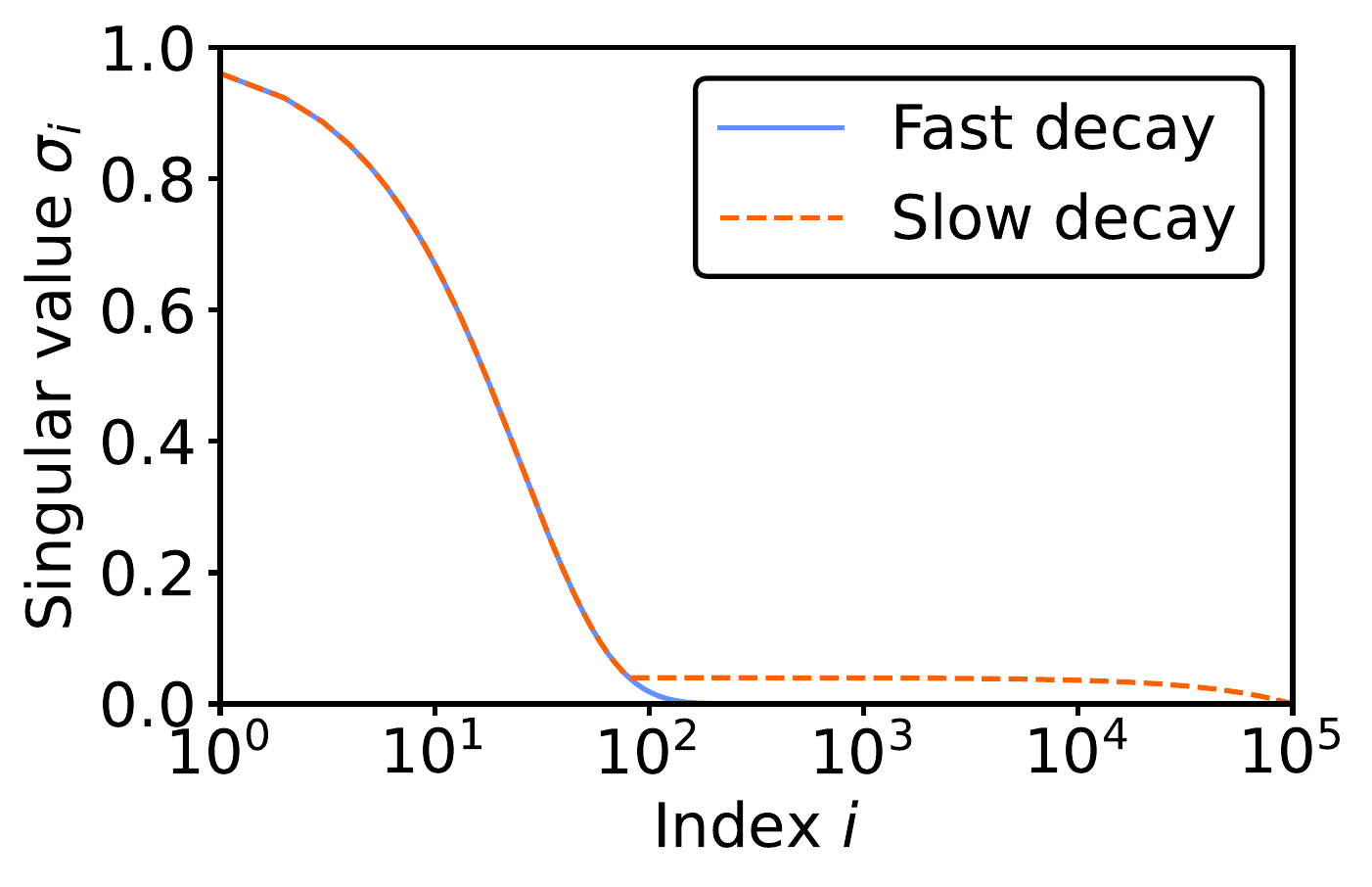}
    \caption{\textbf{(Fast vs. slow singular value decay).} Fast singular value decay of the matrix $\bm{A}$ versus slow singular value decay of the matrix $\bm{B}$, as constructed in \cref{sec:comparison}.}
    \label{fig:synthetic_data}
\end{figure}

As a simple, direct comparison, we apply \texttt{RSVD}, \texttt{RSI}, \texttt{RBKI}, and their Nystr\"om-based variants to approximate two psd matrices $\bm{A}$ and $\bm{B}$ with singular values
\begin{equation*}
    \sigma_i(\bm{A}) = \mathrm{e}^{-i/25}, \qquad 
    \sigma_i(\bm{B}) = \max\biggl\{\mathrm{e}^{-i/25}, \frac{1 - i/10^5}{25} \biggr\}, \qquad
    i = 1, 2, \ldots, 10^5.
\end{equation*}
Matrix $\bm{B}$ models a noisy version of $\bm{A}$, with noise affecting the singular values $\sigma_i(\bm{B})$ for $i \geq 81$ (see \cref{fig:synthetic_data}).
The noise singular values are small, ranging from $0.00$ to $0.04$ in magnitude, but given the high-dimensionality ($L = N = 10^5$), they threaten to drown out the signal.
We construct $\bm{A}$ and $\bm{B}$ as diagonal matrices,
as the choice of singular vectors does not affect the performance of our algorithms in exact arithmetic (\cref{lem:diagonal}).

\begin{figure}[t]
    \centering
    \includegraphics[width = \textwidth]{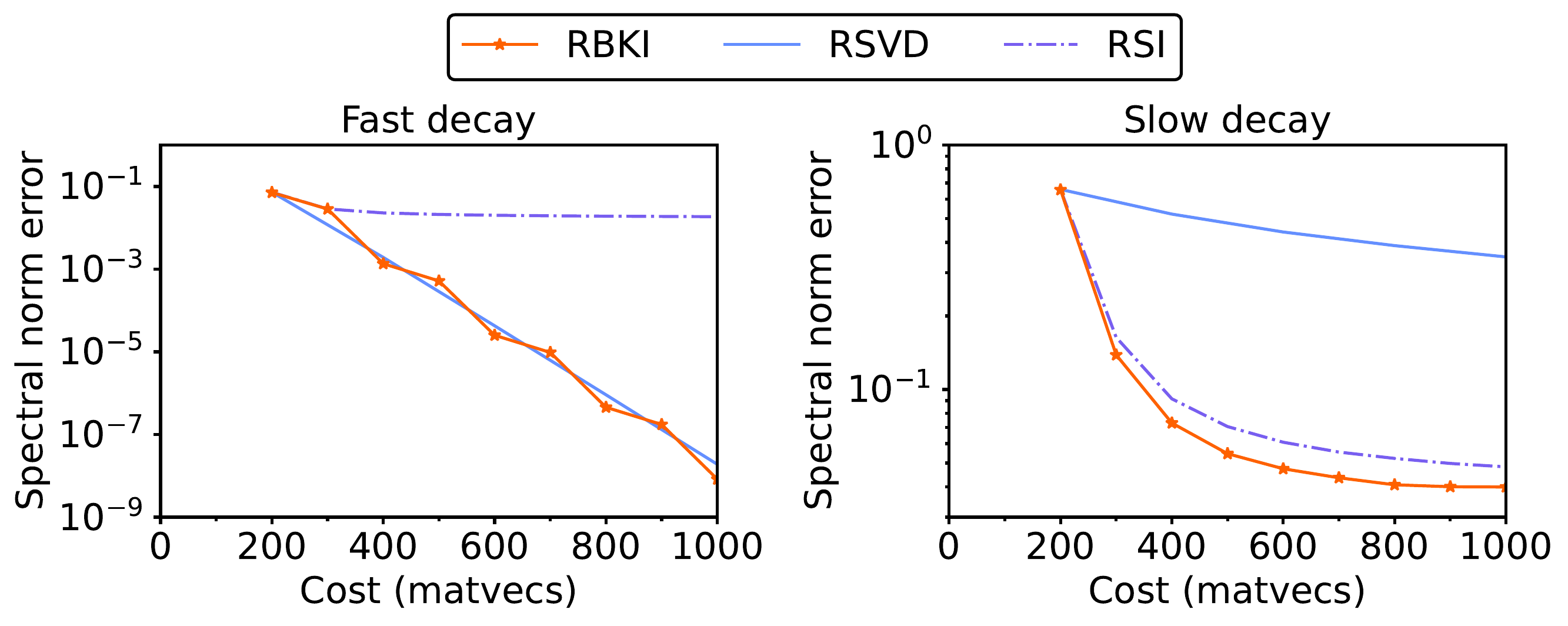}
    \includegraphics[width = \textwidth]{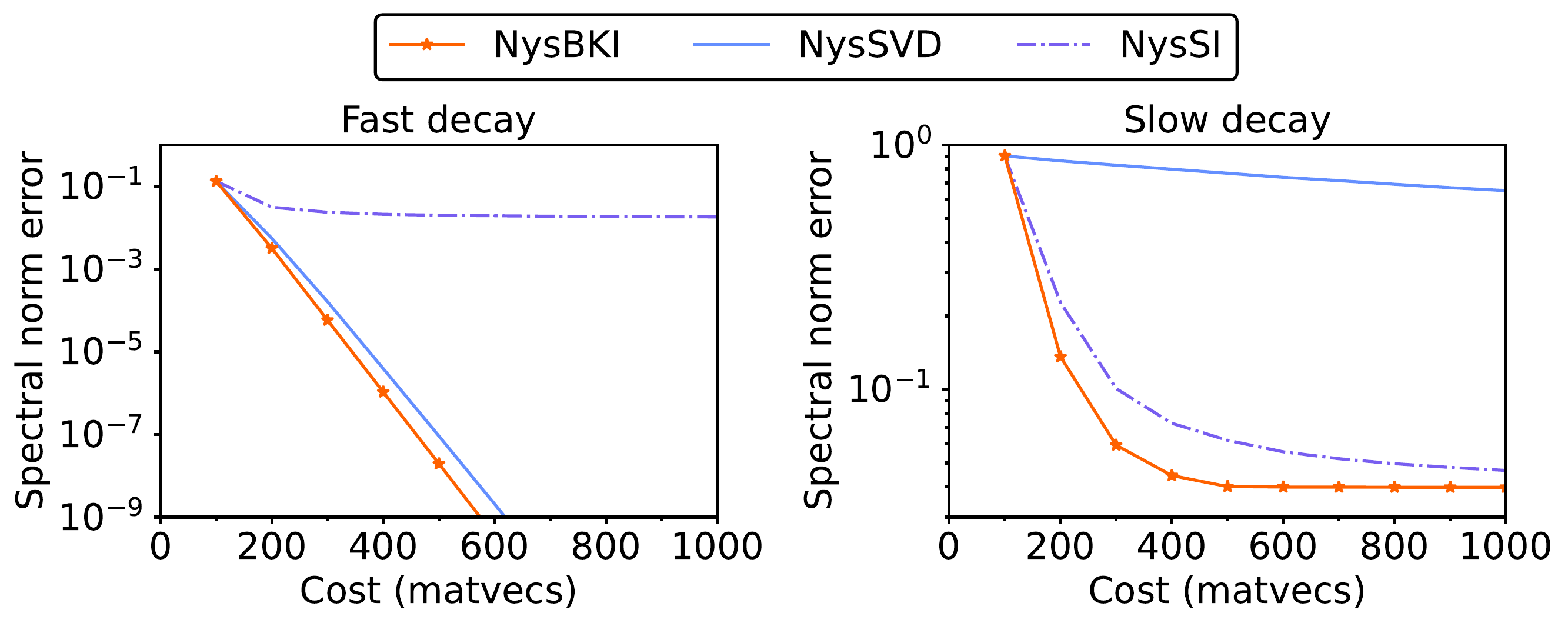}
    \caption{\textbf{(Matrix approximation comparisons).} Matrix approximation error for the matrix $\bm{A}$ with fast singular value decay (left) and the matrix $\bm{B}$ with slow singular value decay (right), as constructed in \cref{sec:comparison}.
    Algorithms for general matrices are on top, algorithms for psd matrices are on bottom.}
    \label{fig:krylov_v_subspace}
\end{figure}

\Cref{fig:krylov_v_subspace} evaluates the computational cost and approximation error of several low-rank approximation algorithms.
The horizontal axis measures the computational cost using the number of matrix--vector products $km$.
We allow the block size $k$ to vary for \texttt{RSVD} and \texttt{NysSVD}.
For the other algorithms, we fix the block size to $k = 100$ and allow the number of matrix--matrix multiplications $m$ to vary.
The vertical axis measures the root-mean-square spectral-norm approximation error
\begin{equation*}
    \bigl(\Expect \bigl\lVert \hat{\bm{A}} - \bm{A} \bigr\rVert^2\bigr)^{1/2}
    \quad \text{or} \quad
    \bigl(\Expect \bigl\lVert \hat{\bm{B}} - \bm{B} \bigr\rVert^2\bigr)^{1/2}.
\end{equation*}
The expectation is approximated empirically over $100$ independent Gaussian initializations.
Note that we are using absolute (not relative) errors here.

The results in \cref{fig:krylov_v_subspace} show that various randomized low-rank approximation algorithms, including \texttt{NysSVD} and \texttt{RSVD}, can produce a high-quality approximation of the matrix $\bm{A}$ with fast singular value decay.
Since \texttt{NysSVD} and \texttt{RSVD} use fewer matrix--matrix multiplications ($m = 1$ for \texttt{NysSVD}, $m = 2$ for \texttt{RSVD}), these strategies are more cost-efficient when approximating such a matrix.
\texttt{NysSVD} is more efficient than \texttt{RSVD} since it achieves the same approximation accuracy using half as many matrix--matrix multiplications.

In contrast, the Krylov methods yield the highest accuracy approximations for the matrix $\bm{B}$ with slowly decaying singular values.
The \texttt{NysBKI} method is the top performer, since it achieves this high approximation accuracy using a factor of $\sqrt{2}$ fewer matrix--matrix multiplications than \texttt{RBKI}.

\begin{figure}[t]
    \centering
    \includegraphics[width = \textwidth]{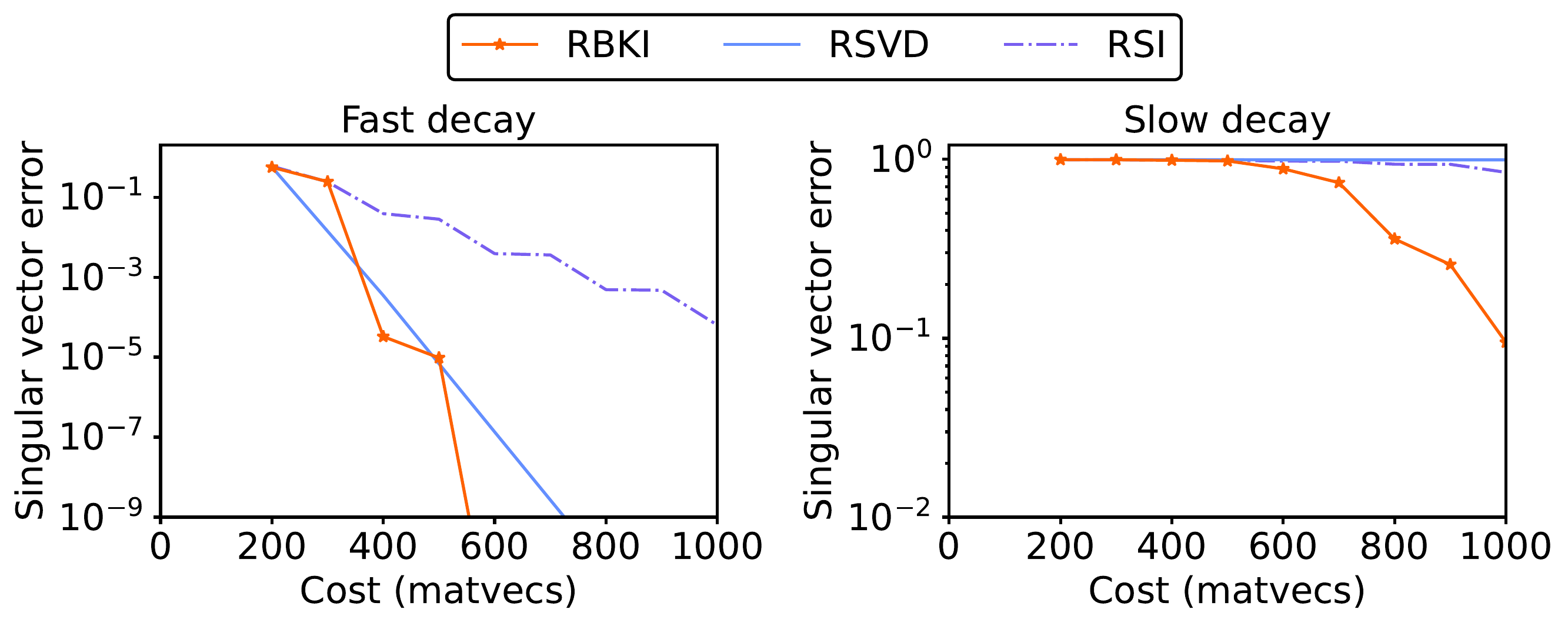}
    \includegraphics[width = \textwidth]{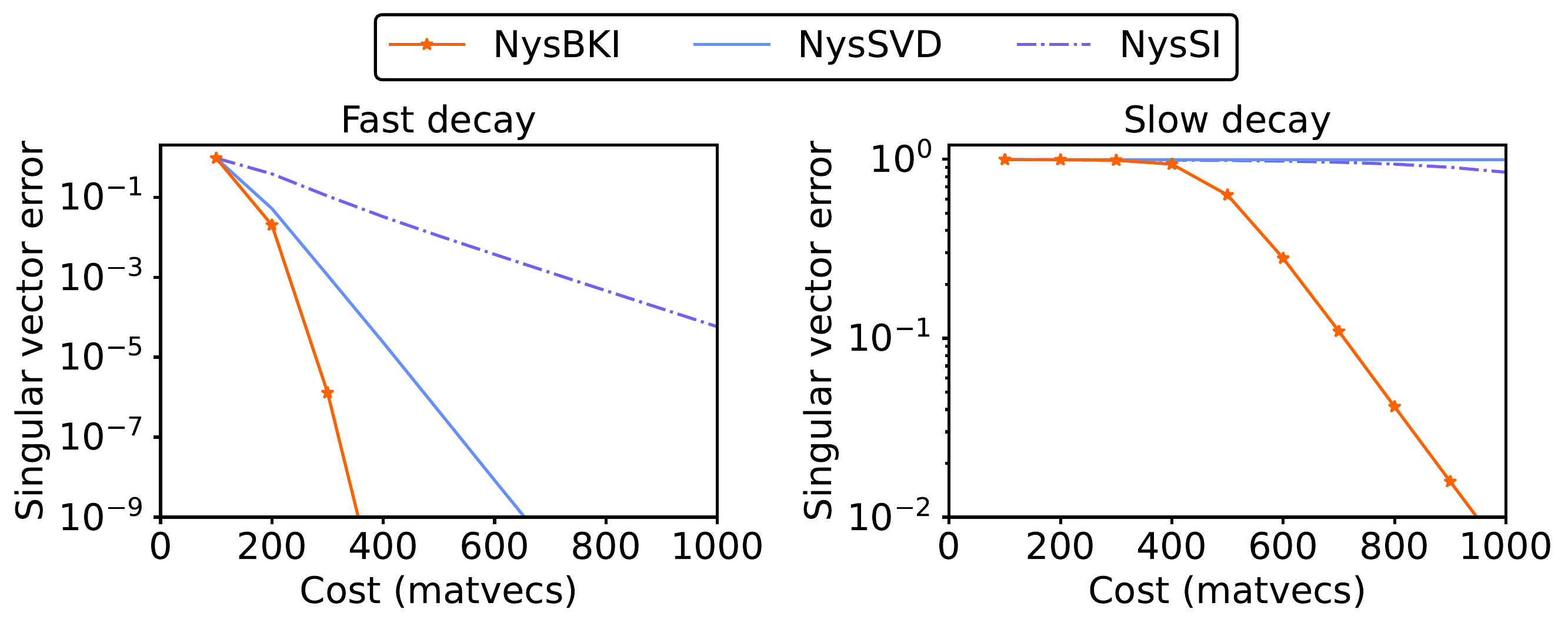}    \caption{\textbf{(Singular vector approximation comparisons).} Error in the dominant 75 singular vectors for the matrix $\bm{A}$ with fast singular value decay (left) and the matrix $\bm{B}$ with slow singular value decay (right), as constructed in \cref{sec:comparison}.
    Algorithms for general matrices are on top, algorithms for psd matrices are on bottom.}
    \label{fig:eigenvectors}
\end{figure}

Next, \cref{fig:eigenvectors} compares the singular vector approximations from several low-rank approximation algorithms.
The singular vectors are vital for principal component analysis and kernel spectral clustering, as we will show through examples in \cref{sec:genetics,sec:clustering}.
We evaluate the error in the estimated singular vectors as
\begin{equation*}
    \bigl(\Expect \bigl\lVert \bm{\Pi}_{75}(\hat{\bm{A}}) - \bm{\Pi}_{75}(\bm{A}) \bigr\rVert^2\bigr)^{1/2}
    \quad \text{or} \quad
    \bigl(\Expect \bigl\lVert \bm{\Pi}_{75}(\hat{\bm{B}}) - \bm{\Pi}_{75}(\bm{B}) \bigr\rVert^2\bigr)^{1/2},
\end{equation*}
where $\bm{\Pi}_{75}(\cdot)$ denotes the orthogonal projection onto the dominant $75$ right singular vectors.
We find that \texttt{RSVD} and \texttt{NysSVD} lead to accurate singular vector approximations for the matrix $\bm{A}$ with rapid singular value decay.
However, \texttt{RBKI} and \texttt{NysBKI} produce \textbf{10$\times$ to 300$\times$ more accurate singular vector approximations} for the matrix $\bm{B}$ than the other low-rank approximation methods.

\subsection{Human genetic diversity data} \label{sec:genetics}

Principal component analysis is frequently applied to large matrices containing single-nucleotide polymorphism (SNP) data, encoding the variations in DNA at specific locations in the genome \cite{abraham2014fast,galinsky2016fast}.
The principal components of the SNP data can be used to cluster individuals that share similar genomes.

\begin{figure}[t]
    \centering
    \includegraphics[width=.7\textwidth]{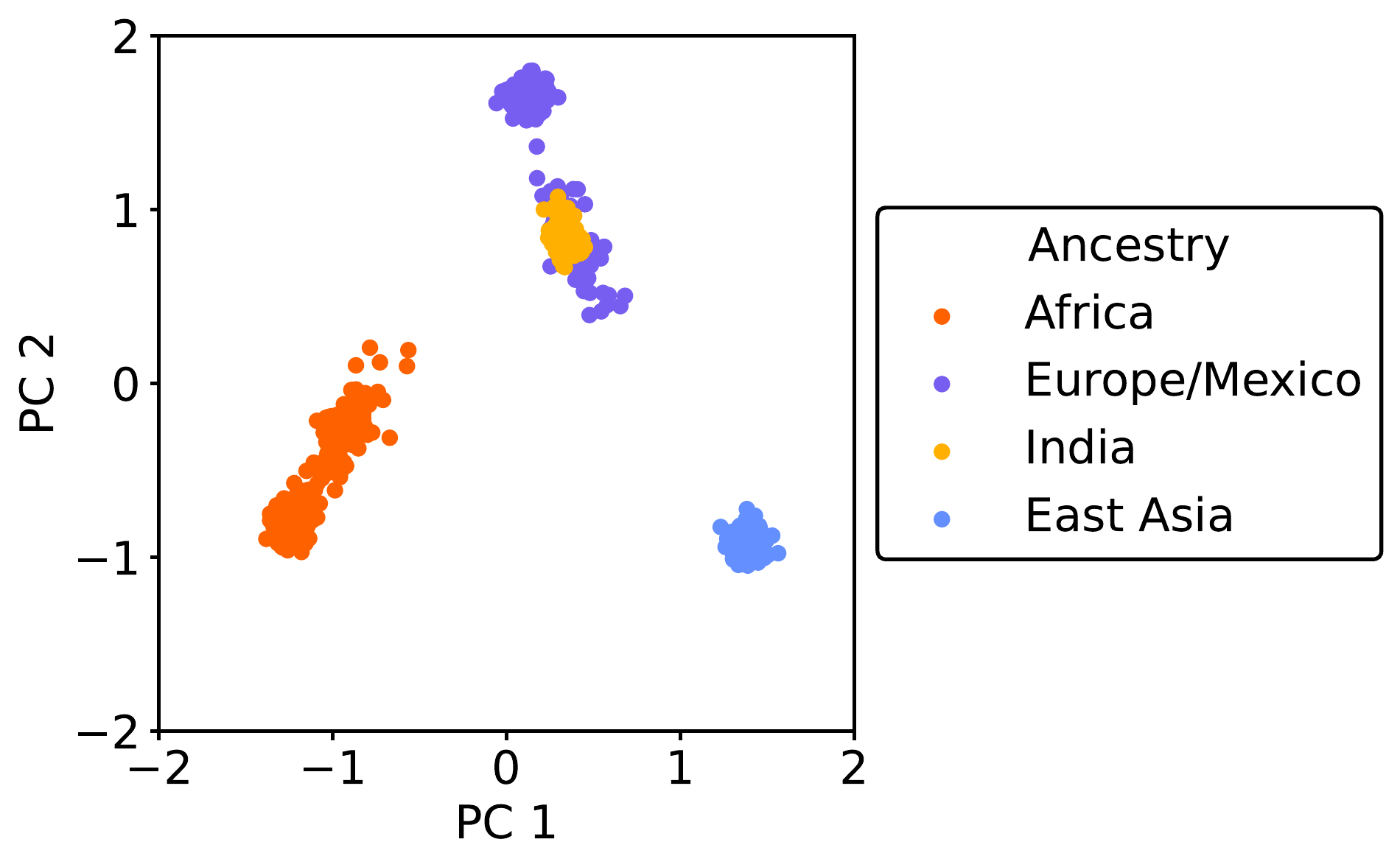}
    \caption{\textbf{(Clustering in genetics).} $k$-means clustering applied to the top 5 principal components of the HapMap3 data set, as described in \cref{sec:genetics}.
    Figure shows the projection onto the top 2 principal components.}
    \label{fig:reference}
\end{figure}

Here, we apply principal component analysis to the HapMap3 data set \cite{international2010integrating}, one of the earliest data sets measuring human genetic diversity, which we downloaded from \cite{flashpca}.
The data is organized into a matrix $\bm{A} \in \mathbb{R}^{957 \times 14,079}$ containing counts of SNPs for 957 individuals in 14,079 chromosomal locations.
The entries of $\bm{A}$ are 0, 1, or 2,
corresponding to the number of affected chromosomes.
To normalize the matrix entries, we apply the transformation
\begin{equation*}
    \bm{B}_{ij} = \frac{\bm{A}_{ij} - \mu_j}{\sqrt{\frac{1}{2}\mu_j ( 1 - \frac{1}{2} \mu_j)}},
    \quad \text{where} \quad
    \mu_j = \frac{1}{957} \sum_{i=1}^{957} \bm{A}_{ij}.
\end{equation*}
Our goal is to extract the principal components, which are defined as the right singular vectors of $\bm{B}$.
Once we have identified the principal components, we can cluster individuals by genetic ancestry as depicted in \cref{fig:reference}.

HapMap3 is a small data set by today's standards,
and we can obtain the principal components using a full singular value decomposition on a laptop.
Yet modern genetic data sets may contain millions of individuals and genetic markers, making a full singular value decomposition infeasible \cite{BKK+19}.
Therefore, we consider a more scalable approach to principal component analysis based on randomized low-rank approximation.

\begin{figure}[t]
    \centering    
    \includegraphics[width=.5\textwidth]{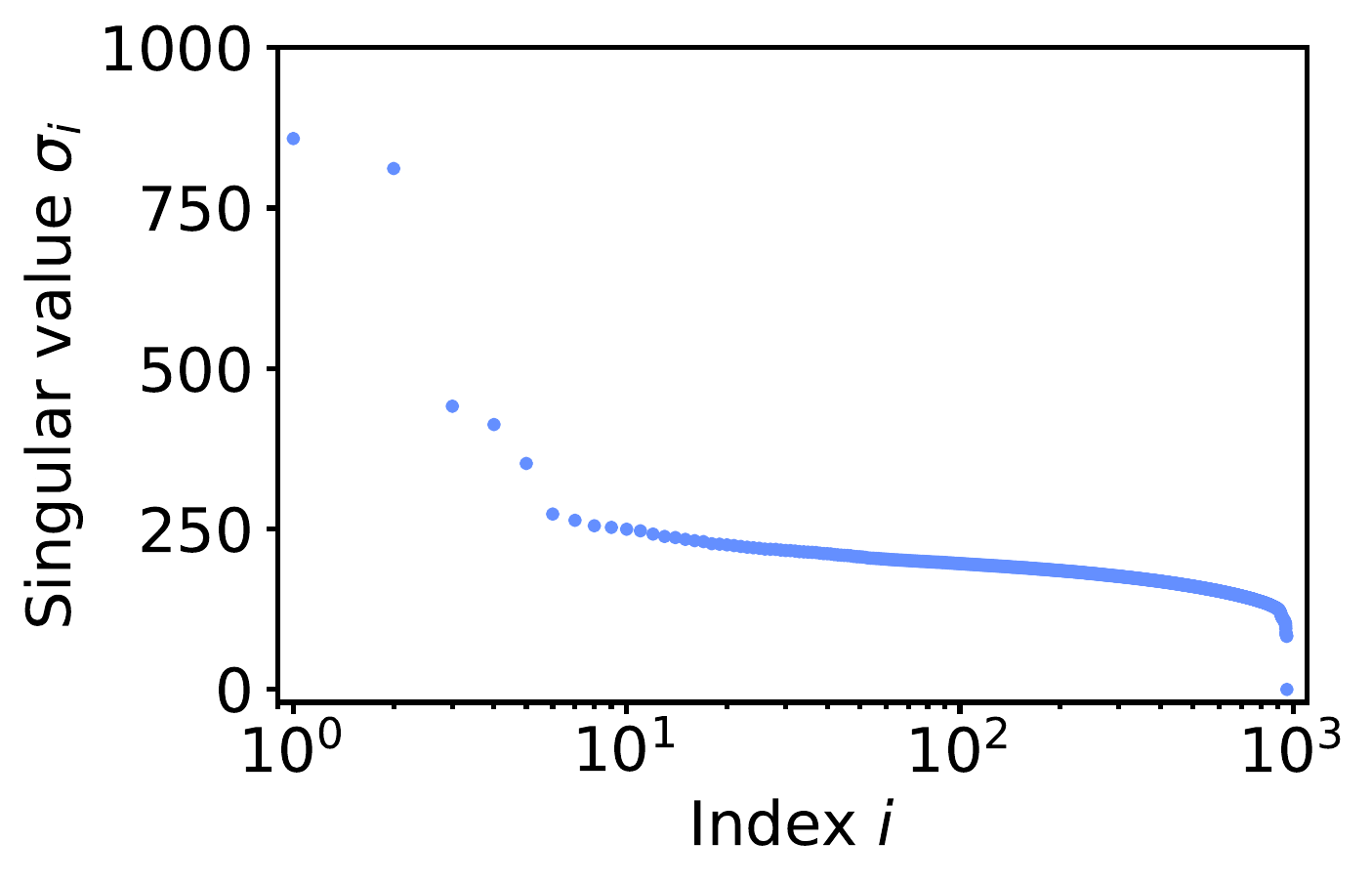}
    \caption{\textbf{(Singular value decay in genetics).} Singular value decay for the HapMap3 data set, as described in \cref{sec:genetics}.}
    \label{fig:hapmap_decay}
\end{figure}

The randomized approach to principal component analysis has been questioned in the past.
In 2013, Chen et al.~wrote \cite{chen2013improved},
\vspace{0.2cm}
\begin{quote}
``In theory, randomized eigenvector approximations [\texttt{RSVD}] can reduce the running time... (Rokhlin et al., 2009). However, a colleague of ours reports that efforts to apply this approach to genetic data have not yet been successful, as in large datasets, eigenvalues may be highly significant (reflecting real population structure in the data) but only slightly larger than background noise eigenvalues, and thus sometimes missed by randomized methods (N. Patterson, personal communication).''
\end{quote}
\vspace{0.2cm}
The authors are concerned that randomized algorithms cannot accurately identify principal components when the singular values are close together.
This would present a major limitation, since genetic data sets often have small singular value gaps, as shown for the HapMap3 data in \cref{fig:hapmap_decay}.

\begin{figure}[t]
    \centering
    \includegraphics[width=\textwidth]{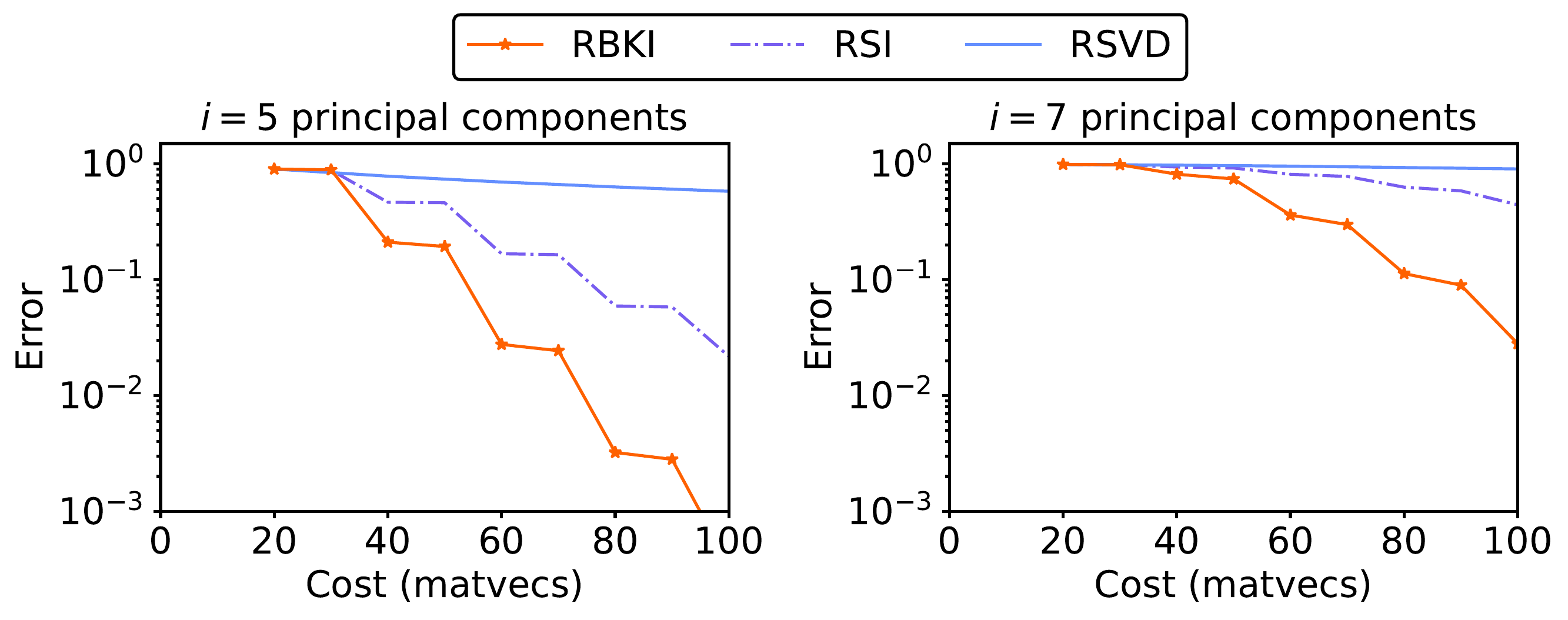}
    \caption{\textbf{(Approximating principal components in genetics).} Error in the top $i = 5$ principal components (left) or top $i = 7$ principal components (right) for the HapMap3 data, as described in \cref{sec:genetics}.}
    \label{fig:hapmap_results}
\end{figure}

\Cref{fig:hapmap_results} compares the accuracy of \texttt{RSVD}, \texttt{RSI}, and \texttt{RBKI} (the latter two with block size $k = 10$) when approximating the principal components of the HapMap3 data.
The vertical axis evaluates error in the principal components using 
\begin{equation*}
    \bigl(\Expect \bigl\lVert \bm{\Pi}_i(\hat{\bm{B}}) - \bm{\Pi}_i(\bm{B}) \bigr\rVert^2\bigr)^{1/2},
\end{equation*}
where $\bm{\Pi}_{i}(\cdot)$ represents the orthogonal projection onto the top $i = 5$ or $i = 7$ right singular vectors and the expectation is measured empirically over $100$ independent Gaussian initializations.
The top $i = 5$ principal components (left column) are separated by a sizable spectral gap,
whereas the top $i = 7$ principal components (right column) are separated by a tiny singular value gap leading to larger errors.
The results in \cref{fig:hapmap_results} confirm that \texttt{RSVD} and \texttt{RSI} struggle to identify the top $i = 7$ principal components in the data.
Only \texttt{RBKI} identifies
all 7 principal components up to a precision of $\varepsilon = 0.1$.

\texttt{RBKI} also performs well at separating individuals into clusters.
With just $km = 40$ matrix--vector products, the algorithm identifies the top $i = 5$ principal components well enough to match the ideal clustering results in \cref{fig:reference}.
In contrast, we would need to perform \texttt{RSVD} with \textbf{20$\times$ as many matvecs} ($km = 800$) to achieve the same clustering accuracy.

These investigations show that \texttt{RBKI} is both fast and accurate when performing principal component analysis with genetic data.
\texttt{RBKI} is faster than traditional singular value decomposition, since we can take $km$ to be a small fraction of the number of rows.
Moreover \texttt{RBKI} is highly accurate, even when the signal barely rises above the noise.

\subsection{Spectral clustering} \label{sec:clustering}

Kernel spectral clustering is a powerful approach for analyzing the dynamics of proteins using data from molecular dynamics simulations \cite{GHR21}.
Spectral clustering is useful for identifying metastable states that the protein occupies for a long time, with rare transitions between states (e.g., folded and unfolded states).
However in the past, kernel spectral clustering has mainly been limited to data sets with $N \leq 10^4$ data points, because it requires calculating the dominant eigenvectors of an $N \times N$ matrix.

Here, we use \texttt{NysBKI} to extend kernel spectral clustering to a large data set with $N = 250{,}000$ points.
The data comes from a 250~ns simulation of alanine dipeptide (\ce{CH_3-CO-NH-C_{\alpha}HCH_3-CO-NH-CH_3}), which we downloaded from \cite{ala2}.
The data points $\bm{x}^{(i)} \in \mathbb{R}^{30}$ identify the spatial $(x, y, z)$-positions of the 10 non-hydrogen atoms at 1~ps intervals.
We will perform kernel spectral clustering as follows.
\vspace{0.2cm}
\begin{enumerate}
\item Form the psd Gaussian kernel matrix $\bm{A} \in \mathbb{R}^{N \times N}$ with entries 
\begin{equation*}
\label{eq:gaussian}
a_{ij} = \exp\biggl(-\frac{1}{2\sigma^2} \lVert \bm{x}^{(i)} - \bm{x}^{(j)} \rVert^2\biggr),
\end{equation*}
where $\sigma > 0$ is a tunable bandwidth.
\item Form the diagonal matrix $\bm{D} \in \mathbb{R}^{N \times N}$ containing the row sums of $\bm{A}$.
\item Calculate the dominant $r$ eigenvectors $\bm{V} = \begin{bmatrix} \bm{v}^{(1)} & \cdots & \bm{v}^{(r)} \end{bmatrix}$ for the psd matrix $\bm{D}^{-1/2} \bm{A} \bm{D}^{-1/2}$,
where $r$ is a tunable parameter.
\item Apply the diagonal rescaling $\bm{V} \leftarrow \bm{D}^{-1/2} \bm{V}$, so that $\bm{V}$ contains right eigenvectors of $\bm{D}^{-1} \bm{A}$.
\item Apply $k$-means clustering to the rows of $\bm{V}$.
\end{enumerate}
\vspace{0.2cm}
\begin{figure}[t]
    \centering
    \includegraphics[width=.35\textwidth]{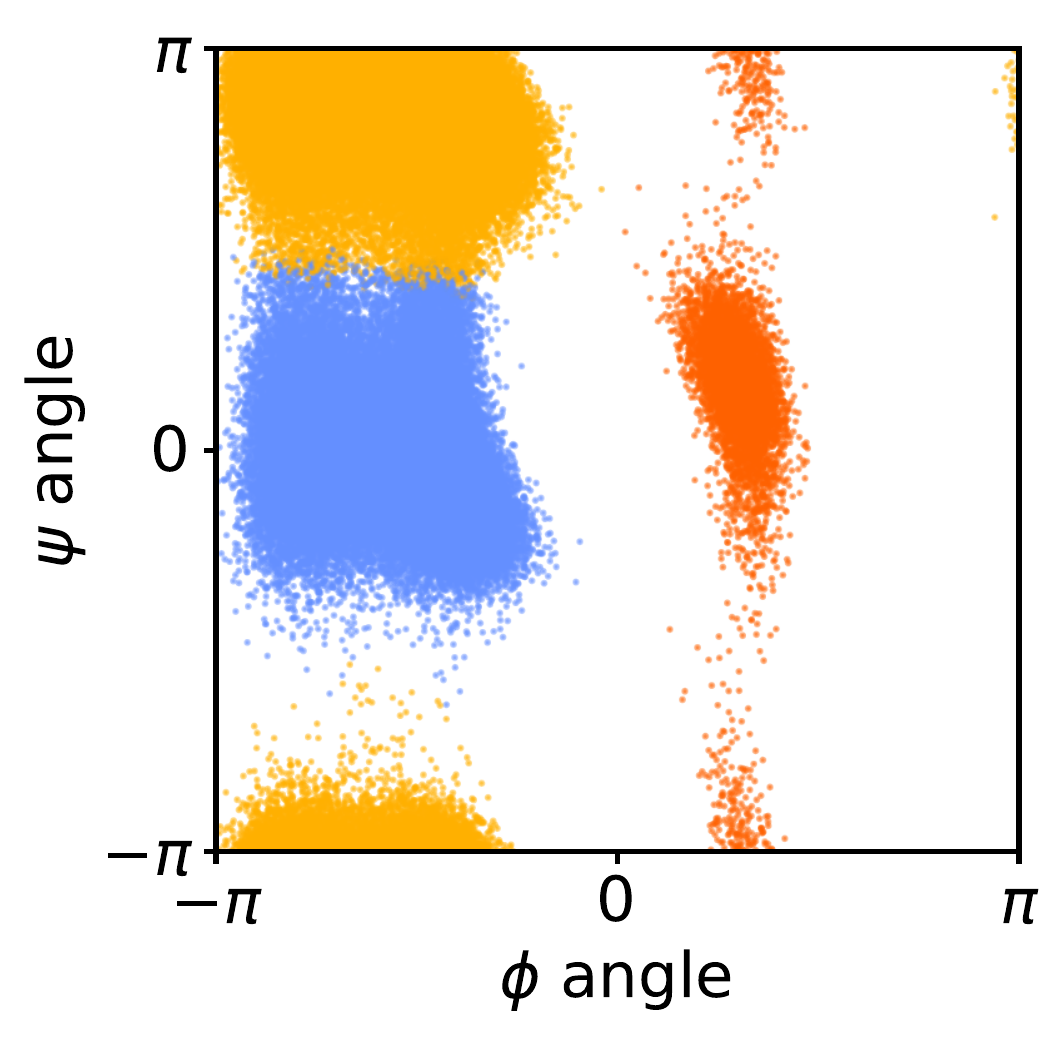}
    \caption{\textbf{(Clustering in biochemistry).} Kernel spectral clustering with $k = 3$ clusters, $r = 3$ eigenvectors, and a bandwidth $\sigma = .05$nm, applied to the alanine dipeptide data set ($N = 250{,}000$, $d = 30$) as described in \cref{sec:clustering}.
    Figure shows the projection onto the $\phi$ and $\psi$ dihedral angles.}
    \label{fig:adp_clusters}
\end{figure}

The spectral clustering  algorithm identifies three clusters in the 30-dimensional alanine dipeptide data space.
These clusters are highly correlated with the dihedral angle $\phi$ between 
\ce{C}, \ce{N}, \ce{C_{\alpha}}, and \ce{C} and $\psi$ 
between \ce{N}, \ce{C_{\alpha}}, \ce{C}, and \ce{N}, as shown in 
\cref{fig:adp_clusters}.
Yet we emphasize that the algorithm does not have access to the dihedral angles: it discovers the clusters organically, reproducing the past work \cite{NWPW+17}.

It is computationally challenging to complete the third step of kernel spectral clustering, which requires calculating the dominant eigenvectors of an $N \times N$ psd matrix.
For the alanine dipeptide problem, the matrix $\bm{B} = \bm{D}^{-1/2} \bm{A} \bm{D}^{-1/2}$ requires 222GB storage and
is therefore too large to fit in working memory on a 64GB laptop.
Therefore, to produce the results in \cref{fig:adp_clusters}, we store the matrix on disk and use \texttt{NysBKI} to approximate the leading eigenvectors.
With parameter settings $k = 100$ and $m = 2$, \texttt{NysBKI} runs in just $11$ minutes and $99.9\%$ of the time is spent loading blocks of the matrix and performing blockwise matrix multiplications.

\begin{figure}[t]
    \centering
    \includegraphics[width = .5\textwidth]{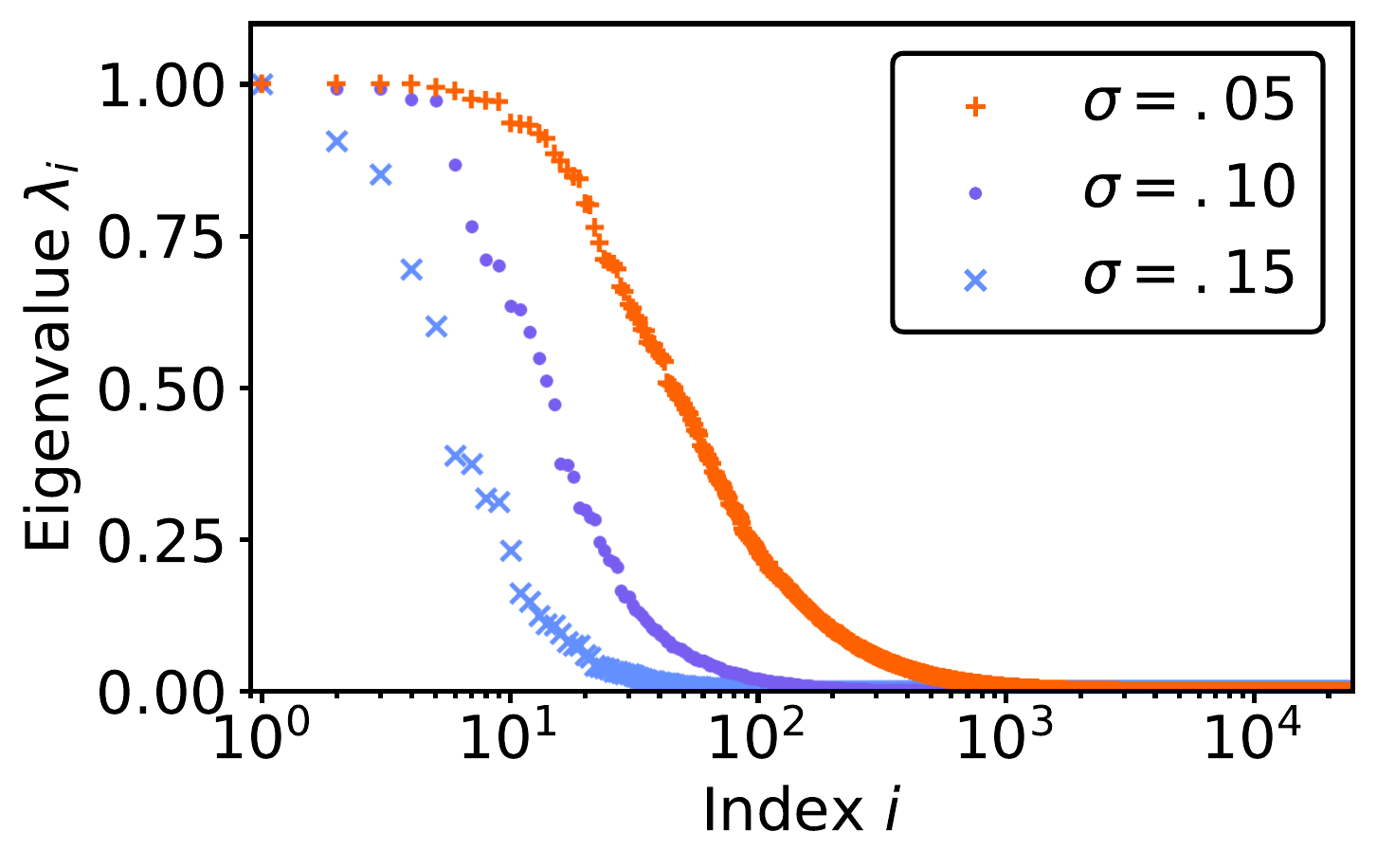}
    \caption{\textbf{(Eigenvalue decay in biochemistry).} Eigenvalues of the kernel matrix for a subset of $N = 25{,}000$ alanine dipeptide data points, as described in \cref{sec:clustering}.
    The kernel matrix is evaluated with bandwidth parameter $\sigma = .05$, $.10$, or $.15$nm.}
    \label{fig:alanine_evals}
\end{figure}

\begin{figure}[t]
    \centering
    \includegraphics[width=\textwidth]{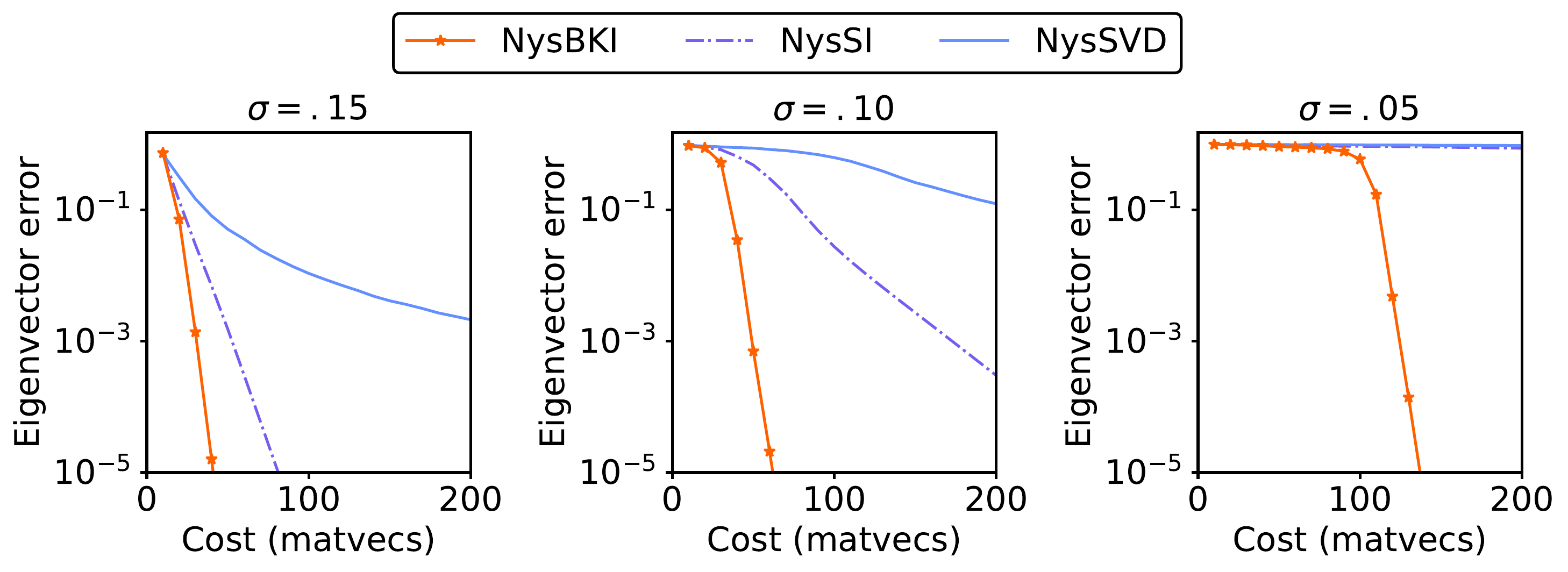}
    \caption{\textbf{(Approximating eigenvectors in biochemistry).} Error in the top 3 eigenvectors of the kernel matrix for a subset of $N = 25{,}000$ alanine dipeptide data points, as described in \cref{sec:clustering}. 
    The kernel matrix is evaluated with bandwidth parameter $\sigma = .15$ (left), $\sigma = .10$ (center) or $\sigma = .05$nm (right).}
    \label{fig:cluster_test}
\end{figure}

To provide additional insight into the difficulties of kernel spectral clustering, we take a subset of $N = 25{,}000$ equally spaced data points (so we can perform a full eigendecomposition as a reference) and consider the impact of varying the bandwidth $\sigma$.
A small bandwidth $\sigma$ is needed to obtain physically relevant clusters with sharp boundaries between clusters.
Unfortunately, however, the eigenvalues decay more slowly when $\sigma$ is small (\cref{fig:alanine_evals}), leading to a difficult approximation problem.

\Cref{fig:cluster_test} presents the eigenvector approximation error when \texttt{NysSVD}, \texttt{NysSI}, and \texttt{NysBKI} (the latter two with block size $k = 10$) are applied to the subset of $N = 25{,}000$ alanine dipeptide data points.
We measure the eigenvector approximation error using
\begin{equation*}
    \bigl(\Expect \bigl\lVert \bm{\Pi}_3(\hat{\bm{B}}) - \bm{\Pi}_3(\bm{B}) \bigr\rVert^2\bigr)^{1/2}
\end{equation*}
where $\bm{\Pi}_3$ is the orthogonal projection onto the dominant $3$ eigenvectors, and the expectation is evaluated empirically over $100$ independent runs.
The results show that a relatively small number of matvecs ($km<100$) are sufficient to approximate the top 3 eigenvectors when $\sigma = .15$ (left) and even when $\sigma = .10$ (center).
However, a larger number of matvecs is needed when $\sigma = .05$ (right), which is the physically relevant parameter setting for alanine dipeptide.
The accuracy is dramatically higher when using \texttt{NysBKI}, instead of \texttt{RSVD} or \texttt{RSI}.

In summary, we find that kernel spectral clustering leads to challenging, high-dimensional eigenvalue problems, especially when the bandwidth $\sigma$ is small.
To solve these eigenvalue problems with high accuracy and scalability, we recommend the \texttt{NysBKI} algorithm.
\texttt{NysBKI} works robustly across different bandwidth settings,
and it leads to dramatic computational speedups.
With $N = 250{,}000$ data points and an approximation rank of $km = 200$, the difference between the traditional $\mathcal{O}(N^3)$ operation count and \texttt{NysBKI}'s $\mathcal{O}(km N^2)$ operation count is \textbf{over 3 orders of magnitude}.

\section{Analysis of \texttt{RSVD} and \texttt{NysSVD}} \label{sec:theory1}

In this section, we derive error bounds for \texttt{RSVD} (\cref{alg:rsvd}) and \texttt{NysSVD} (\cref{alg:nyssvd}), which are the simplest randomized low-rank approximation algorithms.  These results provide the foundation for studying the more sophisticated methods (\texttt{RSI}, \texttt{RBKI}, \texttt{NysSI}, \texttt{NysBKI}).  Although related results are already well-established~\cite{halko2011finding}, we have found opportunities for simplification and improvement.

We compare the \texttt{RSVD} approximation to an optimal rank-$r$ approximation $\lfloor \bm{A} \rfloor_r$,
which arises from an $r$-truncated singular value decomposition, and we establish the following main result.  The proof appears below in \cref{sec:random}.

\begin{theorem}[\texttt{RSVD} error] \label{thm:main_thm}
Fix a matrix $\bm{A} \in \mathbb{R}^{L \times N}$
and a target rank $r \geq 1$.
For any $2 \leq p \leq \infty$ and $u,t \geq 0$,
\texttt{RSVD} with $k \geq r$ Gaussian initialization vectors
generates a random rank-$k$ approximation $\hat{\bm{A}} \in \mathbb{R}^{L \times N}$ that satisfies 
\begin{equation}
\label{eq:aesthetic}
\tag{\texttt{RSVD}1}
    \lVert \bm{A} - \hat{\bm{A}} \rVert_p^2 \leq
    \lVert \bm{A} - \lfloor \bm{A} \rfloor_r \rVert_p^2
    + \frac{utr}{k - r + 1}
    \lVert \bm{A} - \lfloor \bm{A} \rfloor_r \rVert_{\rm F}^2
\end{equation}
with failure probability at most $\mathrm{e}^{-(u-2) \slash 4} + \sqrt{\pi r} (t \slash \mathrm{e})^{-(k-r+1) \slash 2}$.

Additionally, if the block size $k$ is no smaller than $r + 2$, then
\begin{equation}
\label{eq:simple}
\tag{\texttt{RSVD}2}
    \Expect \lVert \bm{A} - \hat{\bm{A}} \rVert_p^2 \leq
    \lVert \bm{A} - \lfloor \bm{A} \rfloor_r \rVert_p^2 + \frac{r}{k - r - 1} \lVert \bm{A} - \lfloor \bm{A} \rfloor_r \rVert_{\rm F}^2.
\end{equation}
\end{theorem}

The result for \texttt{NysSVD} is a corollary
of the result for \texttt{RSVD}.

\begin{corollary}[\texttt{NysSVD} error] \label{cor:main_cor}
Fix a psd matrix $\bm{A} \in \mathbb{R}^{N \times N}$ and
a target rank $r \geq 1$.
For any $1 \leq p \leq \infty$ and $u, t \geq 0$, 
\texttt{NysSVD} with $k \geq r$ Gaussian initialization vectors
generates a random rank-$k$ approximation $\hat{\bm{A}} \in \mathbb{R}^{N \times N}$ that satisfies 
\begin{equation}
\label{eq:aesthetic2}
\tag{\texttt{NysSVD}1}
    \lVert \bm{A} - \hat{\bm{A}} \rVert_p \leq
    \lVert \bm{A} - \lfloor \bm{A} \rfloor_r \rVert_p + \frac{utr}{k - r + 1}
    \lVert \bm{A} - \lfloor \bm{A} \rfloor_r \rVert_{\ast}
\end{equation}
with failure probability at most $\mathrm{e}^{-(u-2) \slash 4} + \sqrt{\pi r} (t \slash \mathrm{e})^{-(k-r+1) \slash 2}$.

Additionally, if the block size $k$ is no smaller than $r + 2$, then
\begin{equation}
\label{eq:simple2}
\tag{\texttt{NysSVD}2}
    \Expect \lVert \bm{A} - \hat{\bm{A}} \rVert_p \leq
    \lVert \bm{A} - \lfloor \bm{A} \rfloor_r \rVert_p + \frac{r}{k - r - 1} \lVert \bm{A} - \lfloor \bm{A} \rfloor_r \rVert_{\ast}.
\end{equation}
\end{corollary}

The error bounds in \cref{thm:main_thm,cor:main_cor} only depend on the singular values $\sigma_i(\bm{A})$ for $i \geq r+1$.
If these singular values are small and rapidly decaying, the theorem immediately establishes the accuracy of the \texttt{RSVD} and \texttt{NysSVD} approximations.

Interpreting these error bounds, we identify two limitations of \texttt{RSVD} and \texttt{NysSVD}.
First, when the singular values decay slowly, the terms $\lVert \bm{A} - \lfloor \bm{A} \rfloor_r \rVert_{\rm F}$ and $\lVert \bm{A} - \lfloor \bm{A} \rfloor_r \rVert_\ast$ are large, in which case \texttt{RSVD} and \texttt{NysSVD} perform poorly, as confirmed by our experiments in \cref{sec:experiments}.
To target matrices with slow singular value decay, we need more powerful algorithms such as \texttt{RBKI} and \texttt{NysBKI}.
As another limitation, the expectation bounds are not applicable when $k = r$ or $k = r + 1$,
and large random errors can occur in these settings.
When $k - r$ is small, there is a fundamental problem that the range of $k$ random Gaussian initialization vectors might barely align with the leading $r$ singular vectors.

\texttt{RSVD} and \texttt{NysSVD} were previously analyzed in \cite{halko2011finding,gu2015subspace,tropp2017fixed}.
However,
by revisiting and reworking the details of these analyses,
we obtain a more flexible framework that delivers many new bounds.
The probability bounds \cref{eq:aesthetic,eq:aesthetic2} are especially novel and informative, since previous analyses \cite{gu2015subspace,tropp2022randomized} separated the cases $k = r$, $k = r + 1$, and $k \geq r + 2$ whereas we provide a unified treatment.
Additionally, our analytic approach leads to a broad range of expectation bounds for \texttt{RSVD} and \texttt{NysSVD}, many of them new, to be presented in \cref{sec:discussion}.

The new technical idea in this section is to decompose the Schatten $p$-norm error
of \texttt{RSVD} using \emph{parallel sums}
of psd matrices \cite{anderson1969series}.
Parallel sums, which generalize the harmonic mean to the case of psd matrices, are sufficiently tractable that they allow us to produce a wide range of probability and expectation bounds.
To our knowledge, parallel sums have not been previously exploited in the randomized matrix approximation literature.

The rest of this section contains a complete and up-to-date analysis of \texttt{RSVD} and \texttt{NysSVD}.
\Cref{sec:diagonal} presents a reduction to the case of diagonal, psd matrices;
\cref{sec:essential} analyzes the error of \texttt{RSVD} using parallel sums;
\cref{sec:random} uses random matrix theory to prove our main error bounds;
and \cref{sec:discussion} derives additional expectation bounds for \texttt{RSVD} and \texttt{NysSVD}.
Afterward, \cref{sec:theory2} deploys these results to
describe the behavior
of \texttt{RSI}, \texttt{NysSI}, \texttt{RBKI}, and \texttt{NysBKI}.

\subsection{Diagonal, psd reduction}
\label{sec:diagonal}

As the first step in our analysis, we show that the error of many randomized low-rank approximation algorithms depends only on the singular values of the target matrix, regardless of the singular vectors.  This reduction is standard.

\begin{lemma}[Diagonal, psd reduction] \label{lem:diagonal}
When we apply \texttt{RSVD}, \texttt{RSI}, \texttt{RBKI}, \texttt{NysSVD}, \texttt{NysSI}, or \texttt{NysBKI} to approximate a matrix $\bm{A} \in \mathbb{R}^{L \times N}$ using a Gaussian initialization matrix $\bm{\Omega} \in \mathbb{R}^{N \times k}$, 
the Schatten $p$-norm error $\lVert \bm{A} - \hat{\bm{A}} \rVert_p$ has a distribution that only depends on the singular values of $\bm{A}$, regardless of the singular vectors, for any $1 \leq p \leq \infty$.
\end{lemma}
\begin{proof}
Let $\bm{U} \bm{\Sigma} \bm{V}^\ast$ be a singular value decomposition for $\bm{A}$ and set $\bm{\Omega}^\prime = \bm{V}^\ast \bm{\Omega}$.
The randomized low-rank approximation $\hat{\bm{A}}$ can be written:
\begin{align*}
\tag{\texttt{RSVD}}
    & \bm{U} \bigl(\bm{\Pi}_{\bm{\Sigma} \bm{\Omega}^\prime} \bm{\Sigma} \bigr) \bm{V}^\ast, \\
\tag{\texttt{RSI}, even $m$}
    & \bm{U} \bigl(\bm{\Pi}_{\bm{\Sigma}^{m - 1} \bm{\Omega}^\prime} \bm{\Sigma} \bigr) \bm{V}^\ast, \\
\tag{\texttt{RSI}, odd $m$}
    & \bm{U} \bigl(\bm{\Sigma} \bm{\Pi}_{\bm{\Sigma}^{m - 1} \bm{\Omega}^\prime} \bigr) \bm{V}^\ast, \\
\tag{\texttt{RBKI}, even $m$}
    & \bm{U} \bigl(\bm{\Pi}_{[\begin{smallmatrix}
    \bm{\Sigma}^{m-1} \bm{\Omega}^\prime & \bm{\Sigma}^{m-3} \bm{\Omega}^\prime & \cdots & \bm{\Sigma} \bm{\Omega}^\prime
    \end{smallmatrix}]
    } \bm{\Sigma} \bigr) \bm{V}^\ast, \\
\tag{\texttt{RBKI}, odd $m$}
    & \bm{U} \bigl(\bm{\Sigma} \bm{\Pi}_{[\begin{smallmatrix}
    \bm{\Sigma}^{m-1} \bm{\Omega}^\prime & \bm{\Sigma}^{m-3} \bm{\Omega}^\prime & \cdots & \bm{\Sigma}^2 \bm{\Omega}^\prime
    \end{smallmatrix}]
    } \bigr) \bm{V}^\ast,
\end{align*}
The Schatten $p$-norm is unitarily invariant, so the $p$-norm error only depends on $\bm{\Sigma}$ and the matrix $\bm{\Omega}^\prime = \bm{V}^\ast \bm{\Omega}$, which is a random Gaussian matrix.

Similarly, if $\bm{A} \in \mathbb{R}^{N \times N}$ is psd, let $\bm{U} \bm{\Sigma} \bm{U}^\ast$ be an eigenvalue decomposition for $\bm{A}$ and set $\bm{\Omega}^\prime = \bm{U}^\ast \bm{\Omega}$.
The randomized low-rank approximation $\hat{\bm{A}}$ can be written:
\begin{align*}
\tag{\texttt{NysSVD}}
    & \bm{U} \bigl(\bm{\Sigma}^{1 \slash 2} \bm{\Pi}_{\bm{\Sigma}^{1 \slash 2} \bm{\Omega}^\prime} 
    \bm{\Sigma}^{1 \slash 2}\bigr) \bm{U}^\ast, \\
\tag{\texttt{NysSI}}
    & \bm{U} \bigl(\bm{\Sigma}^{1 \slash 2} \bm{\Pi}_{\bm{\Sigma}^{m - 1 \slash 2} \bm{\Omega}^\prime} 
    \bm{\Sigma}^{1 \slash 2}\bigr) \bm{U}^\ast, \\
\tag{\texttt{NysBKI}}
    & \bm{U} \bigl(\bm{\Sigma}^{1 \slash 2} \bm{\Pi}_{[\begin{smallmatrix}
    \bm{\Sigma}^{m-1\slash 2} \bm{\Omega}^\prime & \bm{\Sigma}^{m-3\slash 2} \bm{\Omega}^\prime & \cdots & \bm{\Sigma}^{1 \slash 2} \bm{\Omega}^\prime
    \end{smallmatrix}]
    } \bm{\Sigma}^{1 \slash 2}\bigr) \bm{U}^\ast.
\end{align*}
The same conclusion follows as before: the error only depends on $\bm{\Sigma}$ and the Gaussian matrix $\bm{\Omega}^\prime = \bm{U}^\ast \bm{\Omega}$.
\end{proof}
As a consequence of \cref{lem:diagonal} and without loss of generality, the subsequent analysis will assume that we are approximating a matrix $\bm{A} \in \mathbb{R}^{N \times N}$ that is diagonal and psd, with nonincreasing diagonal entries.

\subsection{Deterministic analysis using parallel sums}
\label{sec:essential}

In this section, we develop a simple formula for the error of \texttt{RSVD} when approximating a diagonal, psd matrix.
We fix the target matrix $\bm{A} \in \mathbb{R}^{N \times N}$, fix the initialization matrix $\bm{\Omega} \in \mathbb{R}^{N \times k}$, and analyze the $p$-norm error $\lVert \bm{A} - \bm{\Pi}_{\bm{A} \bm{\Omega}} \bm{A} \rVert_p$.
For this section, it does not matter whether the test matrix $\bm{\Omega}$ is deterministic or random.

Our analysis is based on \emph{parallel sums} of psd matrices.
Parallel sums were introduced by Anderson \& Duffin \cite{anderson1969series}
to give a mathematical framework for analyzing electric networks of capacitors and resistors.
Parallel sums have become a favorite topic in linear algebra textbooks (e.g., \cite[Ch.~4]{bhatia2009positive}), and they satisfy the following simple and natural properties:

\begin{lemma}[Parallel sums] \label{lem:parallel_sum}
For psd matrices $\bm{A}, \bm{B}$ with the same dimension,
define the \emph{parallel sum} $\bm{A}:\bm{B}$ to be the psd matrix
\begin{equation*}
\bm{A}:\bm{B} = \bm{A} - \bm{A}(\bm{A} + \bm{B})^{\dagger} \bm{A}.
\end{equation*}
Then the following properties hold:
\vspace{0.2cm}
\begin{enumerate}
    \item The parallel sum is symmetric; that is, $\bm{A}:\bm{B} = \bm{B}:\bm{A}$.
    \item If $\bm{A}$ and $\bm{B}$ are strictly positive definite, then
    $\bm{A}:\bm{B} = (\bm{A}^{-1} + \bm{B}^{-1})^{-1}$.
    \item The mapping
    $\bm{B} \mapsto \bm{A}:\bm{B}$ is monotone and concave with respect to the psd ordering, and it is bounded above as $\bm{A}:\bm{B} \preccurlyeq \bm{A}$.
    \item
    $\lVert \bm{A}:\bm{B} \rVert_p \leq \lVert \bm{A} \rVert_p : \lVert \bm{B} \rVert_p$ for any Schatten $p$-norm with $1 \leq p \leq \infty$.
\end{enumerate}
\end{lemma}
\begin{proof}
Anderson \& Duffin \cite{anderson1969series} established properties (1)--(3) and property (4) in the special cases $p = 1$ and $p = \infty$.
Ando extended property (4) to the general case $1 \leq p \leq \infty$; see \cite[eq.~3.13]{And94}.
\end{proof}

Next, we use parallel sums to analyze the \texttt{RSVD} error.
Our error decomposition is similar to \cite[Theorem 9.1]{halko2011finding},
but parallel sums allow us to highlight the main reasoning
and shorten the proof.

\begin{proposition}[Error decomposition] \label{prop:ext_analysis}
Fix a diagonal, psd matrix $\bm{A} \in \mathbb{R}^{N \times N}$ and a test matrix $\bm{\Omega} \in \mathbb{R}^{N \times k}$.
Partition the matrices as
\begin{equation*}
    \bm{A} = \begin{bmatrix} \bm{A}_1 & \\
    & \bm{A}_2 \end{bmatrix},
    \qquad 
    \bm{\Omega} = \begin{bmatrix}
    \bm{\Omega}_1 \\ \bm{\Omega}_2
    \end{bmatrix},
\end{equation*}
so that $\bm{A}_1$ is a $r \times r$ matrix and $\bm{\Omega}_1$ is a $r \times k$ matrix with $r \leq k$.  Assume that $\textup{rank}(\bm{\Omega}_1) = r$.
Then the orthogonal projection $\bm{Q} = \bm{\Pi}_{\bm{A} \bm{\Omega} \bm{\Omega}_1^{\dagger}}$ satisfies
\begin{align}
    \lVert \bm{A} - \bm{Q} \bm{A} \rVert_p^2
    & \leq \lVert \bm{A}_2 \rVert_p^2
    + \lVert \bm{A}_1 \rVert_p^2 : \lVert \bm{A}_2 \bm{\Omega}_2 \bm{\Omega}_1^{\dagger} \rVert_p^2 \\
\label{eq:satisfies}
    & \leq \lVert \bm{A}_2 \rVert_p^2
    + \lVert \bm{A}_2 \bm{\Omega}_2 \bm{\Omega}_1^{\dagger} \rVert_p^2.
\end{align}
for any Schatten $p$-norm with $2 \leq p \leq \infty$.
\end{proposition}
\begin{proof}
We make the following calculation:
\begin{align} \allowdisplaybreaks
\label{eq:calc1}
    & \lVert (\mathbf{I} - \bm{Q}) \bm{A} \rVert_p^2
    = \lVert (\mathbf{I} - \bm{Q}) \bm{A}^2 (\mathbf{I} - \bm{Q}) \rVert_{p\slash 2} \\
\label{eq:calc2}
    &\qquad\leq \biggl\lVert (\mathbf{I} - \bm{Q})
    \begin{bmatrix}
    \bm{0} & ~ \\
    ~ & \bm{A}_2^2
    \end{bmatrix} 
    (\mathbf{I} - \bm{Q}) \biggr\rVert_{p\slash 2}
    + \biggl\lVert (\mathbf{I} - \bm{Q}) 
    \begin{bmatrix}
    \bm{A}_1^2 & ~ \\
    ~ & \bm{0}
    \end{bmatrix}
    (\mathbf{I} - \bm{Q}) \biggr\rVert_{p\slash 2} \\
\label{eq:calc3}
    &\qquad\leq 
    \biggl\lVert 
    \begin{bmatrix}
    \bm{0} & ~ \\
    ~ & \bm{A}_2
    \end{bmatrix} \biggr\rVert_p^2
    + \biggl\lVert (\mathbf{I} - \bm{Q}) 
    \begin{bmatrix}
    \bm{A}_1 & ~ \\
    ~ & \bm{0}
    \end{bmatrix} \biggr\rVert_p^2 \\
\label{eq:calc4}
    &\qquad=
    \lVert \bm{A}_2 \rVert_p^2
    + \biggl\lVert 
    \begin{bmatrix}
    \bm{A}_1 & ~ \\
    ~ & \bm{0}
    \end{bmatrix} (\mathbf{I} - \bm{Q}) 
    \begin{bmatrix}
    \bm{A}_1 & ~ \\
    ~ & \bm{0}
    \end{bmatrix} \biggr\rVert_{p/2}.
\end{align}
\Cref{eq:calc1,eq:calc3,eq:calc4} use identities
$\lVert \bm{M} \rVert_p^2 = \lVert \bm{M}^\ast \bm{M} \rVert_{p/2} = \lVert \bm{MM}^\ast \rVert_{p/2}$,
while
\cref{eq:calc2} follows from the triangle inequality, and
\cref{eq:calc3} relies on the property that the orthogonal projection $\mathbf{I} - \bm{Q}$ is a $p$-norm contraction.

Next, the orthogonal projection $\bm{Q}$ is given explicitly by the formula
\begin{equation}
\label{eq:use1}
    \bm{Q} = 
    (\bm{A} \bm{\Omega} \bm{\Omega}_1^{\dagger})
    \bigl[(\bm{A} \bm{\Omega} \bm{\Omega}_1^{\dagger})^{\ast} 
    (\bm{A} \bm{\Omega} \bm{\Omega}_1^{\dagger})\bigr]^{\dagger}
    (\bm{A} \bm{\Omega} \bm{\Omega}_1^{\dagger})^{\ast}.
\end{equation}
By the assumption that $\textup{rank}[\bm{\Omega}_1] = r$,
we have $\bm{\Omega}_1 \bm{\Omega}_1^{\dagger} = \mathbf{I}$.
Consequently,
\begin{equation}
\label{eq:use2}
    \bm{A} \bm{\Omega} \bm{\Omega}_1^{\dagger}
    = \begin{bmatrix}
    \bm{A}_1 \\ 
    \bm{A}_2 \bm{\Omega}_2 \bm{\Omega}_1^{\dagger}
    \end{bmatrix}
    =: \begin{bmatrix}
    \bm{A}_1 \\ 
    \bm{F}
    \end{bmatrix}.
\end{equation}
We have introduced the abbreviation $\bm{F} = \bm{A}_2 \bm{\Omega}_2 \bm{\Omega}_1^{\dagger}$.
A direct calculation involving \cref{eq:use1} and \cref{eq:use2} delivers an expression
involving the parallel sum:
\begin{equation*}
    \begin{bmatrix}
    \bm{A}_1 & ~ \\
    ~ & \bm{0}
    \end{bmatrix} (\mathbf{I} - \bm{Q}) 
    \begin{bmatrix}
    \bm{A}_1 & ~ \\
    ~ & \bm{0}
    \end{bmatrix}
    = \begin{bmatrix}
    \bm{A}_1^2 : (\bm{F}^\ast \bm{F})
    & ~ \\
    ~ & \bm{0}
    \end{bmatrix}.
\end{equation*}
By \cref{lem:parallel_sum} parts (3) and (4), we obtain the bound
\begin{equation*}
    \bigl\lVert \bm{A}_1^2: (\bm{F}^\ast \bm{F})  \bigr\rVert_{p/2}
    \leq \bigl\lVert \bm{A}_1^2 \bigr\rVert_{p/2} : \bigl\lVert \bm{F}^\ast \bm{F} \bigr\rVert_{p/2}
    = \lVert \bm{A}_1 \rVert_p^2 : \lVert \bm{F} \rVert_p^2
    \leq \lVert \bm{F} \rVert_p^2.
\end{equation*}
This statement leads to \cref{eq:satisfies} and completes the proof.
\end{proof}

The error bounds in \cref{prop:ext_analysis} ensure that \texttt{RSVD} gives a high-quality approximation when the quantity $\lVert \bm{A}_2 \bm{\Omega}_2 \bm{\Omega}_1^{\dagger} \rVert_p^2$ is small.
To understand when this quantity is likely to be small, we need additional tools from random matrix theory, to be presented in the next section.

\subsection{Calculations using random matrix theory} \label{sec:random}

To complete the proof of our main result, we apply the following error bounds for Gaussian matrices, which are derived in \cref{app:convenient}.

\begin{proposition}[Gaussian matrix products] \label{prop:convenient}
Fix a matrix $\bm{S} \in \mathbb{R}^{(N - r) \times (N - r)}$ and consider independent Gaussian matrices $\bm{G} \in \mathbb{R}^{(N - r) \times k}$ and $\bm{H} \in \mathbb{R}^{r \times k}$, where $k \geq r$.
Then, for any $u, t \geq 0$,
\begin{equation}
\label{eq:probability_bound}
    \mathbb{P} \left\{\lVert \bm{S} \bm{G} \bm{H}^{\dagger} \rVert^2_{\textup{F}}
    > \frac{utr}{k-r+1} \left\lVert \bm{S} \right\rVert_{\textup{F}}^2\right\}
    \leq \mathrm{e}^{-(u - 2) \slash 4} + \sqrt{\pi r} \biggl(\frac{t}{\mathrm{e}}\biggr)^{-(k-r+1) \slash 2}.
\end{equation}
Additionally, if $k \geq r + 2$,
\begin{equation}
\label{eq:expectation_bound}
    \Expect \lVert \bm{S} \bm{G} \bm{H}^{\dagger} \rVert^2_{\rm F}
    = \frac{r}{k-r-1} \lVert \bm{S} \rVert^2_{\rm F}.
\end{equation}
\end{proposition}

The expectation formula \eqref{eq:expectation_bound} was previously presented in \cite[Thm.~10.5]{halko2011finding}.
It follows from the expectation formula for an inverse Wishart matrix \cite[pg.~119]{press2005applied}.
The probability bound \cref{eq:probability_bound} improves on previous results \cite[Thm.~10.7]{halko2011finding}.
It is derived using a Lipschitz concentration inequality for Gaussian random fields, together with a new concentration inequality for the trace of an inverse Wishart matrix.

The random matrix theory allow us to complete the proof of our main results, \cref{thm:main_thm,cor:main_cor},
as follows.

\begin{proof}[Proof of \cref{thm:main_thm} and \cref{cor:main_cor}]
Without loss of generality, assume that $\bm{A} \in \mathbb{R}^{N \times N}$ is diagonal and psd, with nonincreasing diagonal entries.
We partition the Gaussian test matrix
$\bm{\Omega} = \bigl[\begin{smallmatrix}
\bm{\Omega}_1 \\ \bm{\Omega}_2
\end{smallmatrix}\bigr]$, where $\bm{\Omega}_1$ is a $r \times k$ matrix and $\bm{\Omega}_2$ is a $(N - r) \times k$ matrix.  The two submatrices $\bm{\Omega}_1$ and $\bm{\Omega}_2$ are independent Gaussian matrices, and $\bm{\Omega}_1$ almost surely has rank $r$.

\texttt{RSVD} produces an approximation $\hat{\bm{A}} \in \mathbb{R}^{N \times n}$ that satisfies
\begin{equation*}
    \lVert \bm{A} - \hat{\bm{A}} \rVert_p^2
    = \lVert (\mathbf{I} - \bm{\Pi}_{\bm{A} \bm{\Omega}}) \bm{A} \rVert_p^2
    \leq \lVert (\mathbf{I} - \bm{\Pi}_{\bm{A} \bm{\Omega} \bm{\Omega}_1^\dagger}) \bm{A} \rVert_p^2
\end{equation*}
because $\bm{A} \bm{\Omega}$ has a bigger range than $\bm{A} \bm{\Omega} \bm{\Omega}_1^\dagger$; see \cref{lem:bigger_projection}.
Next, \cref{prop:ext_analysis} shows that, almost surely, 
\begin{equation*}
    \lVert (\mathbf{I} - \bm{\Pi}_{\bm{A} \bm{\Omega} \bm{\Omega}_1^\dagger}) \bm{A} \rVert_p^2
    \leq \lVert \bm{A}_2 \rVert_p^2 
    + \lVert \bm{A}_2 \bm{\Omega}_2 \bm{\Omega}_1^{\dagger} \rVert_p^2
    \leq \lVert \bm{A}_2 \rVert_p^2 
    + \lVert \bm{A}_2 \bm{\Omega}_2 \bm{\Omega}_1^{\dagger} \rVert_{\rm F}^2.
\end{equation*}
Since $\bm{\Omega}_1$ and $\bm{\Omega}_2$ are independent Gaussian matrices, the error bounds from \cref{prop:convenient} deliver the \texttt{RSVD} error bounds \cref{eq:aesthetic,eq:simple}.

As for the corollary, \texttt{NysSVD} produces an
approximation $\hat{\bm{A}} \in \mathbb{R}^{N \times N}$ that satisfies
\begin{equation*}
    \lVert \bm{A} - \hat{\bm{A}} \rVert_p
    = \lVert \bm{A}^{1 \slash 2}(\mathbf{I} - \bm{\Pi}_{\bm{A}^{1 \slash 2} \bm{\Omega}}) \bm{A}^{1 \slash 2} \rVert_p
    = \lVert (\mathbf{I} - \bm{\Pi}_{\bm{A}^{1 \slash 2} \bm{\Omega}}) \bm{A}^{1 \slash 2} \rVert_{2p}^2.
\end{equation*}
We recognize $\lVert (\mathbf{I} - \bm{\Pi}_{\bm{A}^{1 \slash 2} \bm{\Omega}}) \bm{A}^{1 \slash 2} \rVert_{2p}^2$ as the square error of the \texttt{RSVD}, applied to the matrix $\bm{A}^{1 \slash 2}$.
By invoking the \texttt{RSVD} error bounds for $\lVert (\mathbf{I} - \bm{\Pi}_{\bm{A}^{1 \slash 2} \bm{\Omega}}) \bm{A}^{1 \slash 2} \rVert_{2p}^2$, we confirm the \texttt{NysSVD} error bounds \cref{eq:aesthetic2,eq:simple2}.
\end{proof}

\subsection{Additional expectation bounds} \label{sec:discussion}

The analysis of \texttt{RSVD} using parallel sums of psd matrices leads to a family of expectation bounds with respect to various Schatten $p$-norms.
We give a collection of these bounds in \cref{thm:more_bounds}, with the proof appearing in \cref{app:more_bounds}.
The bounds underscore the potential for high errors when $k - r$ is small.

\begin{theorem}[\texttt{RSVD} expectation bounds] \label{thm:more_bounds}
Fix a matrix $\bm{A} \in \mathbb{R}^{L \times N}$, a target rank $r \geq 1$,
and a block size $k \geq r$.  Then \texttt{RSVD} with a block of $k$ Gaussian starting vectors
generates a random approximation $\hat{\bm{A}} \in \mathbb{R}^{L \times N}$ that satisfies
\begin{equation*}
    \frac{\Expect \lVert \bm{A} - \hat{\bm{A}} \rVert_{\rm F}^2}
    {\lVert \bm{A} - \lfloor \bm{A} \rfloor_r \rVert_{\rm F}^2}
    \leq 
    \begin{cases}
    1 + \dfrac{r}{k-r-1},
    & k \geq r + 2; \\[10pt]
    1
    + r \log\biggl(1 + 
    \dfrac{\lVert \lfloor \bm{A} \rfloor_r \rVert_{\rm F}^2}
    {\lVert \bm{A} - \lfloor \bm{A} \rfloor_r \rVert_{\rm F}^2}\biggr),
    & k = r + 1; \\[10pt]
    1
    + r \biggl(
    \dfrac{{\rm \pi} \lVert \lfloor \bm{A} \rfloor_r \rVert_{\rm F}^2}
    {2 \lVert \bm{A} - \lfloor \bm{A} \rfloor_r \rVert_{\rm F}^2} \biggr)^{1 \slash 2},
    & k = r.
    \end{cases}
\end{equation*}
Additionally, if the block size $k \geq r + 2$,
\begin{equation*}
    \Expect \lVert \bm{A} - \hat{\bm{A}} \rVert^2 
    \leq \biggl(1 + \frac{2r}{k - r - 1}\biggr) \biggl(\lVert \bm{A} - \lfloor \bm{A} \rfloor_r \rVert^2
    + \frac{{\rm e}^2}{k-r} \lVert \bm{A} - \lfloor \bm{A} \rfloor_r \rVert^2_{\rm F} \biggr).
\end{equation*}
Last, if the block size $k \geq r + 4$,
\begin{align*}
    & \bigl(\Expect \lVert \bm{A} - \hat{\bm{A}} \rVert_4^4 \bigr)^{1 \slash 2}
    \leq \biggl(1 + \frac{r + 1}{k-r-3}\biggr) \biggl(\lVert \bm{A} - \lfloor \bm{A} \rfloor_r \rVert_4^2 + \frac{1}{\sqrt{k - r}} \lVert \bm{A} - \lfloor \bm{A} \rfloor_r \rVert_{\rm F}^2 \biggr); \\
    & \bigl(\Expect \lVert \bm{A} - \hat{\bm{A}} \rVert^4 \bigr)^{1 \slash 2} 
    \leq \biggl(1 + \frac{2(r+1)}{k - r - 3}\biggr) \biggl(\lVert \bm{A} - \lfloor \bm{A} \rfloor_r \rVert^2
    + \frac{\sqrt{3}\, {\rm e}^2}{k - r} \lVert \bm{A} - \lfloor \bm{A} \rfloor_r \rVert_{\rm F}^2 \biggr).
\end{align*}
\end{theorem}

Once we have derived expectation bounds for \texttt{RSVD}, we obtain matching bounds for \texttt{NysSVD} due to the close relationship between the two algorithms:

\begin{corollary}[\texttt{NysSVD} expectation bounds] \label{cor:nystrom}
Fix a psd matrix $\bm{A} \in \mathbb{R}^{N \times N}$ and a target rank $r \geq 1$.
Then \texttt{NysSVD} with $k \geq r + 2$ Gaussian initialization vectors produces a random approximation $\hat{\bm{A}} \in \mathbb{R}^{N \times N}$ that satisfies
\begin{equation*}
    \frac{\Expect \lVert \bm{A} - \hat{\bm{A}} \rVert_{\ast}}
    {\lVert \bm{A} - \lfloor \bm{A} \rfloor_r \rVert_{\ast}}
    \leq 
    \begin{cases}
    1 + \dfrac{r}{k-r-1},
    & k \geq r + 2; \\[10pt]
    1
    + r \log\biggl(1 + 
    \dfrac{\lVert \lfloor \bm{A} \rfloor_r \rVert_{\ast}}
    {\lVert \bm{A} - \lfloor \bm{A} \rfloor_r \rVert_{\ast}}\biggr),
    & k = r + 1; \\[10pt]
    1
    + r \biggl(
    \dfrac{{\rm \pi} \lVert \lfloor \bm{A} \rfloor_r \rVert_{\ast}}
    {2 \lVert \bm{A} - \lfloor \bm{A} \rfloor_r \rVert_{\ast}} \biggr)^{1 \slash 2},
    & k = r.
    \end{cases}
\end{equation*}
Additionally, if the block size $k \geq r + 2$,
\begin{equation*}
    \Expect \lVert \bm{A} - \hat{\bm{A}} \rVert
    \leq \biggl(1 + \frac{2r}{k - r - 1}\biggr) \biggl(\lVert \bm{A} - \lfloor \bm{A} \rfloor_r \rVert
    + \frac{{\rm e}^2}{k-r} \lVert \bm{A} - \lfloor \bm{A} \rfloor_r \rVert_{\ast} \biggr).
\end{equation*}
Last, if the block size $k \geq r + 4$,
\begin{align*}
    & \bigl(\Expect \lVert \bm{A} - \hat{\bm{A}} \rVert_{\rm F}^2 \bigr)^{1 \slash 2}
    \leq \biggl(1 + \frac{r + 1}{k-r-3}\biggr) \biggl(\lVert \bm{A} - \lfloor \bm{A} \rfloor_r \rVert_{\rm F} + \frac{1}{\sqrt{k - r}} \lVert \bm{A} - \lfloor \bm{A} \rfloor_r \rVert_{\ast} \biggr); \\
    & \bigl(\Expect \lVert \bm{A} - \hat{\bm{A}} \rVert^2 \bigr)^{1 \slash 2} 
    \leq \biggl(1 + \frac{2(r+1)}{k - r - 3}\biggr) \biggl(\lVert \bm{A} - \lfloor \bm{A} \rfloor_r \rVert
    + \frac{\sqrt{3} {\rm e}^2}{k - r} \lVert \bm{A} - \lfloor \bm{A} \rfloor_r \rVert_{\ast} \biggr).
\end{align*}
\end{corollary}
\begin{proof}
\texttt{NysSVD} produces an approximation $\bm{A}\langle \bm{\Omega} \rangle$ satisfying
\begin{equation*}
    \lVert \bm{A} - \bm{A}\langle \bm{\Omega} \rangle \rVert_p
    = \lVert \bm{A}^{1 \slash 2}(\mathbf{I} - \bm{\Pi}_{\bm{A}^{1 \slash 2} \bm{\Omega}}) \bm{A}^{1 \slash 2} \rVert_p
    = \lVert (\mathbf{I} - \bm{\Pi}_{\bm{A}^{1 \slash 2} \bm{\Omega}}) \bm{A}^{1 \slash 2} \rVert_{2p}^2,
\end{equation*}
where $\bm{\Omega} \in \mathbb{R}^{N \times k}$ is a random Gaussian matrix.
We recognize
$\lVert (\mathbf{I} - \bm{\Pi}_{\bm{A}^{1 \slash 2} \bm{\Omega}}) \bm{A}^{1 \slash 2} \rVert_{2p}^2$ as the square error of \texttt{RSVD} applied to the matrix $\bm{A}^{1 \slash 2}$ and apply the error bounds
given in 
\cref{thm:more_bounds}.
\end{proof}

The error bounds in \cref{thm:more_bounds,cor:nystrom} can be used in the analysis of related randomized algorithms.
For example, the recent paper \cite{epperly2023xtrace} uses \cref{thm:more_bounds,cor:nystrom} to control the variance of several randomized trace estimators.

\section{Analysis of Krylov algorithms} \label{sec:theory2}

In this section, we analyze \texttt{RBKI} and \texttt{NysBKI} and provide a matching analysis of \texttt{RSI} and \texttt{NysSI} for comparison.
We develop both \emph{gapless} error bounds that do not require any distance between the singular values and \emph{gapped} error bounds that grow increasingly strong with the size of singular value gaps.

The gapless error bounds are given in \cref{thm:gapless},
and the proof appears in \cref{sec:spectral}.
Instead of bounding the error ratio
\begin{equation*}
    \frac{\Expect \lVert \bm{A} - \hat{\bm{A}} \rVert^2}{\sigma_{r+1}(\bm{A})^2},
\end{equation*}
which is ideally close to one, these bounds apply to the log-error ratio
\begin{equation*}
    \log\biggl(\frac{\Expect \lVert \bm{A} - \hat{\bm{A}} \rVert^2}{\sigma_{r+1}(\bm{A})^2}\biggr),
\end{equation*}
which is ideally close to zero.
We demonstrate that the log-error ratio decreases at a rate $\mathcal{O}(1/m)$ for \texttt{RSI} or \texttt{NysSI} and a faster rate $\mathcal{O}(1/m^2)$ for \texttt{RBKI} or \texttt{NysBKI}.  In each case, $m$ is the total number of multiplications with $\bm{A}$ or $\bm{A}^\ast$.

\begin{theorem}[Gapless error bounds] \label{thm:gapless}
Fix a matrix $\bm{A} \in \mathbb{R}^{L \times N}$ and a target rank $r \geq 1$ such that $\sigma_{r+1}(\bm{A}) > 0$. Then \texttt{RSI} with $k \geq r + 2$ Gaussian initialization vectors and $m$ matrix multiplications generates a random approximation $\hat{\bm{A}} \in \mathbb{R}^{L \times N}$ that satisfies
\begin{equation*}
\tag{\texttt{RSI}}
    \log\biggl(\frac{\Expect \lVert \bm{A} - \hat{\bm{A}} \rVert^2}{\sigma_{r+1}(\bm{A})^2}\biggr)
    \leq \frac{1}{m - 1} \log\biggl(1 + \frac{r}{k - r - 1}
    \sum_{i > r} \frac{\sigma_i(\bm{A})^{2}}
    {\sigma_{r+1}(\bm{A})^{2}}
    \biggr).
\end{equation*}
With the same initialization and the same number of multiplications, \texttt{RBKI} satisfies
\begin{equation*}
\tag{\texttt{RBKI}}
    \log\biggl(\frac{\Expect \lVert \bm{A} - \hat{\bm{A}} \rVert^2}{\sigma_{r+1}(\bm{A})^2}\biggr)
    \leq \frac{1}{4(m - 2)^2} \biggl[\log\biggl(4 + \frac{4r}{k - r - 1}
    \sum_{i > r} \frac{\sigma_i(\bm{A})^2}{\sigma_{r+1}(\bm{A})^2}
    \biggr)\biggr]^2.
\end{equation*}
If $\bm{A} \in \mathbb{R}^{N \times N}$ is a psd matrix, \texttt{NysSI} satisfies
\begin{equation*}
\tag{\texttt{NysSI}}
    \log\biggl(\frac{\Expect \lVert \bm{A} - \hat{\bm{A}} \rVert^2}{\sigma_{r+1}(\bm{A})^2}\biggr)
    \leq \frac{1}{m - \frac{1}{2}} \log\biggl(1 + \frac{r}{k - r - 1}
    \sum_{i > r} \frac{\sigma_i(\bm{A})^{2}}
    {\sigma_{r+1}(\bm{A})^{2}}
    \biggr),
\end{equation*}
and \texttt{NysBKI} satisfies
\begin{equation*}
\tag{\texttt{NysBKI}}
    \log\biggl(\frac{\Expect \lVert \bm{A} - \hat{\bm{A}} \rVert^2}{\sigma_{r+1}(\bm{A})^2}\biggr)
    \leq \frac{1}{8(m - \frac{3}{2})^2} \biggl[\log\biggl(4 + \frac{4r}{k - r - 1}
    \sum_{i > r} \frac{\sigma_i(\bm{A})^2}{\sigma_{r+1}(\bm{A})^2}
    \biggr)\biggr]^2.
\end{equation*}
\end{theorem}

The gapless error bounds lead to algorithms insights.
For example, \texttt{RBKI} and \texttt{NysBKI} converge at a fundamentally faster rate than the other algorithms.
However, these bounds only describe the worst-case performance, which applies when the gaps between singular values are negligibly small.
In more realistic applications, including all the test problems in \cref{sec:experiments}, we find that noticeable gaps appear.
Therefore, we develop a set of gapped error bounds in \cref{thm:gapped} that guarantee faster, exponential convergence rates.  The proof appears in \cref{sec:frobenius}.

\begin{theorem}[Gapped error bounds] \label{thm:gapped}
Fix a matrix $\bm{A} \in \mathbb{R}^{L \times N}$, a Schatten $p$-norm with $2 \leq p \leq \infty$, and two singular values $\sigma_r(\bm{A})$ and $\sigma_s(\bm{A})$ with $r < s$.
Define the singular value gap
\begin{equation*}
    \gamma = \frac{\sigma_r(\bm{A}) - \sigma_s(\bm{A})}{\sigma_r(\bm{A}) + \sigma_s(\bm{A})} \in [0, 1],
\end{equation*}
and let $\lfloor \bm{A} \rfloor_r$ and $\lfloor \bm{A} \rfloor_{s-1}$ denote optimal rank-$r$ and rank-$(s-1)$ approximations of $\bm{A}$.
Then, \texttt{RSI} with $k \geq s + 1$ Gaussian initialization vectors and $m$ matrix multiplications generates a random approximation $\hat{\bm{A}} \in \mathbb{R}^{L \times N}$ that satisfies
\begin{equation*}
\tag{\texttt{RSI}}
    \Expect \lVert \bm{A} - \hat{\bm{A}} \rVert_p^2 \leq \lVert \bm{A} - \lfloor \bm{A} \rfloor_r \rVert_p^2
    + \mathrm{e}^{-4(m - 2) \gamma} \cdot \frac{s - 1}{k - s} \lVert \bm{A} - \lfloor \bm{A}\rfloor_{s-1} \rVert_{\rm F}^2.
\end{equation*}
With the same initialization and the same number of multiplications, \texttt{RBKI} satisfies
\begin{equation*}
\tag{\texttt{RBKI}}
    \Expect \lVert \bm{A} - \hat{\bm{A}} \rVert_p^2 \leq \lVert \bm{A} - \lfloor \bm{A} \rfloor_r \rVert_p^2
    + 4\,\mathrm{e}^{-4(m - 2) \sqrt{\gamma}} \cdot \frac{s - 1}{k - s} \lVert \bm{A} - \lfloor \bm{A} \rfloor_{s-1} \rVert_{\rm F}^2.
\end{equation*}
If $\bm{A} \in \mathbb{R}^{N \times N}$ is a psd matrix and $1 \leq p \leq \infty$, \texttt{NysSI} satisfies
\begin{equation*}
\tag{\texttt{NysSI}}
    \Expect \lVert \bm{A} - \hat{\bm{A}} \rVert_p^2 \leq \lVert \bm{A} - \lfloor \bm{A} \rfloor_r \rVert_p^2
    + \mathrm{e}^{-4(m - 3/2) \gamma} \cdot \frac{s - 1}{k - s} \lVert \bm{A} - \lfloor \bm{A}\rfloor_{s-1} \rVert_{\rm F}^2.
\end{equation*}
and \texttt{NysBKI} satisfies
\begin{equation*}
\tag{\texttt{NysBKI}}
    \Expect \lVert \bm{A} - \hat{\bm{A}} \rVert_p^2 \leq \lVert \bm{A} - \lfloor \bm{A} \rfloor_r \rVert_p^2
    + 4\, \mathrm{e}^{-4(m - 3/2) \sqrt{2\gamma}} \cdot \frac{s - 1}{k - s} \lVert \bm{A} - \lfloor \bm{A} \rfloor_{s-1} \rVert_{\rm F}^2.
\end{equation*}
\end{theorem}

The gapped error bounds show how Krylov methods capitalize upon small singular value gaps $\gamma \ll 1$ to achieve accelerated convergence.
For example, we can ensure an error bound
\begin{equation*}
    \Expect \lVert \bm{A} - \hat{\bm{A}} \rVert_{\rm F}^2 \leq (1 + \varepsilon) \cdot \lVert \bm{A} - \lfloor \bm{A} \rfloor_r \rVert_{\rm F}^2,
\end{equation*}
if we set $m = \mathcal{O}(\log (1 \slash \varepsilon) \slash \gamma)$ in \texttt{RSI} or $m = \mathcal{O}(\log(1 \slash \varepsilon) \slash \sqrt{\gamma})$ in \texttt{RBKI}.
The change from a $1 / \gamma$-dependence to a $1 / \sqrt{\gamma}$-dependence leads to an accelerated convergence rate in many applications.

\Cref{thm:gapless,thm:gapped} improve upon existing results for \texttt{RSI} \cite{halko2011finding,gu2015subspace,bjarkason2019pass}, \texttt{RBKI} \cite{musco2015randomized,wang2015improved,drineas2018structural,tropp2018analysis,BCW22,bakshi2023krylov,meyer2023unreasonable}, \texttt{NysSI}, and \texttt{NysBKI}.
Musco and Musco \cite{musco2015randomized} previously established the $\mathcal{O}(m^{-2})$ gapless convergence rate for \texttt{RBKI}, and we improve this result by obtaining explicit constants.
Our gapped error bounds for \texttt{RBKI} and \texttt{RSI} also feature sharper constants than previous results \cite{gu2015subspace,bjarkason2019pass,drineas2018structural,tropp2018analysis}.
Last, the bounds for \texttt{NysSI} and \texttt{NysBKI} are completely new.

\begin{figure}[t]
    \centering
    \includegraphics[width=\textwidth]{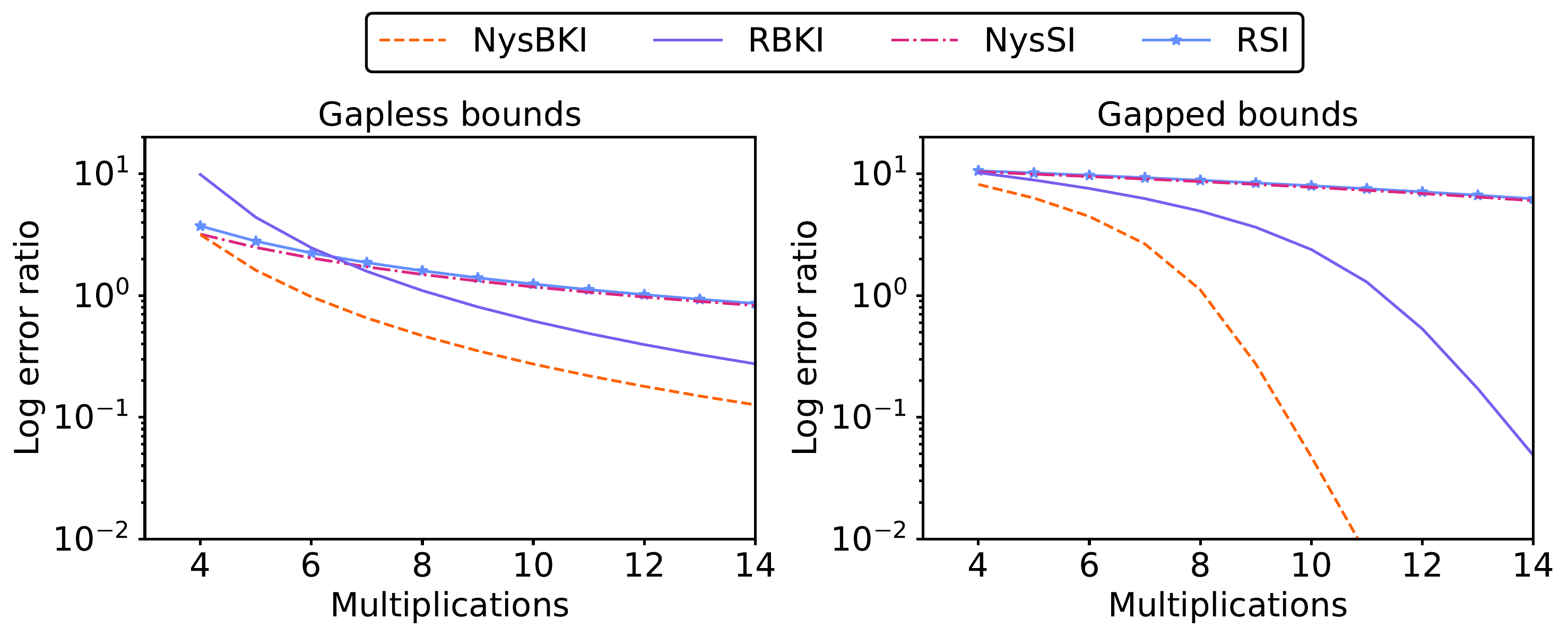}
    \caption{\textbf{(Gapless vs. gapped bounds).} Gapless (left) versus gapped (right) error bounds on the log error ratio $\log \bigl(\Expect \lVert \bm{B} - \hat{\bm{B}} \rVert^2 / \sigma_{r+1}(\bm{B})^2\bigr)$ when approximating the matrix $\bm{B}$ from \cref{sec:comparison} with block size $k = 100$ and target rank $r = 75$.}
    \label{fig:error_bounds}
\end{figure}

Because of the sharp, explicit nature of our estimates,
we are able to make quantitative predictions:
the gapless (\cref{thm:gapless}) and gapped (\cref{thm:gapped}) bounds predict that \texttt{NysSI} outperforms \texttt{RSI} by one-half of a matrix--matrix multiplication while \texttt{NysBKI} outperforms \texttt{RBKI} by a factor of $\sqrt{2}$ matrix--matrix multiplications.
These predictions are in close agreement with numerical experiments (\cref{fig:krylov_v_subspace,fig:eigenvectors}).
To illustrate the quantitative nature of our bounds, we evaluate these bounds for the matrix $\bm{B}$ as constructed in \cref{sec:comparison}
and present the output in \cref{fig:error_bounds}.

Our error bounds represent a step forward in the theoretical understanding of Krylov methods, but they are not the end of the story.
Recent work \cite{bakshi2020robust,bakshi2023krylov,meyer2023unreasonable} has begun to develop gapless error bounds for \texttt{RBKI} in the Frobenius norm, which complement the spectral norm approximation guarantees given here.
Additionally, our main error bounds in \cref{thm:gapless,thm:gapped} only apply when the block size $k$ exceeds the target rank $r$.
Although it is impossible to establish any \emph{gapless} error bounds in the case $k < r$,
Meyer and coauthors \cite{meyer2023unreasonable} (also see earlier work of \cite{yuan2018superlinear}) have established a \emph{gapped} error bound for \texttt{RBKI} for $k < r$ \cite[Thm.~1]{meyer2023unreasonable}.
Their result is more complicated than \cref{thm:gapless,thm:gapped} and lacks explicit constants, but it provides additional evidence of the computational advantages of block Krylov methods.

The rest of this section contains an overview of our technical approach (\cref{sec:overview}) and background material on Chebyshev polynomials (\cref{sec:chebyshev}), followed by a derivation of the the gapped error bounds (\cref{sec:frobenius}) and the gapless error bounds (\cref{sec:spectral}).

\subsection{Overview of technical approach} \label{sec:overview}

To prove the gapped and gapless error bounds, our main approach is to apply a polynomial filter $\phi(\bm{A})$ that increases the top singular values and decreases the bottom singular values.
As we choose a filter, we need to ensure that $\phi(\bm{A}) \bm{\Omega}$ lies inside the approximation space.
Therefore, we use power function filters to analyze subspace iteration methods, and we use Chebyshev polynomial filters to analyze Krylov methods.
The enhanced ability of Chebyshev polynomials to separate the top singular values from the bottom singular values is why we obtain such powerful Krylov bounds.  Chebyshev polynomials have long been used in this context; for example, see Kuczy{\'n}ski and Wo{\'z}niakowski~\cite{kuczynski1992estimating}.

Our filtering analysis relies on the fact that
\begin{equation}
\label{eq:to_relate}
    \lVert (\mathbf{I} - \bm{\Pi}_{\phi(\bm{A}) \bm{\Omega}}) \phi(\bm{A}) \rVert,
\end{equation}
is the error of \texttt{RSVD} applied to a matrix $\phi(\bm{A})$, and the \texttt{RSVD} error is bounded by \cref{thm:main_thm}.
The main difficulty is relating the \texttt{RSVD} error \cref{eq:to_relate} to the quantity
\begin{equation*}
    \lVert (\mathbf{I} - \bm{\Pi}_{\phi(\bm{A}) \bm{\Omega}}) \bm{A} \rVert,
\end{equation*}
the error actually attained by \texttt{RSI} and \texttt{RBKI} algorithms.

\begin{figure}[t]
    \centering
    \includegraphics[width = .5\textwidth]{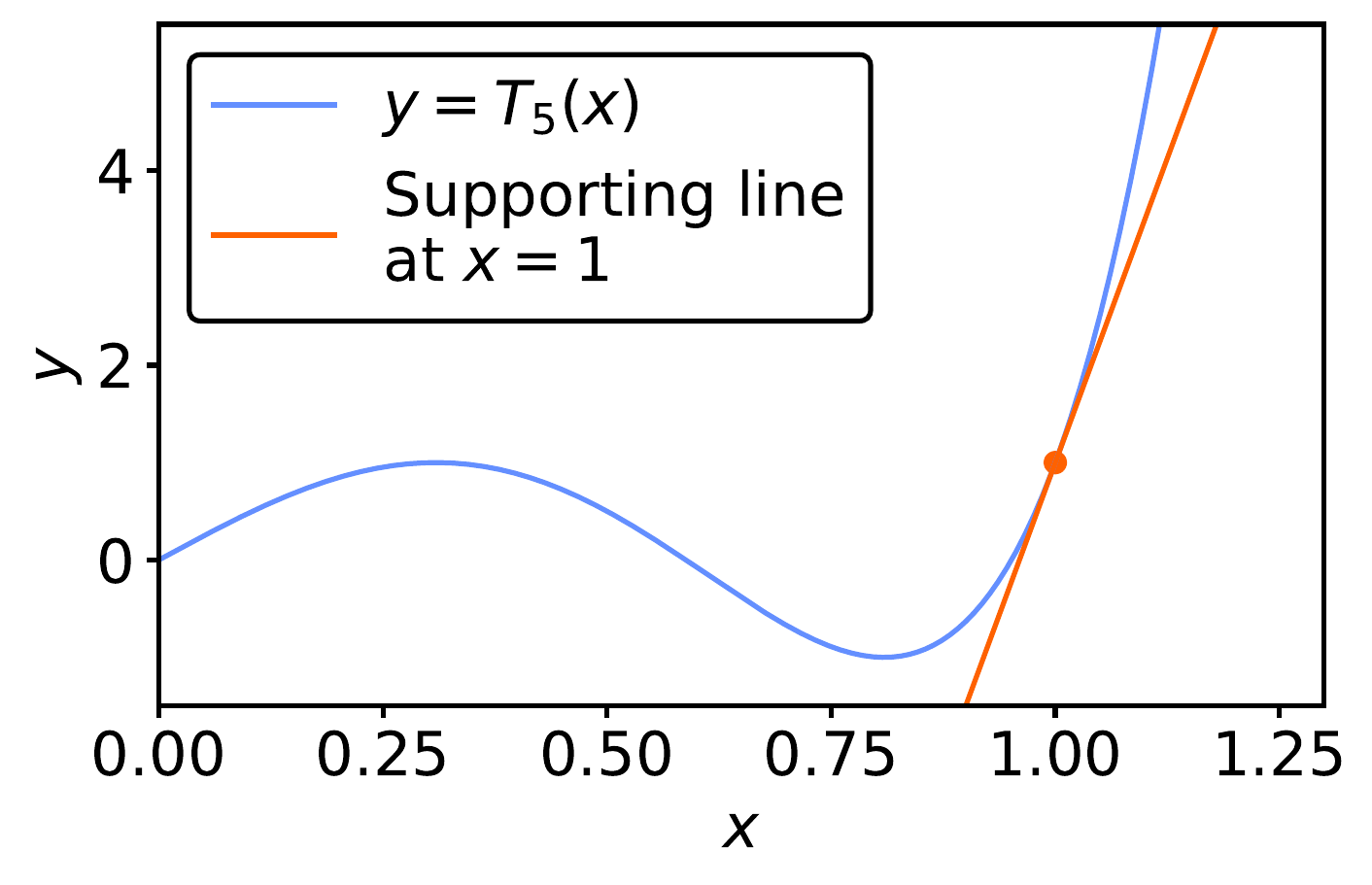}
    \caption{\textbf{(Chebyshev polynomials are almost convex).} The Chebyshev polynomial $T_5(x)$ admits supporting lines on the range $x \geq 1.$}
    \label{fig:quasi_convex}
\end{figure}

We present two mathematical arguments that relate $\lVert (\mathbf{I} - \bm{\Pi}_{\phi(\bm{A}) \bm{\Omega}}) \phi(\bm{A}) \rVert$ with $\lVert (\mathbf{I} - \bm{\Pi}_{\phi(\bm{A}) \bm{\Omega}}) \bm{A} \rVert$.
First, we use convexity to establish \cref{lem:quasi_convex}, which allows us to prove the gapless bounds.
For \texttt{RSI}, the convexity argument is already visible in~\cite[Prop.~8.6]{halko2011finding}.  For \texttt{RBKI}, the crucial new observation 
is that Chebyshev polynomials are not globally convex, yet they admit supporting lines on the range $x \geq 1$. The supporting lines are all that we need to extend the convexity-based arguments.
See \cref{fig:quasi_convex} for an illustration.
Second, as a more standard mathematical result, we use matrix submultiplicativity to establish \cref{lem:enhancing}, which allows us to prove the gapped bounds.

The gapless and gapped error bounds apply to the \emph{untruncated} approximation $\hat{\bm{A}}$, which is the approximation returned by our pseudocode and recommended for most practical applications.
However, similar bounds extend to the approximation $\lfloor \hat{\bm{A}} \rfloor_r$, which comes from applying an $r$-truncated singular value decomposition to $\hat{\bm{A}}$.
The gapped error bounds (\cref{thm:gapped}) remain valid even when we replace $\hat{\bm{A}}$ with $\lfloor \hat{\bm{A}} \rfloor_r$, due to the reverse Eckart--Young inequality \cite[Thm.~3.4]{gu2015subspace}.
Likewise, we can prove gapless error bounds for $\lfloor \hat{\bm{A}} \rfloor_r$ by pursuing Musco and Musco's observation \cite{musco2015randomized} that either there is a singular value gap and the gapped error bounds apply or else there is no gap and the approximation is close to optimal. We omit a detailed treatment.

\subsection{Chebyshev polynomials} \label{sec:chebyshev}

The $q$th Chebyshev polynomial $x \mapsto T_q(x)$ is a polynomial of degree $q$ that shares some properties with the $q$th power function $x \mapsto x^q$.
The values of both polynomials lie inside $[-1, 1]$ for $0 \leq x \leq 1$,
and both polynomials exhibit a sharp rate of increase for $x \geq 1$.
\Cref{lem:chebyshev} provides more fine-grained estimates that will be used in the analysis.

\begin{lemma}[Chebyshev polynomials] \label{lem:chebyshev}
For $q \in \{0, 1, 2, \ldots\}$ and $x \geq 1$, define the Chebyshev polynomial $T_q$ of the first kind by
\begin{equation}
    T_q(x) = 
    \begin{cases}
    \cos(q \arccos x), & 0 \leq x \leq 1, \\
    \cosh (q \arcosh x), & x \geq 1.
    \end{cases}
\end{equation}
Then, the following properties hold:
\begin{enumerate}
\item $T_q$ is a polynomial of degree $q$ with the same parity as $q$.
If $q$ is odd, then $T_q$ is odd.  If $q$ is even, then $T_q$ is even.
\item $T_q(x)$ is bounded above by
\begin{equation*}
|T_q(x)| \leq 1,
\qquad 0 \leq x \leq 1.
\end{equation*}
\item $T_q(x)$ is increasing on the interval $x \geq 1$, and it is bounded below as
\begin{equation*}
    T_q\biggl(\frac{1 + r}{1 - r}\biggr) \geq \frac{1}{2} \mathrm{e}^{2q\sqrt{r}}, \qquad 0 \leq r < 1.
\end{equation*}
In contrast, the quantity $x^q$ exhibits a slower rate of increase:
\begin{equation*}
    \biggl(\frac{1 + r}{1 - r}\biggr)^q \geq \mathrm{e}^{2 q r}, \qquad 0 \leq r < 1.
\end{equation*}
\item
The functional inverse of the Chebyshev polynomial satisfies
\begin{equation*}
    T_q^{-1}(x) \leq \exp\biggl(\frac{1}{2} \biggl[\frac{\log(2x)}{q}\biggl]^2\biggl), \qquad x \geq 1.
\end{equation*}
\item
The function $\phi(x) = x T_q(x^{1 \slash 2})^2$ has a supporting line at every point $(x, \phi(x))$ with $x \geq 1$, as does the function $\phi(x) = x T_q(x^{1 \slash 4})^2$.
\end{enumerate}
\end{lemma}
\begin{proof}
Parts 1 and 2 are standard facts presented in \cite[Ch.~2]{mason2002chebyshev}. Parts 3--5 are technical results proved in \cref{app:chebyshev}.
\end{proof}

\subsection{Gapless error bounds} \label{sec:spectral}

To prove the gapless error bounds, we will need the following version of Jensen's inequality.
Compare this result with~\cite[Prop.~8.6]{halko2011finding}.

\begin{lemma}[Jensen's inequality with ``almost'' convex functions] \label{lem:quasi_convex}
    Consider a diagonal, psd matrix $\bm{A} \in \mathbb{R}^{N \times N}$, a random rank-$r$ orthogonal projection $\bm{Q} \in \mathbb{R}^{N \times N}$, and a function $f:[0, \infty) \rightarrow \mathbb{R}$ that has a supporting line at every point $(x, f(x))$ for $x \geq \sigma_{r+1}(\bm{A})^2$. Then
    \begin{equation}
    \label{eq:supporting_line}
        f(\Expect \lVert (\mathbf{I} - \bm{Q}) \bm{A} \rVert^2) \leq \Expect \lVert (\mathbf{I} - \bm{Q}) f(\bm{A}^2) (\mathbf{I} - \bm{Q}) \rVert.
    \end{equation}
\end{lemma}

\begin{proof}
We will repeatedly use the fact that $\bm{Q} \bm{A}$ is a rank-$r$ approximation of $\bm{A}$ and hence $\lVert \bm{A} - \bm{Q} \bm{A} \rVert \geq \sigma_{r+1}(\bm{A})$ by Weyl's inequality \cite[Thm.~4.3.1]{horn2012matrix}.
First, we set $x_0 = \Expect \lVert \bm{A} - \bm{Q} \bm{A} \rVert^2 \geq \sigma_{r+1}(\bm{A})^2$ and use the supporting line property to argue
\begin{equation*}
    f(x_0)
    = \Expect[f(x_0)
    + f^\prime(x_0) (\lVert \bm{A} - \bm{Q} \bm{A} \rVert^2 - x_0)]
    \leq \Expect f(\lVert \bm{A} - \bm{Q} \bm{A} \rVert^2)
\end{equation*}
Next, we define the random vector $\bm{v}$ to be a dominant eigenvector of $(\mathbf{I} - \bm{Q}) \bm{A}^2 (\mathbf{I} - \bm{Q})$, lying in the range of 
$\mathbf{I} - \bm{Q}$ and
normalized so that $\lVert \bm{v} \rVert = 1$.
We set $X = \lVert (\mathbf{I} - \bm{Q}) \bm{A} \rVert^2$ and observe that
\begin{equation}
\label{eq:averaging}
    X = \bm{v}^\ast (\mathbf{I} - \bm{Q}) \bm{A}^2 (\mathbf{I} - \bm{Q}) \bm{v}^\ast
    = \bm{v}^\ast \bm{A}^2 \bm{v}^\ast
    = \sum\nolimits_{i=1}^N a_{ii}^2 v_i^2,
\end{equation} 
where $a_{ii}$ denotes the $i$th diagonal entry of $\bm{A}$ (recall that $\bm{A}$ is a diagonal matrix) and $v_i$ is the $i$th entry of $\bm{v}$.
Last, we use \cref{eq:averaging} and the supporting line property to argue
\begin{align*}
    f(X) &= \sum\nolimits_{i=1}^N v_i^2 \bigl[f(X)
    + f^\prime(X) (a_{ii}^2 - X)\bigr] \\
    &\leq \sum\nolimits_{i=1}^N v_i^2 f(a_{ii}^2) \\
    & = \bm{v}^{\ast} (\mathbf{I} - \bm{Q}) f(\bm{A}^2) (\mathbf{I} - \bm{Q}) \bm{v} \\
    &\leq \lVert (\mathbf{I} - \bm{Q}) f(\bm{A}^2) (\mathbf{I} - \bm{Q}) \rVert,
\end{align*}
thus verifying the stated inequality \cref{eq:supporting_line}.
\end{proof}

We use \cref{lem:quasi_convex} to prove the main result, \cref{thm:gapless}.

\begin{proof}[Proof of \cref{thm:gapless}]
We assume that $\bm{A} \in \mathbb{R}^{N \times N}$ is diagonal and psd, with nonincreasing diagonal entries (\cref{lem:diagonal}).
By replacing $\bm{A}$ with $\bm{A} / \sigma_{r+1}(\bm{A})$ if necessary, we also assume that $\sigma_{r+1}(\bm{A}) = 1$.  
Last, we partition the test matrix
$\bm{\Omega} = \bigl[\begin{smallmatrix}
\bm{\Omega}_1 \\ \bm{\Omega}_2
\end{smallmatrix}\bigr]$ 
into a $r \times k$ matrix $\bm{\Omega}_1$ and a $(N - r) \times k$ matrix $\bm{\Omega}_2$.

To analyze \texttt{RSI}, we introduce the orthogonal projection $\bm{Q} = \bm{\Pi}_{\bm{A}^{m - 1} \bm{\Omega} \bm{\Omega}_1^{\dagger}}$ and apply \cref{lem:quasi_convex} with the convex function $f(x) = x^{m-1}$ to verify
\begin{equation*}
    \Expect\lVert (\mathbf{I} - \bm{Q}) \bm{A} \rVert^{2(m-1)}
    \leq \Expect\lVert (\mathbf{I} - \bm{Q}) \bm{A}^{m-1} \rVert^2.
\end{equation*}
We identify $\Expect\lVert (\mathbf{I} - \bm{Q}) \bm{A}^{m-1} \rVert^2$ as the mean-square error of \texttt{RSVD} applied to the matrix $\bm{A}^{m-1}$ with a test matrix $\bm{\Omega} \bm{\Omega}_1^{\dagger}$, and we mimic the proof of \cref{thm:main_thm} to show that
\begin{equation*}
    \Expect\lVert (\mathbf{I} - \bm{Q}) \bm{A}^{m-1} \rVert^2
    \leq 1 + \frac{r}{k-r-1} \sum\nolimits_{i = r+1}^N \sigma_i(\bm{A})^{2(m-1)}.
\end{equation*}
We take logarithms and divide through by $m - 1$ to obtain
\begin{equation*}
    \log\bigl(\Expect \lVert (\mathbf{I} - \bm{Q}) \bm{A} \rVert^2\bigr)
    \leq \frac{1}{m - 1} \log\biggl(1 + \frac{r}{k - r - 1}
    \sum\nolimits_{i = r + 1}^N \sigma_i(\bm{A})^{2(m-1)}
    \biggr),
\end{equation*}
which is a stronger version of the error bound for \texttt{RSI}.

To analyze \texttt{RBKI}, we define $\bm{Q} = \bm{\Pi}_{\bm{A} T_{m-2}(\bm{A}) \bm{\Omega} \bm{\Omega}_1^{\dagger}}$.
Then we apply \cref{lem:chebyshev} part (5) to identify the ``almost'' convex function $f(x) = x \bigl[T_{m-2}\bigl(x^{1 \slash 2}\bigr)\bigr]^2$, and we invoke \cref{lem:quasi_convex} to obtain
\begin{equation*}
    \bigl[T_{m-2}\bigl(\bigl(\Expect\lVert (\mathbf{I} - \bm{Q}) \bm{A} \rVert^2\bigr)^{1 \slash 2}\bigr)\bigr]^2
    \leq f(\Expect \lVert (\mathbf{I} - \bm{Q}) \bm{A} \rVert^2) 
    \leq \Expect \lVert (\mathbf{I} - \bm{Q}) \bm{A} T_{m-2}(\bm{A}) \rVert^2.
\end{equation*}
We recognize $\Expect \lVert (\mathbf{I} - \bm{Q}) \bm{A} T_{m-2}(\bm{A}) \rVert^2$ as the mean-square error of \texttt{RSVD} applied to the matrix $\bm{A} T_{m-2}(\bm{A})$.
Hence, by mimicking the proof of \cref{thm:main_thm}, we obtain
\begin{equation*}
    \Expect \lVert (\mathbf{I} - \bm{Q}) \bm{A} T_{m-2}(\bm{A}) \rVert^2
    \leq 1
    + \frac{r}{k-r-1} \sum\nolimits_{i = r+1}^N \sigma_i(\bm{A})^2.
\end{equation*}
Last, we take square roots, apply the inverse Chebyshev polynomial $T_{m-2}^{-1}$, take logarithms, and multiply through by $2$ to obtain
\begin{align*}
    \log\bigl(\Expect\lVert (\mathbf{I} - \bm{Q}) \bm{A} \rVert^2\bigr) 
    &\leq 2\log\biggl(T_{m-2}^{-1} \biggl(\sqrt{1
    + \frac{r}{k-r-1} \sum\nolimits_{i = r+1}^N \sigma_i(\bm{A})^2}\biggr)\biggr) \\
    &\leq \frac{1}{4(m-2)^2} \biggl[\log\biggl(4
    + \frac{4r}{k-r-1} \sum\nolimits_{i = r+1}^N \sigma_i(\bm{A})^2\biggr)\biggr]^2,
\end{align*}
which confirms the error bound for \texttt{RBKI}.
The last line uses the bound for $T_{m-2}^{-1}$ which was stated in \cref{lem:chebyshev} part 4.

We confirm the error bounds for \texttt{NysSI} and \texttt{NysBKI} in a similar fashion.
For \texttt{NysSI}, we define the orthogonal projection $\bm{Q} = \bm{\Pi}_{\bm{A}^{m-1/2} \bm{\Omega} \bm{\Omega}_1^{\dagger}}$ and calculate
\begin{multline*}
    \lVert \bm{A} - \hat{\bm{A}} \rVert \leq
    \lVert \bm{A} - \bm{A}\langle \bm{A}^{m-1} \bm{\Omega} \bm{\Omega}_1^{\dagger} \rangle \rVert
    = \lVert \bm{A}^{1/2}\bigl(\mathbf{I} - \bm{Q} \bigr) \bm{A}^{1/2} \rVert \\
    = \lVert \bigl(\mathbf{I} - \bm{Q} \bigr) \bm{A}^{1/2} \rVert^2
    = \lVert \bigl(\mathbf{I} - \bm{Q} \bigr) \bm{A} \bigl(\mathbf{I} - \bm{Q} \bigr) \rVert
    \leq \lVert \bigl(\mathbf{I} - \bm{Q} \bigr) \bm{A} \rVert.
\end{multline*}
By applying \cref{lem:quasi_convex} with $f(x) = x^{m - 1 \slash 2}$,
we obtain
\begin{equation*}
    \Expect\lVert (\mathbf{I} - \bm{Q}) \bm{A} \rVert^{2(m-1/2)}
    \leq \Expect\lVert (\mathbf{I} - \bm{Q}) \bm{A}^{m-1/2} \rVert^2,
\end{equation*}
which leads to an error bound for \texttt{NysSI}.

Last, for \texttt{NysBKI}, we define the orthogonal projection $\bm{Q} = \bm{\Pi}_{\bm{A} T_{2m-1}(\bm{A}^{1/2}) \bm{\Omega} \bm{\Omega}_1^{\dagger}}$ so that $\lVert \bm{A} - \hat{\bm{A}} \rVert
\leq \lVert \bigl(\mathbf{I} - \bm{Q} \bigr) \bm{A} \rVert$.
By applying \cref{lem:quasi_convex} with $f(x) = x \bigl[T_{2m-1}\bigl(x^{1 \slash 4}\bigr)\bigr]^2$, we obtain
\begin{equation*}
    \bigl[T_{2m-1}\bigl(\bigl(\Expect\lVert (\mathbf{I} - \bm{Q}) \bm{A} \rVert^2\bigr)^{1 \slash 4}\bigr)\bigr]^2
    \leq \Expect \lVert (\mathbf{I} - \bm{Q}) \bm{A} T_{2m-1}(\bm{A}^{1/2}) \rVert^2.
\end{equation*}
which leads to an error bound for \texttt{NysBKI} and completes the proof of \cref{thm:gapless}.
\end{proof}

\subsection{Gapped error bounds} \label{sec:frobenius}

To establish the gapped error bounds, we will need to apply a polynomial $\phi$ that enhances the singular value gap.
To that end, we establish the following filtering lemma.

\begin{lemma}[Enhancing the gap] \label{lem:enhancing}
    Choose a diagonal, psd matrix with nonincreasing entries: $\bm{A} = \diag(a_{11}, a_{22}, \dots, a_{NN})$.
    Fix two indices $1 \leq r < s \leq N$ and a function $\phi(x)$ that is positive for $x \geq a_{rr}$, and partition the matrix $\bm{A}$ as
    \begin{equation*}
        \bm{A} = \begin{bmatrix} \bm{A}_1 & & \\ & \bm{A}_2 & \\ & & \bm{A}_3 \end{bmatrix},
    \end{equation*}
    where $\bm{A}_1$, $\bm{A}_2$, and $\bm{A}_3$ are square matrices with dimensions $r$, $s - 1 - r$, and $N - s + 1$.
    If $\bm{\Omega} \in \mathbb{R}^{N \times k}$ is a Gaussian matrix with block size $k \geq s + 1$, then
    \begin{equation*}
        \Expect \lVert (\mathbf{I} - \bm{\Pi}_{\bm{A} \phi(\bm{A}) \bm{\Omega}}) \bm{A} \rVert_p^2
        \leq \Biggl\lVert 
        \begin{bmatrix}
        \bm{0} & ~ & ~ \\
        ~ & \bm{A}_2 & ~ \\
        ~ & ~ & \bm{A}_3
        \end{bmatrix} \Biggr\rVert_p^2
        + \frac{s - 1}{k - s} \cdot \frac{\lVert \bm{A}_3 \phi(\bm{A}_3) \rVert_p^2}{\phi(a_{rr})^2}.
    \end{equation*}
for any Schatten $p$-norm with $2 \leq p \leq \infty$.
\end{lemma}
\begin{proof}
We partition the matrix $\bm{\Omega}$ as
$\bm{\Omega} = \bigl[ \begin{smallmatrix} \bm{\Omega}_\alpha \\ \bm{\Omega}_\beta \end{smallmatrix}\bigr]$, where $\bm{\Omega}_\alpha$ is $(s - 1) \times k$ and $\bm{\Omega}_\beta$ is $(N - s + 1) \times k$.
Then we set $\bm{Q} = \bm{\Pi}_{\bm{A} \phi(\bm{A}) \bm{\Omega} \bm{\Omega}_\alpha^{\dagger}}$ and apply the coarse upper bound
\begin{equation*}
    \lVert (\mathbf{I} - \bm{Q}) \bm{A} \rVert_p^2
    \leq 
    \Biggl\lVert 
    \begin{bmatrix}
    \bm{0} & ~ & ~ \\
    ~ & \bm{A}_2 & ~ \\
    ~ & ~ & \bm{A}_3
    \end{bmatrix} \Biggr\rVert_p^2
    + \Biggl\lVert (\mathbf{I} - \bm{Q}) 
    \begin{bmatrix}
    \bm{A}_1 & ~ & ~ \\
    ~ & \bm{0} & ~ \\
    ~ & ~ & \bm{0}
    \end{bmatrix} \Biggr\rVert_p^2 \\
\end{equation*}
We can bound the second term using
\begin{align}
\label{eq:new_calc_2}
    & \Biggl\lVert (\mathbf{I} - \bm{Q}) 
    \begin{bmatrix}
    \bm{A}_1 \phi(\bm{A}_1) & ~ & ~ \\
    ~ & \bm{A}_2 \phi(\bm{A}_2) & ~ \\
    ~ & ~ & \bm{0}
    \end{bmatrix} 
    \begin{bmatrix}
    \phi(\bm{A}_1)^{-1} & ~ & ~ \\
    ~ & \bm{0} & ~ \\
    ~ & ~ & \bm{0}
    \end{bmatrix} \Biggr\rVert_p^2 \\
\label{eq:new_calc_3}
    &\qquad\leq 
    \frac{1}{\phi(a_{rr})^2}
    \Biggl\lVert (\mathbf{I} - \bm{Q}) 
    \begin{bmatrix}
    \bm{A}_1 \phi(\bm{A}_1) & ~ & ~ \\
    ~ & \bm{A}_2 \phi(\bm{A}_2) & ~ \\
    ~ & ~ & \bm{0}
    \end{bmatrix} 
    \Biggr\rVert_p^2 \\
\label{eq:new_calc_4}
    &\qquad\leq \frac{1}{\phi(a_{rr})^2}
    \lVert \bm{A}_3 \phi(\bm{A}_3) \bm{\Omega}_\beta \bm{\Omega}_\alpha^{\dagger}
    \rVert_p^2.
\end{align}
The inequality \cref{eq:new_calc_3} depends on the fact that the operator norm of the right-most matrix is $1 / \phi(a_{rr})$, while inequality \cref{eq:new_calc_4} follows from \cref{prop:ext_analysis}.
Last, we take expectations and apply the exact calculations for Gaussian matrices from \cref{prop:convenient} to show
\begin{equation*}
    \Expect\lVert \bm{A}_3 \phi(\bm{A}_3) \bm{\Omega}_\beta \bm{\Omega}_\alpha^{\dagger}
    \rVert_p^2
    \leq \Expect\lVert \bm{A}_3 \phi(\bm{A}_3) \bm{\Omega}_\beta \bm{\Omega}_\alpha^{\dagger}
    \rVert_{\rm F}^2
    = \frac{s - 1}{k - s} \Expect\lVert \bm{A}_3 \phi(\bm{A}_3) \rVert_{\rm F}^2,
\end{equation*}
which completes the proof.
\end{proof}

We apply \cref{lem:enhancing} to establish the main result, \cref{thm:gapped}.

\begin{proof}[Proof of \cref{thm:gapped}]
We assume without loss of generality that $\bm{A} \in \mathbb{R}^{N \times N}$ is diagonal and psd, with nonincreasing diagonal entries (\cref{lem:diagonal}).
We then apply \cref{lem:enhancing} with $\phi(x) = x^{m-2}$ to bound the \texttt{RSI} error as follows:
\begin{equation*}
    \Expect \lVert (\mathbf{I} - \bm{\Pi}_{\bm{A}^{m-1} \bm{\Omega}}) \bm{A} \rVert_p^2 \\
    \leq \lVert \bm{A} - \lfloor \bm{A} \rfloor_r \rVert_p^2
    + \frac{s-1}{k-s} \sum_{i=s}^N
    \frac{\sigma_i(\bm{A})^{2(m-1)}}{\sigma_r(\bm{A})^{2(m-2)}}.
\end{equation*}
\Cref{lem:chebyshev} part 3 gives the numerical identity
\begin{equation*}
    \biggl(\frac{\sigma_s(\bm{A})}{\sigma_r(\bm{A})}\biggr)^{2(m - 2)}
    = \biggl(\frac{1 - \gamma}{1 + \gamma} \biggr)^{2(m - 2)} 
    \leq \mathrm{e}^{-4(m - 2) \gamma},
\end{equation*}
which confirms the \texttt{RSI} error bound.
The \texttt{NysSI} error bound is proved similarly.
Observe that
\begin{multline*}
    \lVert \bm{A} - \bm{A}\langle \bm{A}^{m - 1} \bm{\Omega} \rangle \rVert_p 
    = \lVert \bm{A}^{1/2}(\mathbf{I} - \bm{\Pi}_{\bm{A}^{m - 1/2} \bm{\Omega}}) \bm{A}^{1/2} \rVert_p
    = \lVert (\mathbf{I} - \bm{\Pi}_{\bm{A}^{m - 1/2} \bm{\Omega}}) \bm{A}^{1/2} \rVert_{2p}^2 \\
    = \lVert (\mathbf{I} - \bm{\Pi}_{\bm{A}^{m - 1/2} \bm{\Omega}}) \bm{A}
    (\mathbf{I} - \bm{\Pi}_{\bm{A}^{m - 1/2} \bm{\Omega}}) \rVert_p
    \leq \lVert (\mathbf{I} - \bm{\Pi}_{\bm{A}^{m-1/2} \bm{\Omega}}) \bm{A} \rVert_p.
\end{multline*}
We bound $\Expect \lVert (\mathbf{I} - \bm{\Pi}_{\bm{A}^{m-1/2} \bm{\Omega}}) \bm{A} \rVert_p^2$ by applying \cref{lem:enhancing} with $\phi(x) = x^{m - 3/2}$.

To control the \texttt{RBKI} error, we apply \cref{lem:enhancing} with $\phi(x) = T_{m-2}(x \slash \sigma_s(\bm{A}))$, which gives
\begin{equation*}
    \Expect \lVert \bm{A} - \bm{\Pi}_{\bm{A} \phi(\bm{A}) \bm{\Omega}} \bm{A} \rVert_p^2
    \leq \lVert \bm{A} - \lfloor \bm{A} \rfloor_r \rVert_p^2
    + \frac{s - 1}{k - s} \sum_{i = s}^N \sigma_i(\bm{A})^2 \biggl(\frac{\phi(\sigma_i(\bm{A}))}
    {\phi(\sigma_r(\bm{A}))}\biggr)^2.
\end{equation*}
Then we use the properties of Chebyshev polynomials (\cref{lem:chebyshev} parts (2) and (3)) to check that $|\phi(\sigma_i(\bm{A}))| \leq 1$ for $i = s, \ldots, N$ and
\begin{equation*}
    \phi(\sigma_r(\bm{A}))
    = T_{m-2}\biggl(\frac{\sigma_{r}(\bm{A})}{\sigma_s(\bm{A})}\biggr)
    = T_{m-2}\biggl(\frac{1 + \gamma}{1 - \gamma}\biggr)
    \geq \frac{1}{2} \mathrm{e}^{2(m-2) \sqrt{\gamma}}.
\end{equation*}
Similarly, to bound the \texttt{NysBKI} error, we define $\phi(x) = T_{2m - 3}(\sqrt{x \slash \sigma_s(\bm{A})})$ so that $\lVert \bm{A} - \hat{\bm{A}} \rVert_p \leq \lVert \bm{A} - \bm{\Pi}_{\bm{A} \phi(\bm{A}) \bm{\Omega}} \bm{A} \rVert_p$.
We apply \cref{lem:enhancing}, which gives
\begin{equation*}
    \Expect \lVert \bm{A} - \bm{\Pi}_{\bm{A} \phi(\bm{A}) \bm{\Omega}} \bm{A} \rVert_p^2 
    \leq \lVert \bm{A} - \lfloor \bm{A} \rfloor_r \rVert_p^2
    + \frac{s - 1}{k - s} \sum_{i = s}^N \sigma_i(\bm{A})^2 \biggl(\frac{\phi(\sigma_i(\bm{A}))}
    {\phi(\sigma_r(\bm{A}))}\biggr)^2.
\end{equation*}
We bound $|\phi(\sigma_i(\bm{A}))| \leq 1$ for $i = s, \ldots, N$ and
\begin{equation*}
    \phi(\sigma_r(\bm{A}))
    = T_{2m-3}\biggl(\sqrt{\frac{1 + \gamma}{1 - \gamma}}\biggr)
    \geq T_{2m-3}\biggl(\frac{1 + \gamma \slash 2}{1 - \gamma \slash 2}\biggr)
    \geq \frac{1}{2} \mathrm{e}^{2(m-3 \slash 2) \sqrt{2 \gamma}},
\end{equation*}
where we have used the numerical identity
\begin{equation*}
    \sqrt{\frac{1 + \gamma}{1 - \gamma}} \geq 
    \frac{1 + \gamma \slash 2}{1 - \gamma \slash 2}, \qquad \text{for } \gamma \geq 0.
\end{equation*}
This completes the proof.
\end{proof}

\section{Conclusion} \label{sec:conclusion}

In this work, we have described randomized low-rank approximation algorithms and provided user recommendations regarding which algorithms are the most efficient.
We have focused on computations involving high-dimensional matrices, in which matrix multiplications are the main computational bottleneck.
We have established simple, explicit bounds which allow for precise descriptions of the differences between algorithms.
Our results demonstrate that \texttt{RSVD} and \texttt{NysSVD} are fast and accurate when the singular values of the target matrix decay quickly.
However, in settings of slow singular value decay,
\texttt{RBKI} and \texttt{NysBKI} are the available algorithms with the greatest speed and robustness.
We have presented numerical tests demonstrating the utility of these Krylov methods for principal component analysis and kernel spectral clustering.

The widespread adoption of \texttt{RBKI} and \texttt{NysBKI} will require changes in software, since \texttt{RSI} is currently the default randomized low-rank matrix approximation algorithm in Matlab \cite{matlab2020version} and sci-kit learn \cite{scikit-learn}.
Yet, \texttt{RBKI} and \texttt{NysBKI} are of great benefit to computational scientists who have found \texttt{RSI} to be expensive while resulting in large random errors \cite{chen2013improved,kumar2019target,yang2021low,bjarkason2018randomized,bousserez2020enhanced,lee2014large}.
\texttt{RBKI} and \texttt{NysBKI} are more accurate and scalable algorithms.
To support the broader use of these Krylov methods, we have provided simple and stable pseudocode that cuts down on the cost of \texttt{RBKI} by roughly $33\%$ of the number of matrix--vector products, and we have provided the first pseudocode for \texttt{NysBKI}, which is faster than \texttt{RBKI} by a factor of $\sqrt{2}$ matrix--vector products.

\section*{Acknowledgments}
We acknowledge many helpful conversations with Tyler Chen, Mateo D{\'i}az, Barbara Engelhardt, Ethan Epperly, Maksim Melnichenko, Riley Murray, and Christopher Musco.

\appendix

\section{Quality assurance for singular vectors}

\Cref{alg:rbki_quality,alg:nysbki_quality} are adaptive versions of \texttt{RBKI} and \texttt{NysBKI} that we developed with assistance from Maksim Melnichenko and Riley Murray.
The algorithms calculate the singular vector residuals and stop when each of the first $r$ singular vector residuals falls below a threshold $\varepsilon$.

In these algorithms, we use several shortcuts to evaluate the residual error formula \cref{eq:residual}.
If the low-rank approximation takes the form $\hat{\bm{A}} = \bm{\Pi}_{\bm{M}} \bm{A}$, the first term in \cref{eq:residual} is zero and the residual is $r_i = \lVert (\bm{A} - \hat{\bm{A}}) \hat{\bm{v}}_i \rVert.$
If instead $\hat{\bm{A}} = \bm{A} \bm{\Pi}_{\bm{M}}$, the second term in \cref{eq:residual} is zero and $r_i = \lVert (\bm{A} - \hat{\bm{A}})^\ast \hat{\bm{u}}_i \rVert.$
Last, if
$\hat{\bm{A}} = \bm{A}\langle \bm{M} \rangle$ is a Nystr\"om approximation, the left and right singular vectors are the same, and $r_i = \sqrt{2} \lVert (\bm{A} - \hat{\bm{A}}) \hat{\bm{v}}_i \rVert$.

\begin{algorithm}[t]
\caption{\texttt{RBKI} with quality assurance \label{alg:rbki_quality}}
\begin{algorithmic}[1]
\Require Matrix $\bm{A} \in \mathbb{R}^{L \times N}$; block size $k$; number of singular vectors $r$; tolerance $\varepsilon_{\textup{tol}}$
\Ensure Orthogonal $\bm{U}$,
orthogonal $\bm{V}$, and diagonal $\bm{\Sigma}$ such that $\bm{A} \approx \bm{U} \bm{\Sigma} \bm{V}^{\ast}$
\State Generate a random matrix $\bm{Y}_0 \in \mathbb{R}^{N \times k}$
\State $i = 1$
\State $\bm{X}_1 = \bm{A} \bm{Y}_0$
\State $[\bm{X}_1, \sim] = \textup{qr\_econ}(\bm{X}_1)$
\State $\bm{W}_2 = \bm{A}^\ast \bm{X}_1$
\State $\bm{R} = [\,]; \bm{S} = [\,]$; $\bm{E} = [\,]$
\While{$i = 1$ or one of the first $r$ columns of $\bm{E}$ has norm $> \varepsilon$}
\State $i = i + 1$
\If{$i$ is even}
    \State $\bm{Y}_i = \bm{W}_i$
    \State $\bm{R}_{\bullet i} = \begin{bmatrix} \bm{Y}_2 & \bm{Y}_4 & \cdots & \bm{Y}_{i-2} \end{bmatrix}^\ast \bm{Y}_i$
    \State $\bm{Y}_i = \bm{Y}_i - \sum_{\textup{even } j < i} \bm{Y}_j
    (\bm{Y}_j^\ast \bm{Y}_i)$
    \Comment{Orthog. w.r.t. even iterates (2$\times$)} 
    \State $[\bm{Y}_i, \bm{R}_{ii}] = \textup{qr\_econ}(\bm{Y}_i)$
    \Comment{Optional: stabilized QR \cref{eq:stabilized_qr}}
    \State $\bm{R} = \begin{bmatrix} \bm{R} & \bm{R}_{\bullet i} \\ \bm{0} & \bm{R}_{ii} \end{bmatrix}$
    \State $\bm{Z}_{i+1} = \bm{A} \bm{Y}_i$
    \State $[\hat{\bm{U}}, \bm{\Sigma}, \hat{\bm{V}}] = 
    \textup{svd\_econ} (\bm{R}^\ast)$
    \State $\bm{E} = \begin{bmatrix} \bm{Z}_3 & \bm{Z}_5 & \cdots \end{bmatrix} \hat{\bm{V}} 
    - \begin{bmatrix} \bm{X}_1 & \bm{X}_3 & \cdots \end{bmatrix} \hat{\bm{U}} \bm{\Sigma}$
    \Comment{Residual evaluation}
\Else
    \State $\bm{X}_i = \bm{Z}_i$
    \State $\bm{S}_{\bullet i} = \begin{bmatrix} \bm{X}_1 & \bm{X}_3 & \cdots & \bm{X}_{i-2} \end{bmatrix}^\ast \bm{X}_i$
    \State $\bm{X}_i = \bm{X}_i - \sum_{\textup{odd } j < i} \bm{X}_j
    (\bm{X}_j^\ast \bm{X}_i)$
    \Comment{Orthog. w.r.t. odd iterates (2$\times$)} 
    \State $[\bm{X}_i, \bm{S}_{ii}] = \textup{qr\_econ}(\bm{X}_i)$
    \Comment{Optional: stabilized QR \cref{eq:stabilized_qr}}
    \State $\bm{S} = \begin{bmatrix} \bm{S} & \bm{S}_{\bullet i} \\ \bm{0} & \bm{S}_{ii} \end{bmatrix}$
    \State $\bm{W}_{i+1} = \bm{A}^\ast \bm{X}_i$
    \State $[\hat{\bm{U}}, \bm{\Sigma}, \hat{\bm{V}}] = 
    \textup{svd\_econ} (\bm{S})$
    \State $\bm{E} = \begin{bmatrix} \bm{W}_2 & \bm{W}_4 & \cdots \end{bmatrix} \hat{\bm{U}} 
    - \begin{bmatrix} \bm{Y}_2 & \bm{Y}_4 & \cdots \end{bmatrix} \hat{\bm{V}} \bm{\Sigma}$
    \Comment{Residual evaluation}
\EndIf
\EndWhile
\State $\bm{U} = \begin{bmatrix} \bm{X}_1 & \bm{X}_3 & \cdots \end{bmatrix} \hat{\bm{U}}$
\State $\bm{V} = \begin{bmatrix} \bm{Y}_2 & \bm{Y}_4 & \cdots \end{bmatrix} \hat{\bm{V}}$
\end{algorithmic}
\end{algorithm}

\begin{algorithm}[t]
\caption{\texttt{NysBKI} with quality assurance \label{alg:nysbki_quality}}
\begin{algorithmic}[1]
\Require Psd matrix $\bm{A} \in \mathbb{R}^{N \times N}$; block size $k$; number of eigenvectors $r$; tolerance $\varepsilon_{\textup{tol}}$; shift $\varepsilon > 0$
\Ensure Orthogonal $\bm{U}$ and diagonal $\bm{\Lambda}$ such that $\bm{A} \approx \bm{U} \bm{\Lambda} \bm{U}^{\ast}$
\State Generate a random matrix $\bm{Y}_0 \in \mathbb{R}^{N \times k}$
\State $[\bm{X}_0, \sim] = \textup{qr\_econ}(\bm{Y}_0)$
\State $i = 1$
\State $\bm{Y}_1 = \bm{A} \bm{X}_0$
\State $\bm{S}_0 = [\,]$; $\bm{E} = [\,]$
\While{$i = 1$ or one of the first $r$ columns of $\bm{E}$ has norm $> \varepsilon$}
\State $\bm{X}_i = \bm{Y}_i + \varepsilon \bm{X}_{i-1}$
\State $\bm{R}_{\bullet i}
= \begin{bmatrix} \bm{0} & \cdots & \bm{0} & \bm{X}_{i-2} & \bm{X}_{i-1} \end{bmatrix}^\ast \bm{X}_i$
\Comment{$\bm{R}_{\bullet i} = \bm{X}_{i-1}^\ast \bm{X}_i$ if $i = 1$}
\State $\bm{X}_i = \bm{X}_i - \sum_{j < i} \bm{X}_j (\bm{X}_j^{\ast} \bm{X}_i)$
\Comment{Orthog.~w.r.t.~past iterates (2$\times$)}
\State $[\bm{X}_i, \bm{R}_{ii}] = \textup{qr\_econ}(\bm{X}_i)$
\Comment{Optional: stabilized QR \cref{eq:stabilized_qr}}
\State $\bm{C} = \textup{chol}\left(\begin{bmatrix} \bm{S}_{i-1} & \bm{R}_{\bullet i} \end{bmatrix}
\right)$
\State $\bm{S}_i = \begin{bmatrix} \bm{S}_{i-1} & \bm{R}_{\bullet i} \\
\bm{0} & \bm{R}_{ii}
\end{bmatrix}$
\Comment{$\bm{A} \begin{bmatrix} \bm{X}_0 & \cdots & \bm{X}_{i-1} \end{bmatrix}
= \begin{bmatrix} \bm{X}_0 & \cdots & \bm{X}_i \end{bmatrix} \bm{S}_i$}
\State $\bm{Z} = \bm{S}_i \bm{C}^{-1}$
\State $[\hat{\bm{U}}, \bm{\Sigma}, \sim] = \textup{svd\_econ}(\bm{Z})$
\State $\bm{U} = \begin{bmatrix} \bm{X}_0 & \cdots & \bm{X}_i \end{bmatrix} \hat{\bm{U}}$
\State $\bm{\Lambda} = \max\{\bm{0}, \bm{\Sigma}^2 - \varepsilon \mathbf{I}\}$
\State $\bm{Y}_{i+1} = \bm{A} \bm{X}_i$
\State
$\bm{E} = \begin{bmatrix} \bm{Y}_1 & \cdots & \bm{Y}_{i+1} \end{bmatrix} \hat{\bm{U}}
- \bm{U} \bm{\Sigma}$
\Comment{Residual evaluation}
\State $i = i + 1$
\EndWhile
\end{algorithmic}
\end{algorithm}

\section{Random matrix theory} \label{app:random}
In this section, we bound the moments (\cref{lem:moment}) and inverse moments (\cref{lem:inverse_moment,lem:moment_spectral}) of Gaussian random matrices.
These are the sharpest available expectation bounds and concentration inequalities of their type, and they are helpful for analyzing randomized matrix computations.
We use these results to prove \cref{prop:convenient} in \cref{app:convenient} and 
\cref{thm:more_bounds} in \cref{app:more_bounds}.

\begin{lemma}[Moment bounds] \label{lem:moment}
Consider fixed matrices $\bm{S} \in \mathbb{R}^{L \times N}$ and $\bm{T} \in \mathbb{R}^{P \times Q}$
and a Gaussian matrix $\bm{G} \in \mathbb{R}^{N \times P}$. Then, $\bm{S} \bm{G} \bm{T}$ satisfies the equalities
\begin{align*}
    & \Expect\lVert \bm{S} \bm{G} \bm{T} \rVert_{\rm F}^2 = \lVert \bm{S} \rVert_{\rm F}^2 \lVert \bm{T} \rVert_{\rm F}^2, \\
    & \Expect\lVert \bm{S} \bm{G} \bm{T} \rVert_4^4 = \lVert \bm{S} \rVert_4^4 \lVert \bm{T} \rVert_{\rm F}^4
    + \lVert \bm{T} \rVert_4^4 \lVert \bm{S} \rVert_{\rm F}^4
    + \lVert \bm{S} \rVert_4^4 \lVert \bm{T} \rVert_4^4,
\end{align*}
as well as the upper bounds
\begin{align*}
    & (\Expect\lVert \bm{S} \bm{G} \bm{T} \rVert^2)^{1 \slash 2} \leq \lVert \bm{S} \rVert \lVert \bm{T} \rVert_{\rm F} +
    \lVert \bm{T} \rVert \lVert \bm{S} \rVert_{\rm F}, \\
    & (\Expect\lVert \bm{S} \bm{G} \bm{T} \rVert^4)^{1 \slash 4} 
    \leq \lVert \bm{S} \rVert 
    (\lVert \bm{T} \rVert_{\rm F}^4
    + 2\lVert \bm{T} \rVert_4^4)^{1 \slash 4}
    + \lVert \bm{T} \rVert 
    (\lVert \bm{S} \rVert_{\rm F}^4
    + 2\lVert \bm{S} \rVert_4^4)^{1 \slash 4}.
\end{align*}
\end{lemma}
\begin{proof}
The equalities come from direct calculations, e.g.,
\begin{equation*}
    \Expect\lVert \bm{S} \bm{G} \bm{T} \rVert_{\rm F}^2
    = \Expect \textup{tr} \bigl(\bm{S} \bm{G} \bm{T} \bm{T}^\ast \bm{G}^\ast \bm{S}^\ast\bigr)
    = \Expect \sum_{ijkl} \bm{S}_{ij} \bm{G}_{jk} \bm{T}_{kl} \bm{T}_{kl} \bm{G}_{jk} \bm{S}_{ij}
    = \sum_{ij} \bm{S}_{ij}^2 \sum_{kl} \bm{T}_{kl}^2.
\end{equation*}
We have used the fact that unmatched Gaussians $\bm{G}_{ab} \bm{G}_{cd}$ for $(a,b) \neq (c,d)$ vanish in expectation, so it is only necessary to consider the matched Gaussian terms $\bm{G}_{ab} \bm{G}_{ab}$.

The inequalities come from a comparison principle for Gaussian random variables.
We compare two Gaussian random fields
\begin{align*}
    X_{\bm{u} \bm{v}} &= \langle \bm{S}^\ast \bm{u}, \bm{G} \bm{T} \bm{v} \rangle + \lVert \bm{S}^\ast \bm{u} \rVert \lVert \bm{T} \bm{v} \rVert \gamma, \\
    Y_{\bm{u} \bm{v}} &= \lVert \bm{S}^\ast \bm{u} \rVert \langle \bm{h}, \bm{T} \bm{v} \rangle
    + \lVert \bm{T} \bm{v} \rVert \langle \bm{g}, \bm{S}^\ast \bm{u} \rangle,
\end{align*}
where $\gamma \in \mathbb{R}$ is a Gaussian random variable
and $\bm{h} \in \mathbb{R}^p$ and $\bm{g} \in \mathbb{R}^n$ are Gaussian random vectors,
independent from one another and independent from $\bm{G}$.
The indices $\bm{u} \in \mathbb{R}^L$ and $\bm{v} \in \mathbb{R}^Q$ are unit-length vectors satisfying $\lVert \bm{u} \rVert = \lVert \bm{v} \rVert = 1$.
A short calculation verifies that
\begin{align*}
    \Expect [X_{\bm{u} \bm{v}} X_{\bm{u}^\prime \bm{v}^\prime}]
    &= \langle \bm{S}^\ast \bm{u},
    \bm{S}^\ast \bm{u}^\prime \rangle
    \langle \bm{T} \bm{v}, \bm{T} \bm{v}^\prime\rangle
    + \lVert \bm{S}^\ast \bm{u} \rVert 
    \lVert \bm{S}^\ast \bm{u}^\prime \rVert
    \lVert \bm{T} \bm{v} \rVert
    \lVert \bm{T} \bm{v}^\prime \rVert, \\
    \Expect [Y_{\bm{u} \bm{v}} Y_{\bm{u}^\prime \bm{v}^\prime}] 
    &= \lVert \bm{S}^\ast \bm{u} \rVert \lVert \bm{S}^\ast \bm{u}^\prime \rVert 
    \langle \bm{T} \bm{v}, \bm{T} \bm{v}^\prime\rangle
    + \langle \bm{S}^\ast \bm{u},
    \bm{S}^\ast \bm{u}^\prime \rangle
    \lVert \bm{T} \bm{v} \rVert
    \lVert \bm{T} \bm{v}^\prime \rVert.
\end{align*}
Hence,
\begin{multline*}
    \Expect [X_{\bm{u} \bm{v}} X_{\bm{u}^\prime \bm{v}^\prime}]
    - \Expect [Y_{\bm{u} \bm{v}} Y_{\bm{u}^\prime \bm{v}^\prime}] \\
    = (\lVert \bm{S}^\ast \bm{u} \rVert \lVert \bm{S}^\ast \bm{u}^\prime \rVert - 
    \langle \bm{S}^\ast \bm{u},
    \bm{S}^\ast \bm{u}^\prime \rangle)
    (\lVert \bm{T} \bm{v} \rVert
    \lVert \bm{T} \bm{v}^\prime \rVert -
    \langle \bm{T} \bm{v}, \bm{T} \bm{v}^\prime\rangle) \geq 0.
\end{multline*}
Consequently, the increments of $X$ and $Y$ satisfy
\begin{equation*}
    \Expect | X_{\bm{u} \bm{v}} - X_{\bm{u}^\prime \bm{v}^\prime} |^2
    \leq \Expect | Y_{\bm{u} \bm{v}} - Y_{\bm{u}^\prime \bm{v}^\prime} |^2,
\end{equation*}
and Slepian's lemma \cite[Sec.~7.2]{vershynin2018high} gives the result
\begin{equation*}
    \Expect | \max_{\bm{u},\bm{v}} X_{\bm{u} \bm{v}} |^p
    \leq \Expect | \max_{\bm{u},\bm{v}} Y_{\bm{u} \bm{v}} |^p
\end{equation*}
for any $p \geq 1$.
Using Jensen's inequality, we find also
\begin{equation*}
    \Expect \lVert \bm{S} \bm{G} \bm{T} \rVert^p
    = \Expect \max_{\bm{u}, \bm{v}}
    \bigl| \bm{u}^\ast \bm{S} \bm{G} \bm{T} \bm{v}\bigr|^p
    \leq \Expect | \max_{\bm{u},\bm{v}} X_{\bm{u} \bm{v}} |^p \leq \Expect | \max_{\bm{u},\bm{v}} Y_{\bm{u} \bm{v}} |^p.
\end{equation*}
To bound $\Expect | \max_{\bm{u},\bm{v}} Y_{\bm{u} \bm{v}} |^p$ in the case $p = 2$, we use the triangle inequality and a direct calculation involving second moments of Gaussians to calculate
\begin{align*}
    \bigl(\Expect | \max_{\bm{u},\bm{v}} Y_{\bm{u} \bm{v}} |^2\bigr)^{1 \slash 2}
    &\leq \bigl( \Expect \bigl| \lVert \bm{S} \rVert \lVert \bm{T}^\ast \bm{h} \rVert
    + \lVert \bm{T} \rVert 
    \lVert \bm{S} \bm{g} \rVert \bigr|^2 \bigr)^{1 \slash 2} \\
    &\leq \bigl( \Expect \bigl| \lVert \bm{S} \rVert \lVert \bm{T}^\ast \bm{h} \rVert \bigr|^2 \bigr)^{1 \slash 2} 
    + \bigl( \Expect \bigl| \lVert \bm{T} \rVert 
    \lVert \bm{S} \bm{g} \rVert \bigr|^2 \bigr)^{1 \slash 2} \\
    &= \lVert \bm{S} \rVert \lVert \bm{T} \rVert_{\rm F}
    + \lVert \bm{T} \rVert \lVert \bm{S} \rVert_{\rm F}.
\end{align*}
Similarly, to bound $\Expect | \max_{\bm{u},\bm{v}} Y_{\bm{u} \bm{v}} |^p$ in the case $p = 4$, we calculate
\begin{align*}
    \bigl(\Expect | \max_{\bm{u},\bm{v}} Y_{\bm{u} \bm{v}} |^4\bigr)^{1 \slash 4}
    &\leq \bigl( \Expect \bigl| \lVert \bm{S} \rVert \lVert \bm{T}^\ast \bm{h} \rVert \bigr|^4 \bigr)^{1 \slash 4} 
    + \bigl( \Expect \bigl| \lVert \bm{T} \rVert 
    \lVert \bm{S} \bm{g} \rVert \bigr|^4 \bigr)^{1 \slash 4} \\
    &= \lVert \bm{S} \rVert 
    (\lVert \bm{T} \rVert_{\rm F}^4
    + 2\lVert \bm{T} \rVert_4^4)^{1 \slash 4}
    + \lVert \bm{T} \rVert 
    (\lVert \bm{S} \rVert_{\rm F}^4
    + 2\lVert \bm{S} \rVert_4^4)^{1 \slash 4},
\end{align*}
where the last equality comes from a direct calculation involving fourth moments of Gaussians.
\end{proof}

\begin{lemma}[Inverse moment bounds]
\label{lem:inverse_moment}
Consider a Gaussian matrix $\bm{G} \in \mathbb{R}^{r \times k}$ with $r \leq k$.
Then, $(\bm{G} \bm{G}^\ast)^{-1}$ satisfies the moment formulas
\begin{align*}
    & \Expect (\bm{G} \bm{G}^\ast)^{-1} = \frac{1}{k-r-1} \mathbf{I}, && r \leq k - 2 \\
    & \Expect\lVert (\bm{G} \bm{G}^\ast)^{-1} \rVert_{\rm \ast}^2 = 
    \frac{r^2(k-r) - 2r(r-1)}{(k-r)(k-r-1)(k-r-3)}, && r \leq k - 4 \\
    & \Expect\lVert (\bm{G} \bm{G}^\ast)^{-1} \rVert_{\rm F}^2 = 
    \frac{r(k-1)}{(k-r)(k-r-1)(k-r-3)}, && r \leq k - 4,
\end{align*}
and the upper bounds
\begin{align}
\label{eq:upper_1}
    \Expect\bigl[a:\textup{tr} (\bm{G} \bm{G}^\ast)^{-1}\bigr]
    &\leq  \frac{r}{2} \log(1 + 2a), && r = k - 1, \\
\label{eq:upper_2}
    \Expect\bigl[a:\textup{tr} (\bm{G} \bm{G}^\ast)^{-1}\bigr]
    &\leq  \sqrt{\frac{\pi a}{2}}, && r = k.
\end{align}
Last, for any $r \leq k$ and $t \geq 0$,
\begin{equation}
\label{eq:upper_3}
    \mathbb{P}\biggl\{\textup{tr} (\bm{G} \bm{G}^\ast)^{-1}
    > \frac{{\rm e}tr}{k-r+1}\biggr\}
    \leq \sqrt{\frac{\pi r}{k - r + 1}} t^{-(k-r+1) \slash 2}.
\end{equation}
\end{lemma}
\begin{proof}
Explicit formulas for $\Expect (\bm{G} \bm{G}^\ast)^{-1}$, $\Expect\lVert (\bm{G} \bm{G}^\ast)^{-1} \rVert_{\rm \ast}^2$,
and $\Expect\lVert (\bm{G} \bm{G}^\ast)^{-1} \rVert_{\rm F}^2$
are given in \cite[pg.~119]{press2005applied}.  The remaining formulas are new,
but the expectation bounds are related to the computations
in~\cite{tropp2022randomized}.

To establish these formulas, we begin by calculating
\begin{equation*}
    \Expect\bigl[a: \textup{tr} (\bm{G} \bm{G}^\ast)^{-1}\bigr]
    = \Expect\biggl[a: \sum_{i=1}^r (\bm{G} \bm{G}^\ast)^{-1}_{ii}\biggr]
    \leq \sum_{i=1}^r \Expect\biggl[a: (\bm{G} \bm{G}^\ast)^{-1}_{ii}\biggr],
\end{equation*}
where we have used the subadditivity of $x \mapsto a:x$ and the linearity of expectation.
The inverse diagonal elements $1 / (\bm{G} \bm{G}^\ast)^{-1}_{ii}$ for $1 \leq i \leq r$ have chi-squared distributions with $p = k-r+1$ degrees of freedom \cite[pg.~119]{press2005applied}. 
Hence, their density function is
\begin{equation*}
    f(x) = \frac{1}
    {2^{p \slash 2} \Gamma(p \slash 2)} x^{p \slash 2 - 1} \mathrm{e}^{-x \slash 2}, \qquad x > 0.
\end{equation*}
In the case $p = 2$, this allows us to evaluate
\begin{equation*}
    \Expect\biggl[a: (\bm{G} \bm{G}^\ast)^{-1}_{ii}\biggr]
    = \frac{1}{2} \int_0^\infty \frac{\mathrm{e}^{-x \slash 2}}{a^{-1} + x} \mathop{dx}
    = \int_{1 \slash (2a)}^\infty \frac{\mathrm{e}^{1 \slash (2a) - y}}{y} \mathop{dy}
    \leq \frac{1}{2} \log(1 + 2a),
\end{equation*}
where we have substituted $y = x \slash 2 + 1 \slash (2a)$
and applied the bound \cite[Sec.~5.1.20]{abramowitz1988handbook}
\begin{equation*}
    \int_x^{\infty} \frac{\mathrm{e}^{-t}}{t} \mathop{dt} \leq \log\biggl(1 + \frac{1}{x}\biggr) \mathrm{e}^{-x}, \qquad x > 0.
\end{equation*}
In the case $p = 1$, we evaluate
\begin{equation*}
    \Expect\biggl[a: (\bm{G} \bm{G}^\ast)^{-1}_{ii}\biggr]
    \leq
    \frac{1}{2 \sqrt{\pi}} \int_0^\infty \frac{x^{-1 \slash 2}}{a^{-1} + x} \mathop{dx}
    = \sqrt{\frac{a}{\pi}} \int_0^\infty \frac{1}{1 + y^2} \mathop{dy}
    = \sqrt{\frac{\pi a}{2}},
\end{equation*}
where we have substituted $x = a^{-1} y^2$ and used the fact that $1 \slash (1 + y^2)$ is the derivative of $\arctan y$.
Thus we establish the bounds \cref{eq:upper_1,eq:upper_2}.

Last, to establish the concentration inequality \cref{eq:upper_3},
we set $u = {\rm e} t r \slash p$ and argue using Markov's inequality
\begin{equation*}
    \mathbb{P}\bigl\{\textup{tr} (\bm{G} \bm{G}^\ast)^{-1}
    > u \bigr\}
    = \mathbb{P}\bigl\{u^p : \lVert (\bm{G} \bm{G}^\ast)^{-1} \rVert_{\rm \ast}^p
    > u^p : u^p \biggr\}
    \leq \frac{\Expect \bigl[u^p : \lVert (\bm{G} \bm{G}^\ast)^{-1} \rVert_{\rm \ast}^p \bigr]}{u^p : u^p}.
\end{equation*}
We directly evaluate $u^p : u^p = \frac{1}{2} u^p$,
and we apply the monotonicity and subadditivity of $x \mapsto u^p:x$ to obtain the upper bound
\begin{align*}
    u^p : \lVert (\bm{G} \bm{G}^\ast)^{-1} \rVert_{\rm \ast}^p
    &= u^p : \biggl(\sum_{i=1}^r (\bm{G} \bm{G}^\ast)^{-1}_{ii} \biggr)^{p} 
    \leq u^p : r^{p - 1} \sum_{i=1}^r \bigl( (\bm{G} \bm{G}^\ast)^{-1}_{ii} \bigr)^{p} \\
    & \leq \sum_{i=1}^r u^p : r^{p - 1} \bigl( (\bm{G} \bm{G}^\ast)^{-1}_{ii} \bigr)^{p}
    = r^{p - 1} \sum_{i=1}^r \frac{u^p}{r^{p-1}} : \bigl( (\bm{G} \bm{G}^\ast)^{-1}_{ii} \bigr)^p.
\end{align*}
Hence, we have shown
\begin{equation*}
    \mathbb{P}\bigl\{\textup{tr} (\bm{G} \bm{G}^\ast)^{-1}
    > u \bigr\} \leq \frac{2r^{p-1}}{u^p} \sum_{i=1}^r \Expect\biggl[\frac{u^p}{r^{p-1}} : \bigl( (\bm{G} \bm{G}^\ast)^{-1}_{ii} \bigr)^{p}\biggr],
\end{equation*}
Last, we set $a = u^p \slash r^{p-1}$ and for $1 \leq i \leq r$ we evaluate
\begin{align*}
    \Expect\biggl[a : \bigl( (\bm{G} \bm{G}^\ast)^{-1}_{ii} \bigr)^{p}\biggr]
    &\leq
    \frac{1}{2^{p \slash 2} \Gamma(p \slash 2)} 
    \int_0^\infty \frac{x^{p \slash 2 - 1}}{a^{-1} + x^p} \mathop{dx} \\
    &= \frac{\pi \sqrt{a}}{2^{p \slash 2 + 1} \Gamma(p \slash 2 + 1)}
    \leq \frac{1}{2} \sqrt{\frac{\pi a}{p}} \biggl(\frac{\rm e}{p}\biggr)^{p \slash 2},
\end{align*}
where the last inequality follows from Stirling's approximation \cite[Sec.~6.1.38]{abramowitz1988handbook}.
We obtain
\begin{equation*}
    \mathbb{P}\bigl\{\textup{tr} (\bm{G} \bm{G}^\ast)^{-1}
    > u \bigr\} 
    \leq 
    \sqrt{\frac{\pi r}{p}} \biggl(\frac{{\rm e} r}{up}\biggr)^{p \slash 2},
\end{equation*}
and substituting $u = {\rm e} tr \slash p$ completes the proof.
\end{proof}

\begin{lemma}[Inverse moment bounds, spectral norm]
\label{lem:moment_spectral}
For $0 \leq p \leq 18$ and $k \geq r + 2m$, the 
Gaussian matrix $\bm{G} \in \mathbb{R}^{r \times k}$ satisfies
\begin{equation}
\label{eq:expectation}
    \bigl( \Expect\lVert (\bm{G} \bm{G}^\ast)^{-1}\rVert^p \bigr)^{1 \slash p}
    \leq \frac{{\rm e}^2(k + r)}{2(k - r)^2}.
\end{equation}
\end{lemma}
\begin{proof}
Our starting point is a calculation due to Edelman \cite[Prop.~5.1]{edelman1988eigenvalues} that bounds the density $f_{\lambda_{\min}}$ for $\lambda_{\min}(\bm{G} \bm{G}^{\ast})$ as
\begin{equation*}
    f_{\lambda_{\min}}(\lambda)
    \leq \frac{\sqrt{\pi} \Gamma\left(\frac{k+1}{2}\right)}
    {2^{\frac{k - r + 1}{2}} \Gamma\left(\frac{r}{2}\right)
    \Gamma\left(\frac{k-r+1}{2}\right)
    \Gamma\left(\frac{k-r+2}{2}\right)}
    \lambda^\frac{k - r - 1}{2} {\rm e}^{-\frac{\lambda}{2}},
    \qquad \lambda > 0.
\end{equation*}
Applying the Legendre duplication formula $\Gamma(z)\Gamma\bigl(z + \tfrac{1}{2}\bigr) = 2^{1 - 2z} \sqrt{\pi} \Gamma(2z)$, we obtain
\begin{equation}
\label{eq:legendre}
    f_{\lambda_{\min}}(\lambda) \leq \frac{2^{\frac{k-r-1}{2}} \Gamma\left(\frac{k+1}{2}\right)}
    {\Gamma\left(\frac{r}{2}\right) \Gamma\left(k-r+1\right)}
    \lambda^\frac{k - r - 1}{2} \mathrm{e}^{-\frac{\lambda}{2}},
    \qquad \lambda > 0.
\end{equation}
We can simplify \cref{eq:legendre} further by neglecting the $\mathrm{e}^{-\lambda \slash 2}$ factor and noting
\begin{equation*}
    \frac{2^{\frac{k-r+1}{2}} \Gamma\left(\frac{k+1}{2}\right)}
    {\Gamma\left(\frac{r}{2}\right)}
    \leq 2^{\frac{k-r+1}{2}}
    \prod_{i=0}^{k-r}
    \sqrt{\frac{r+i}{2}}
    = \prod_{i=0}^{k-r} \sqrt{r + i}
    \leq \left(\frac{k + r}{2}\right)^{\frac{k-r+1}{2}},
\end{equation*}
where we have used $\Gamma\bigl(z + \tfrac{1}{2}\bigr) \leq z^{1 \slash 2} \Gamma(z)$ \cite{kershaw1983some}
and $(r + i)(k - i) \leq \bigl(\frac{k + r}{2}\bigr)^2$.
We find
\begin{equation}
\label{eq:density}
    f_{\lambda_{\min}}(\lambda) \leq \frac{\lambda^{\frac{k-r-1}{2}}}
    {2\Gamma(k - r + 1)}
    \left(\frac{k + r}{2}\right)^{\frac{k-r+1}{2}},
    \qquad \lambda > 0.
\end{equation}
Integrating \cref{eq:density} with respect to $\lambda$ gives the probability bound
\begin{equation}
\label{eq:prob_bound}
    P\{\lambda_{\min}(\bm{G} \bm{G}^{\ast}) \leq t\} \leq
    \frac{1}
    {\Gamma(k - r + 2)}
    \left(\frac{t(k + r)}{2}\right)^{\frac{k-r+1}{2}}, \qquad t > 0.
\end{equation}
which is nonvacuous for
\begin{equation*}
    t \leq t_{\ast} = \Gamma(k - r + 2)^{\frac{2}{k - r + 1}} \Bigl(\frac{2}{k + r}\Bigr).
\end{equation*}
By an explicit integration of \cref{eq:prob_bound}, we obtain the moment bound
\begin{align*}
    \Expect\lVert (\bm{G} \bm{G}^\ast)^{-1}\rVert^p
    &= p \int_0^{\infty} s^{p-1} P\{\lambda_{\min}(\bm{G} \bm{G}^{\ast}) \leq s^{-1} \} \mathop{ds} \\
    &\leq t_{\ast}^{-p} +  p \int_{t_{\ast}}^{\infty} s^{p-1} P\{\lambda_{\min}(\bm{G} \bm{G}^{\ast}) \leq s^{-1} \} \mathop{ds} \\
    &\leq \Bigl(1 + \frac{2p}{k-r+1-2p}\Bigr) t_{\ast}^{-p} \\
    &= \Bigl(1 + \frac{2p}{k-r+1-2p}\Bigr) \Bigl(\frac{1}{\Gamma(k - r + 2)}\Bigr)^{\frac{2p}{k - r + 1}} \Bigl(\frac{k+r}{2}\Bigr)^p.
\end{align*}
The theorem is complete as soon as we can verify the numerical identity
\begin{equation}
\label{eq:verify}
    \Bigl(1 + \frac{2p}{x+1-2p}\Bigr) \Bigl(\frac{1}{\Gamma(x + 2)}\Bigr)^{\frac{2p}{x + 1}}
    \leq \Bigl(\frac{\mathrm{e}}{x}\Bigr)^{2p},
    \qquad 0 \leq p \leq 18, \quad x \geq 2p.
\end{equation}
To establish \cref{eq:verify}, we briefly argue that
\begin{equation}
\label{eq:apply_to}
    x \mapsto \mathrm{e}\, \Gamma(x + 2)^{\frac{1}{x+1}} \Bigl(1 - \frac{2 p}{x+1}\Bigr)^{2p} - x
\end{equation}
is non-decreasing for $x \geq 2 p$, which can be shown by explicitly examining its derivative and applying
\begin{equation*}
    \mathrm{e}\, \Gamma(x + 2)^{\frac{1}{x+1}} \geq x + 1, 
    \qquad \frac{\mathop{d}}{\mathop{dx}}\Bigl(\mathrm{e}\, \Gamma(x + 2)^{\frac{1}{x+1}}\Bigr) \geq 1.
\end{equation*}
Next, using the fact that \eqref{eq:apply_to} is non-decreasing and applying Stirling's approximation \cite[Sec.~6.1.38]{abramowitz1988handbook}, we calculate the lower bound
\begin{multline*}
    \mathrm{e}\, \Gamma(x + 2)^{\frac{1}{x+1}} \Bigl(1 - \frac{2p}{x+1}\Bigr)^{2p} - x
    \geq \mathrm{e}\, \Gamma(2p + 2)^{\frac{1}{2p+1}} \Bigl(1 - \frac{2p}{2p+1}\Bigr)^{2p} - 2p \\
    \geq \Biggl(\frac{2\pi}{(2p + 1)^{1 + \frac{1}{p}}}\Biggr)^{\frac{1}{4p + 2}} (2p + 1) - 2p,
\end{multline*}
The right-hand side is non-negative provided that
\begin{equation*}
    \frac{2\pi}{(2p + 1)^{1 + \frac{1}{p}}} \geq \mathrm{e}^{-2},
\end{equation*}
which holds for any $p \geq 18$.
We conclude that
\begin{equation*}
    \mathrm{e}\, \Gamma(x + 2)^{\frac{1}{x+1}} \Bigl(1 - \frac{2 p}{x+1}\Bigr)^{2 p} - x \geq 0,
    \qquad 0 \leq p \leq 18, \quad x \geq 2 p,
\end{equation*}
which is equivalent to \eqref{eq:verify}, thus completing the proof.
\end{proof}

\subsection{Proof of \texorpdfstring{\cref{prop:convenient}}{Proposition 4.5}} \label{app:convenient}
When we condition on $\bm{H}$, the mapping
$\bm{G} \mapsto \lVert \bm{S} \bm{G} \bm{H}^{\dagger} \rVert_{\rm F}$
is Lipschitz with Lipschitz constant $L = \lVert \bm{S} \rVert_{\rm F} \lVert \bm{H}^{\dagger} \rVert_{\rm F}$, and the conditional expectation satisfies
\begin{equation*}
    \Expect \bigl[ \lVert \bm{S} \bm{G} \bm{H}^{\dagger} \rVert_{\rm F} \,|\, \bm{H} \bigr] \leq 
    \Expect \bigl[ \lVert \bm{S} \bm{G} \bm{H}^{\dagger} \rVert_{\rm F}^2 \,|\, \bm{H} \bigr]^{1 \slash 2}
    = \lVert \bm{S} \rVert_{\rm F} \lVert \bm{H}^{\dagger} \rVert_{\rm F}.
\end{equation*}
Using a concentration inequality for Lipschitz functions of Gaussian random fields \cite[Thm.~4.5.7]{bogachev1998gaussian}, we find for any $u \geq 1$,
\begin{equation*}
    \mathbb{P}\bigl\{ \lVert \bm{S} \bm{G} \bm{H}^{\dagger} \rVert_{\rm F} > \sqrt{u} \lVert \bm{S} \rVert_{\rm F} \lVert \bm{H}^{\dagger} \rVert_{\rm F} \mid \bm{H} \bigr\}
    \leq \mathrm{e}^{-(\sqrt{u}-1)^2 \slash 2},
\end{equation*}
and the right-hand side is bounded by $\mathrm{e}^{-(u - 2) \slash 4}$ for any $u \geq 0$.
Hence, we obtain the conditional probability bound
\begin{equation}
\label{eq:bound_1}
    \mathbb{P} \left\{\lVert \bm{S} \bm{G} \bm{H}^{\dagger} \rVert_{\rm F}^2
    > \frac{utr}{k-r+1} \left\lVert \bm{S}
    \right\rVert_{\textup{F}}^2
    \,\bigg|\,
    \lVert \bm{H}^{\dagger} \rVert_{\rm F}^2 \leq \frac{tr}{k - r + 1}
    \right\} \leq \mathrm{e}^{-(u - 2) \slash 4},
\end{equation}
while \cref{lem:inverse_moment} implies the probability bound
\begin{equation}
\label{eq:bound_2}
    \mathbb{P}\left\{ \lVert \bm{H}^{\dagger} \rVert_{\rm F}^2 > \frac{tr}{k - r + 1} \right\} \leq \sqrt{\pi r} \biggl(\frac{t}{\mathrm{e}}\biggr)^{-(k-r+1) \slash 2}.
\end{equation}
Combining \cref{eq:bound_1,eq:bound_2} establishes the probability bound \cref{eq:probability_bound} for $\lVert \bm{S} \bm{G} \bm{H}^{\dagger} \rVert^2_{\rm F}$.
Additionally, \cref{lem:moment,lem:inverse_moment} allow us to verify the expectation bound \cref{eq:expectation_bound} for $\lVert \bm{S} \bm{G} \bm{H}^{\dagger} \rVert^2_{\rm F}$, namely,
\begin{align*}
    \Expect \lVert \bm{S} \bm{G} \bm{H}^{\dagger} \rVert_{\rm F}^2
    &= \Expect \Expect \bigl[ \lVert \bm{S} \bm{G} \bm{H}^{\dagger} \rVert_{\rm F}^2 \mid \bm{H}^{\dagger} \bigr]
    = \lVert \bm{S} \rVert_{\rm F}^2 \Expect \lVert \bm{H}^{\dagger} \rVert_{\rm F}^2 \\
    &= \lVert \bm{S} \rVert_{\rm F}^2 \Expect \textup{tr}( \bm{H}^{\dagger} (\bm{H}^{\dagger})^\ast)
    = \lVert \bm{S} \rVert_{\rm F}^2 \textup{tr}(\Expect (\bm{H} \bm{H}^\ast)^{-1}) \\
    &= \lVert \bm{S} \rVert_{\rm F}^2 \frac{r}{k-r-1}.
\end{align*}
In the above display, we have applied conditioning and then exchanged the expectation and the trace.

\subsection{Proof of \texorpdfstring{\cref{thm:more_bounds}}{Theorem 4.6}} \label{app:more_bounds}

Using \cref{lem:diagonal}, we assume $\bm{A}$ is diagonal and psd, with non-increasing diagonal entries.
Then, we partition the matrices as
\begin{equation*}
    \bm{A} = \begin{bmatrix} \bm{A}_1 & \\
    & \bm{A}_2 \end{bmatrix},
    \qquad 
    \bm{\Omega} = \begin{bmatrix}
    \bm{\Omega}_1 \\ \bm{\Omega}_2
    \end{bmatrix},
\end{equation*}
so that $\bm{A}_1$ is a $r \times r$ matrix and $\bm{\Omega}_1$ is a $r \times k$ matrix.
We will systematically derive expectation bounds in the Frobenius norm, spectral norm, and Schatten-4 norm.

To begin,
\cref{prop:ext_analysis} guarantees the error bound
\begin{equation*}
    \Expect\lVert \bm{A} - \bm{\Pi}_{\bm{A} \bm{\Omega}} \bm{A} \rVert_{\rm F}^2 \leq \lVert \bm{A}_2 \rVert_{\rm F}^2
    + \Expect \bigl[\lVert \bm{A}_1 \rVert_{\rm F}^2:
    \lVert \bm{A}_2 \bm{\Omega}_2 \bm{\Omega}_1^{\dagger} \rVert_{\rm F}^2\bigr].
\end{equation*}
Using the concavity of the parallel sum (\cref{lem:parallel_sum} part 3) and \cref{lem:moment}, we argue
\begin{align*}
    \Expect [\lVert \bm{A}_1 \rVert_{\rm F}^2:
    \lVert \bm{A}_2 \bm{\Omega}_2 \bm{\Omega}_1^{\dagger} \rVert_{\rm F}^2]
    &= \Expect \bigl[
    \Expect \bigl[\lVert \bm{A}_1 \rVert_{\rm F}^2:
    \lVert \bm{A}_2 \bm{\Omega}_2 \bm{\Omega}_1^{\dagger} \rVert_{\rm F}^2 \,\big|\, \bm{\Omega}_1 \bigr]\bigr] \\
    &\leq \Expect
    \bigl[\lVert \bm{A}_1 \rVert_{\rm F}^2:
    \Expect \bigl[ \lVert \bm{A}_2 \bm{\Omega}_2 \bm{\Omega}_1^{\dagger} \rVert_{\rm F}^2 \,\big|, \bm{\Omega}_1 \bigr]\bigr] \\
    & = \Expect
    \bigl[\lVert \bm{A}_1 \rVert_{\rm F}^2:
    \lVert \bm{A}_2 \rVert_{\rm F}^2 \lVert \bm{\Omega}_1^{\dagger} \rVert_{\rm F}^2 \bigr] \\
    & = \lVert \bm{A}_2 \rVert_{\rm F}^2\cdot
    \Expect
    \biggl[ \frac{\lVert \bm{A}_1 \rVert_{\rm F}^2}{\lVert \bm{A}_2 \rVert_{\rm F}^2}:
    \lVert \bm{\Omega}_1^{\dagger} \rVert_{\rm F}^2 \biggr].
\end{align*}
Using \cref{lem:inverse_moment}, we bound the last quantity as
\begin{equation*}
    \Expect
    \biggl[ \frac{\lVert \bm{A}_1 \rVert_{\rm F}^2}{\lVert \bm{A}_2 \rVert_{\rm F}^2}:
    \lVert \bm{\Omega}_1^{\dagger} \rVert_{\rm F}^2 \biggr]
    \leq \begin{cases}
    \dfrac{r}{k - r - 1}, & k \geq r + 2, \\[10pt]
    r \log\biggl(1 + 
    \dfrac{\lVert \bm{A}_1 \rVert_{\rm F}^2}
    {\lVert \bm{A}_2 \rVert_{\rm F}^2}\biggr),
    & k = r + 1, \\[10pt]
    r \biggl(
    \dfrac{{\rm \pi} \lVert \bm{A}_1 \rVert_{\rm F}^2}
    {2 \lVert \bm{A}_2 \rVert_{\rm F}^2} \biggr)^{1 \slash 2},
    & k = r,
    \end{cases}
\end{equation*}
which confirms the Frobenius--norm error bounds for \texttt{RSVD}.

Next, working in the spectral norm, \cref{prop:ext_analysis} guarantees
\begin{equation*}
    \Expect\lVert \bm{A} - \bm{\Pi}_{\bm{A} \bm{\Omega}} \bm{A} \rVert^2 \leq \lVert \bm{A}_2 \rVert^2
    + \Expect\lVert \bm{A}_2 \bm{\Omega}_2 \bm{\Omega}_1^{\dagger} \rVert^2
\end{equation*}
We use \cref{lem:moment,lem:inverse_moment,lem:moment_spectral} to calculate
\begin{align*}
    \bigl(\Expect\lVert \bm{A}_2 \bm{\Omega}_2 \bm{\Omega}_1^{\dagger} \rVert^2\Expect\bigr)^{1 \slash 2}
    &= 
    \bigl(\Expect\bigl[ \Expect\bigl[\lVert \bm{A}_2 \bm{\Omega}_2 \bm{\Omega}_1^{\dagger} \rVert^2 \,\big|\, \bm{\Omega}_1 \bigr]\bigr]\bigr)^{1 \slash 2} \\
    &= \bigl(\Expect \bigl| \lVert \bm{A}_2 \rVert \lVert \bm{\Omega}_1^{\dagger} \rVert_{\rm F}
    + \lVert \bm{\Omega}_1^{\dagger} \rVert \lVert \bm{A}_2 \rVert_{\rm F}  |^2  \bigr)^{1 \slash 2} \\
    &\leq \lVert \bm{A}_2 \rVert \bigl(\Expect \lVert \bm{\Omega}_1^{\dagger} \rVert_{\rm F}^2\bigr)^{1 \slash 2}
    + \bigl(\Expect \lVert \bm{\Omega}_1^{\dagger} \rVert^2\bigr)^{1 \slash 2} \lVert \bm{A}_2 \rVert_{\rm F} \\
    &\leq \biggl(\frac{r}{k - r - 1}\biggr)^{1 \slash 2} \lVert \bm{A}_2 \rVert
    + \biggl(\frac{\mathrm{e}^2 (k + r)}{2 (k - r)^2}\biggr)^{1 \slash 2} \lVert \bm{A}_2 \rVert_{\rm F},
\end{align*}
whence
\begin{align*}
    \Expect\lVert \bm{A} - \bm{\Pi}_{\bm{A} \bm{\Omega}} \bm{A} \rVert^2
    &\leq \lVert \bm{A}_2 \rVert^2
    + \frac{2r}{k - r - 1} 
    \lVert \bm{A}_2 \rVert^2
    + \frac{\mathrm{e}^2 (k + r)}{(k - r)^2} 
    \lVert \bm{A}_2 \rVert_{\rm F}^2 \\
    &\leq \biggl(1 + \frac{2r}{k - r - 1}\biggr) \biggl(\lVert \bm{A}_2 \rVert^2
    + \frac{{\rm e}^2}{k-r} \lVert \bm{A}_2 \rVert^2_{\rm F} \biggr),
\end{align*}
which confirms the spectral--norm error bound for \texttt{RSVD}.

Last, to obtain fourth moment expectation bounds, we apply \cref{prop:convenient} and the triangle inequality to yield
\begin{align*}
    \bigl(\Expect\lVert \bm{A} - \bm{\Pi}_{\bm{A} \bm{\Omega}} \bm{A} \rVert_p^4\bigr)^{1 \slash 2} 
    &\leq \bigl(\Expect \bigl| \lVert \bm{A}_2 \rVert_p^2
    + \lVert \bm{A}_2 \bm{\Omega}_2 \bm{\Omega}_1^{\dagger} \rVert_p^2 \bigr|^2 \bigr)^{1 \slash 2} \\
    &\leq \lVert \bm{A}_2 \rVert_p^2 + \bigl(\Expect\lVert \bm{A}_2 \bm{\Omega}_2 \bm{\Omega}_1^{\dagger} \rVert_p^4\bigr)^{1 \slash 2}
\end{align*}
and we proceed to bound $\bigl(\Expect\lVert \bm{A}_2 \bm{\Omega}_2 \bm{\Omega}_1^{\dagger} \rVert_p^4\bigr)^{1 \slash 2}$ for $p = 4, \infty$.
We use \cref{lem:moment,lem:inverse_moment} to calculate
\begin{align*}
    \bigl(\Expect\lVert \bm{A}_2 \bm{\Omega}_2 \bm{\Omega}_1^{\dagger} \rVert_4^4\bigr)^{1 \slash 2}
    &= \bigl(\Expect\bigl[ \Expect\bigl[ \lVert \bm{A}_2 \bm{\Omega}_2 \bm{\Omega}_1^{\dagger} \rVert_4^4\, \big|\, \bm{\Omega}_1 \bigr]\bigr]\bigr)^{1 \slash 2} \\
    &= \bigl(\Expect\bigl[ 
    \lVert \bm{A}_2 \rVert_4^4 \lVert \bm{\Omega}_1^{\dagger} \rVert_{\rm F}^4
    + \lVert \bm{A}_2 \rVert_4^4 \lVert \bm{\Omega}_1^{\dagger} \rVert_4^4
    + \lVert \bm{A}_2 \rVert_{\rm F}^4 \lVert \bm{\Omega}_1^{\dagger} \rVert_4^4
    \bigr]\bigr)^{1 \slash 2} \\
    &\leq 
    \bigl(\Expect\bigl[ 
    \lVert \bm{\Omega}_1^{\dagger} \rVert_{\rm F}^4
    + \lVert \bm{\Omega}_1^{\dagger} \rVert_4^4\bigr] \bigr)^{1 \slash 2}
    \lVert \bm{A}_2 \rVert_4^2
    + \bigl(\Expect\lVert \bm{\Omega}_1^{\dagger} \rVert_4^4\bigr)^{1 \slash 2}
    \lVert \bm{A}_2 \rVert_{\rm F}^2 \\
    &\leq
    \frac{r+1}{k - r - 3}
    \lVert \bm{A}_2 \rVert_4^2
    + 
    \biggl(\frac{r(k-1)}{(k-r)(k-r-1)(k-r-3)} \biggr)^{1 \slash 2}
    \lVert \bm{A}_2 \rVert_{\rm F}^2 \\
    &\leq
    \frac{r+1}{k - r - 3}
    \lVert \bm{A}_2 \rVert_4^2
    + 
    \frac{k-2}{(k - r - 3) \sqrt{k - r}}
    \lVert \bm{A}_2 \rVert_{\rm F}^2,
\end{align*}
where the last line follows because $r \leq k - 4$.
Thus, we confirm the Schatten-4 norm error bound for \texttt{RSVD}.
We use \cref{lem:moment,lem:inverse_moment,lem:moment_spectral} to calculate
\begin{align*}
    & \bigl(\Expect\lVert \bm{A}_2 \bm{\Omega}_2 \bm{\Omega}_1^{\dagger} \rVert^4\bigr)^{1 \slash 4} \\
    &= \bigl(\Expect\bigl[ \Expect\bigl[ \lVert \bm{A}_2 \bm{\Omega}_2 \bm{\Omega}_1^{\dagger} \rVert^4 \,\big|\, \bm{\Omega}_1 \bigr]\bigr]\bigr)^{1 \slash 4} \\
    &\leq \bigl(\Expect\bigl|  
    \lVert \bm{A}_2 \rVert 
    (\lVert \bm{\Omega}_1^{\dagger} \rVert_{\rm F}^4
    + 2\lVert \bm{\Omega}_1^{\dagger} \rVert_4^4)^{1 \slash 4}
    + \lVert \bm{\Omega}_1^{\dagger} \rVert 
    (\lVert \bm{A}_2 \rVert_{\rm F}^4
    + 2\lVert \bm{A}_2 \rVert_4^4)^{1 \slash 4}
     \bigr|^4\bigr)^{1 \slash 4} \\
    &\leq \bigl(\Expect\bigl[ \lVert \bm{\Omega}_1^{\dagger} \rVert_{\rm F}^4
    + 2\lVert \bm{\Omega}_1^{\dagger} \rVert_4^4 \bigr])^{1 \slash 4} \lVert \bm{A}_2 \rVert
    + 
    3^{1 \slash 4} \bigl(\Expect \lVert \bm{\Omega}_1^{\dagger} \rVert^4 \bigr)^{1 \slash 4} \lVert \bm{A}_2 \rVert_{\rm F} \\
    &\leq \biggl( \frac{r + 2}{k - r - 3}  \biggr)^{1 \slash 2} \lVert \bm{A}_2 \rVert
    + 
    \biggl( \frac{\sqrt{3} \mathrm{e}^2(k + r)}{2(k - r)^2} \biggr)^{1 \slash 2} \lVert \bm{A}_2 \rVert_{\rm F},
\end{align*}
whence
\begin{equation*}
    \bigl(\Expect\lVert \bm{A}_2 \bm{\Omega}_2 \bm{\Omega}_1^{\dagger} \rVert^4\bigr)^{1 \slash 2}
    \leq \frac{2(r + 2)}{k - r - 3} \lVert \bm{A}_2 \rVert^2
    + 
    \frac{\sqrt{3} \mathrm{e}^2(k + r)}{(k - r)^2} \lVert \bm{A}_2 \rVert_{\rm F}^2,
\end{equation*}
which confirms the last error bound for \texttt{RSVD} and completes the proof of \cref{thm:more_bounds}.

\section{Proof of \texorpdfstring{\cref{lem:chebyshev}}{Lemma 5.7}} \label{app:chebyshev}
To confirm part 3 of \cref{lem:chebyshev},
we first observe
\begin{equation}
    2r = \int_0^r 2 \mathop{ds} \leq \int_0^r \frac{2 \mathop{ds}}{1 - s^2}
    = \log\biggl(\frac{1 + r}{1 - r}\biggr),
\end{equation}
and by exponentiating both sides we obtain
$\mathrm{e}^{2r} \leq (1 + r)/(1 - r)$.
Next, for $x \geq 1$, we express the Chebyshev polynomial $T_q(x)$ as
\begin{equation}
\label{eq:explicit}
    2 T_q(x) =
    2 \cosh(\log(x + \sqrt{x^2 - 1}))
    =
    \bigl(x + \sqrt{x^2 - 1}\bigr)^q + 
    \bigl(x + \sqrt{x^2 - 1}\bigr)^{-q},
\end{equation}
where we have used the logarithmic representation $\arcosh x = \log(x + \sqrt{x^2 - 1})$ \cite[Sec.~4.6.21]{abramowitz1988handbook}.
By substituting $x = (1 + r) \slash (1-r)$ into \cref{eq:explicit}, we establish
\begin{equation*}
    \mathrm{e}^{2q\sqrt{r}}
    \leq \biggl(\frac{1 + \sqrt{r}}{1 - \sqrt{r}} \biggr)^q
    \leq \biggl(\frac{1 + \sqrt{r}}{1 - \sqrt{r}}\biggr)^q
    + \biggl(\frac{1 + \sqrt{r}}{1 - \sqrt{r}}\biggr)^{-q}
    = 2 T_q\biggl(\frac{1 + r}{1 - r}\biggr),
\end{equation*}
which gives the stated lower bound for the Chebyshev polynomial $T_q$.

To confirm part 4 of \cref{lem:chebyshev}, we write
\begin{equation}
    T_q^{-1}(x) = \cosh\bigl(\tfrac{1}{q} \arcosh x\bigr)
\end{equation}
and apply the inequalities $\cosh(x) \leq \exp(x^2 \slash 2)$ \cite{abramowitz1988handbook}
and $\arcosh x \leq \log(2x)$,
which are both valid for $x \geq 1$.

Last, to confirm part 5 of \cref{lem:chebyshev}, set $\phi(x) = x T_q(x^p)^2$
for $p = 1 \slash 2$ or $p = 1 \slash 4$.
A direct calculation shows the derivative $\phi^{\prime}(x)$ is increasing for $x \geq 1$
and attains its largest value at the right endpoint of the interval $0 \leq x \leq 1$, which suffices to
show
\begin{equation*}
    \phi(y) \geq \phi(x) + \phi^\prime(x) (y - x)
\end{equation*}
for any $x \geq 1$ and any $y \geq 0$.

\bibliographystyle{siamplain}
\bibliography{references}

\begin{thebibliography}{100}

\bibitem{abraham2014fast}
{\sc G.~Abraham and M.~Inouye}, {\em Fast principal component analysis of
  large-scale genome-wide data}, PLOS ONE, 9 (2014), pp.~1--5,
  \url{https://doi.org/10.1371/journal.pone.0093766}.

\bibitem{flashpca}
{\sc G.~Abraham, Y.~Qiu, and M.~Inouye}, {\em Flash{PCA}2: {P}rincipal
  component analysis of {B}iobank-scale genotype datasets}.
\newblock \url{https://github.com/gabraham/flashpca}, 2017.

\bibitem{abramowitz1988handbook}
{\sc M.~Abramowitz and I.~A. Stegun}, {\em Handbook of Mathematical Functions
  with Formulas, Graphs, and Mathematical Tables}, National Bureau of
  Standards, 1972.

\bibitem{anderson1969series}
{\sc W.~Anderson and R.~Duffin}, {\em Series and parallel addition of
  matrices}, Journal of Mathematical Analysis and Applications, 26 (1969),
  pp.~576--594, \url{https://doi.org/10.1016/0022-247X(69)90200-5}.

\bibitem{And94}
{\sc T.~Ando}, {\em Majorizations and inequalities in matrix theory}, Linear
  Algebra and its Applications, 199 (1994), pp.~17--67,
  \url{https://doi.org/10.1016/0024-3795(94)90341-7}.

\bibitem{bakshi2020robust}
{\sc A.~Bakshi, N.~Chepurko, and D.~P. Woodruff}, {\em Robust and sample
  optimal algorithms for {PSD} low rank approximation}, in IEEE 61st Annual
  Symposium on Foundations of Computer Science, 2020,
  \url{https://doi.org/10.1109/FOCS46700.2020.00054}.

\bibitem{BCW22}
{\sc A.~Bakshi, K.~L. Clarkson, and D.~P. Woodruff}, {\em Low-rank
  approximation with $1/\epsilon^{1/3}$ matrix-vector products}, in Proceedings
  of the 54th Annual ACM SIGACT Symposium on Theory of Computing, 2022,
  \url{https://doi.org/10.1145/3519935.3519988}.

\bibitem{bakshi2023krylov}
{\sc A.~Bakshi and S.~Narayanan}, {\em Krylov methods are (nearly) optimal for
  low-rank approximation}, arXiv:2304.03191,  (2023),
  \url{https://arxiv.org/abs/2304.03191}.

\bibitem{balabanov2020randomized}
{\sc O.~Balabanov and L.~Grigori}, {\em Randomized gram--schmidt process with
  application to gmres}, SIAM Journal on Scientific Computing, 44 (2022),
  pp.~A1450--A1474, \url{https://doi.org/10.1137/20M138870X}.

\bibitem{balabanov2021randomized}
{\sc O.~Balabanov and L.~Grigori}, {\em Randomized block {Gram-Schmidt} process
  for solution of linear systems and eigenvalue problems}, arXiv:2111.14641,
  (2023), \url{https://arxiv.org/abs/2111.14641}.

\bibitem{berkooz1993proper}
{\sc G.~Berkooz, P.~Holmes, and J.~L. Lumley}, {\em The proper orthogonal
  decomposition in the analysis of turbulent flows}, Annual Review of Fluid
  Mechanics, 25 (1993), pp.~539--575,
  \url{https://doi.org/10.1146/annurev.fl.25.010193.002543}.

\bibitem{bhatia2013matrix}
{\sc R.~Bhatia}, {\em Matrix Analysis}, Springer, 1997,
  \url{https://doi.org/10.1007/978-1-4612-0653-8}.

\bibitem{bhatia2009positive}
{\sc R.~Bhatia}, {\em Positive Definite Matrices}, Princeton University Press,
  2007, \url{https://doi.org/10.1515/9781400827787}.

\bibitem{bjarkason2019pass}
{\sc E.~K. Bjarkason}, {\em Pass-efficient randomized algorithms for low-rank
  matrix approximation using any number of views}, SIAM Journal on Scientific
  Computing, 41 (2019), pp.~A2355--A2383,
  \url{https://doi.org/10.1137/18M118966X}.

\bibitem{bjarkason2018randomized}
{\sc E.~K. Bjarkason, O.~J. Maclaren, J.~P. O'Sullivan, and M.~J. O'Sullivan},
  {\em Randomized truncated {SVD} {L}evenberg-{M}arquardt approach to
  geothermal natural state and history matching}, Water Resources Research, 54
  (2018), pp.~2376--2404, \url{https://doi.org/10.1002/2017WR021870}.

\bibitem{bogachev1998gaussian}
{\sc V.~Bogachev}, {\em Gaussian Measures}, American Mathematical Society,
  1998, \url{https://doi.org/10.1090/surv/062}.

\bibitem{BKK+19}
{\sc A.~Bose, V.~Kalantzis, E.-M. Kontopoulou, M.~Elkady, P.~Paschou, and
  P.~Drineas}, {\em {TeraPCA}: {A} fast and scalable software package to study
  genetic variation in tera-scale genotypes}, Bioinformatics, 35 (2019),
  pp.~3679--3683, \url{https://doi.org/10.1093/bioinformatics/btz157}.

\bibitem{BDN15}
{\sc J.~Bourgain, S.~Dirksen, and J.~Nelson}, {\em Toward a unified theory of
  sparse dimensionality reduction in euclidean space}, Geometric and Functional
  Analysis, 25 (2015), pp.~1009--1088,
  \url{https://doi.org/10.1007/s00039-015-0332-9}.

\bibitem{bousserez2020enhanced}
{\sc N.~Bousserez, J.~J. Guerrette, and D.~K. Henze}, {\em Enhanced
  parallelization of the incremental {4D-V}ar data assimilation algorithm using
  the randomized incremental optimal technique}, Quarterly Journal of the Royal
  Meteorological Society, 146 (2020), pp.~1351--1371,
  \url{https://doi.org/10.1002/qj.3740}.

\bibitem{boutsidis2015spectral}
{\sc C.~Boutsidis, A.~Gittens, and P.~Kambadur}, {\em Spectral clustering via
  the power method - {P}rovably}, in Proceedings of the 32nd International
  Conference on Machine Learning, 2015,
  \url{https://dl.acm.org/doi/10.5555/3045118.3045124}.

\bibitem{BMD09}
{\sc C.~Boutsidis, M.~W. Mahoney, and P.~Drineas}, {\em An improved
  approximation algorithm for the column subset selection problem}, in
  Proceedings of the 2009 Annual ACM-SIAM Symposium on Discrete Algorithms,
  2009, \url{https://doi.org/10.1137/1.9781611973068.105}.

\bibitem{bui2012extreme}
{\sc T.~Bui-Thanh, C.~Burstedde, O.~Ghattas, J.~Martin, G.~Stadler, and L.~C.
  Wilcox}, {\em Extreme-scale {UQ} for {B}ayesian inverse problems governed by
  {PDE}s}, in Proceedings of the International Conference on High Performance
  Computing, Networking, Storage and Analysis, 2012,
  \url{https://doi.org/10.1109/SC.2012.56}.

\bibitem{chen2013improved}
{\sc C.-Y. Chen, S.~Pollack, D.~J. Hunter, J.~N. Hirschhorn, P.~Kraft, and
  A.~L. Price}, {\em Improved ancestry inference using weights from external
  reference panels}, Bioinformatics, 29 (2013), pp.~1399--1406,
  \url{https://doi.org/10.1093/bioinformatics/btt144}.

\bibitem{chen2022randomly}
{\sc Y.~Chen, E.~N. Epperly, J.~A. Tropp, and R.~J. Webber}, {\em Randomly
  pivoted {C}holesky: {P}ractical approximation of a kernel matrix with few
  entry evaluations}, arXiv:2207.06503,  (2023),
  \url{https://arxiv.org/abs/2207.06503}.

\bibitem{cheng2016fast}
{\sc J.~Cheng and M.~D. Sacchi}, {\em Fast dual-domain reduced-rank algorithm
  for {3D} deblending via randomized {QR} decomposition}, Geophysics, 81
  (2016), pp.~V89--V101, \url{https://doi.org/10.1190/geo2015-0292.1}.

\bibitem{ClWo09}
{\sc K.~L. Clarkson and D.~P. Woodruff}, {\em Numerical linear algebra in the
  streaming model}, in Proceedings of the Forty-First Annual ACM Symposium on
  Theory of Computing, 2009, \url{https://doi.org/10.1145/1536414.1536445}.

\bibitem{CMSW79}
{\sc A.~K. Cline, C.~B. Moler, G.~W. Stewart, and J.~H. Wilkinson}, {\em An
  estimate for the condition number of a matrix}, SIAM Journal on Numerical
  Analysis, 16 (1979), pp.~368--375, \url{http://www.jstor.org/stable/2156842}
  (accessed 2023-06-12).

\bibitem{cohen2016nearly}
{\sc M.~B. Cohen}, {\em Nearly tight oblivious subspace embeddings by trace
  inequalities}, in Proceedings of the 2016 Annual ACM-SIAM Symposium on
  Discrete Algorithms, 2016,
  \url{https://doi.org/10.1137/1.9781611974331.ch21}.

\bibitem{constantine2015active}
{\sc P.~G. Constantine}, {\em Active Subspaces}, Society for Industrial and
  Applied Mathematics, 2015, \url{https://doi.org/10.1137/1.9781611973860}.

\bibitem{DV06}
{\sc A.~Deshpande and S.~Vempala}, {\em Adaptive sampling and fast low-rank
  matrix approximation}, in Approximation, Randomization, and Combinatorial
  Optimization. Algorithms and Techniques, 2006,
  \url{https://doi.org/10.1007/11830924_28}.

\bibitem{dixon1983estimating}
{\sc J.~D. Dixon}, {\em Estimating extremal eigenvalues and condition numbers
  of matrices}, SIAM Journal on Numerical Analysis, 20 (1983), pp.~812--814,
  \url{https://doi.org/10.1137/0720053}.

\bibitem{DFK+99}
{\sc P.~Drineas, A.~Frieze, R.~Kannan, S.~Vempala, and V.~Vinay}, {\em
  Clustering in large graphs and matrices}, in Proceedings of the Tenth Annual
  ACM-SIAM Symposium on Discrete Algorithms, 1999,
  \url{https://dl.acm.org/doi/10.5555/314500.314576}.

\bibitem{drineas2018structural}
{\sc P.~Drineas, I.~C.~F. Ipsen, E.-M. Kontopoulou, and M.~Magdon-Ismail}, {\em
  Structural convergence results for approximation of dominant subspaces from
  block {K}rylov spaces}, SIAM Journal on Matrix Analysis and Applications, 39
  (2018), pp.~567--586, \url{https://doi.org/10.1137/16M1091745}.

\bibitem{drineas2006fast2}
{\sc P.~Drineas, R.~Kannan, and M.~W. Mahoney}, {\em Fast {M}onte {C}arlo
  algorithms for matrices {II}: {C}omputing a low-rank approximation to a
  matrix}, SIAM Journal on Computing, 36 (2006), pp.~158--183,
  \url{https://doi.org/10.1137/S0097539704442696}.

\bibitem{DMM06}
{\sc P.~Drineas, M.~W. Mahoney, and S.~Muthukrishnan}, {\em Subspace sampling
  and relative-error matrix approximation: {C}olumn-based methods}, in
  Approximation, Randomization, and Combinatorial Optimization. Algorithms and
  Techniques, 2006, \url{https://doi.org/10.1007/11830924_30}.

\bibitem{edelman1988eigenvalues}
{\sc A.~Edelman}, {\em Eigenvalues and condition numbers of random matrices},
  SIAM Journal on Matrix Analysis and Applications, 9 (1988), pp.~543--560,
  \url{https://doi.org/10.1137/0609045}.

\bibitem{epperly2023efficient}
{\sc E.~N. Epperly and J.~A. Tropp}, {\em Efficient error and variance
  estimation for randomized matrix computations}, arXiv:2207.06342,  (2023),
  \url{https://arxiv.org/abs/2207.06342}.

\bibitem{epperly2023xtrace}
{\sc E.~N. Epperly, J.~A. Tropp, and R.~J. Webber}, {\em {XT}race: {M}aking the
  most of every sample in stochastic trace estimation}, arXiv:2301.07825,
  (2023), \url{https://arxiv.org/abs/2301.07825}.

\bibitem{Foho13}
{\sc S.~Foucart and H.~Rauhut}, {\em A Mathematical Introduction to Compressive
  Sensing}, Springer, 2013, \url{https://doi.org/10.1007/978-0-8176-4948-7}.

\bibitem{frieze1998fast}
{\sc A.~Frieze, R.~Kannan, and S.~Vempala}, {\em Fast {M}onte-{C}arlo
  algorithms for finding low-rank approximations}, in Proceedings of the 39th
  Annual Symposium on Foundations of Computer Science, 1998,
  \url{https://doi.org/10.1109/SFCS.1998.743487}.

\bibitem{galinsky2016fast}
{\sc K.~J. Galinsky, G.~Bhatia, P.-R. Loh, S.~Georgiev, S.~Mukherjee, N.~J.
  Patterson, and A.~L. Price}, {\em Fast principal-component analysis reveals
  convergent evolution of {ADH1B} in {E}urope and {E}ast {A}sia}, The American
  Journal of Human Genetics, 98 (2016), pp.~456--472,
  \url{https://doi.org/10.1016/j.ajhg.2015.12.022}.

\bibitem{gilmore2008lie}
{\sc R.~Gilmore}, {\em Lie Groups, Physics, and Geometry: An Introduction for
  Physicists, Engineers and Chemists}, Cambridge University Press, 2008,
  \url{https://doi.org/10.1017/CBO9780511791390}.

\bibitem{GHR21}
{\sc A.~Glielmo, B.~E. Husic, A.~Rodriguez, C.~Clementi, F.~Noé, and A.~Laio},
  {\em Unsupervised learning methods for molecular simulation data}, Chemical
  Reviews, 121 (2021), pp.~9722--9758,
  \url{https://doi.org/10.1021/acs.chemrev.0c01195}.

\bibitem{Gol65}
{\sc G.~Golub}, {\em Numerical methods for solving linear least squares
  problems}, Numerische Mathematik, 7 (1965), pp.~206--216,
  \url{https://doi.org/10.1007/bf01436075}.

\bibitem{GoKa65}
{\sc G.~Golub and W.~Kahan}, {\em Calculating the singular values and
  pseudo-inverse of a matrix}, Journal of the Society for Industrial and
  Applied Mathematics Series B Numerical Analysis, 2 (1965), pp.~205--224,
  \url{https://doi.org/10.1137/0702016}.

\bibitem{golub2013matrix}
{\sc G.~Golub and C.~F.~V. Loan}, {\em Matrix Computations}, Johns Hopkins
  University Press, 4~ed., 2013, \url{https://doi.org/10.56021/9781421407944}.

\bibitem{GLO81}
{\sc G.~H. Golub, F.~T. Luk, and M.~L. Overton}, {\em A block {L}anczos method
  for computing the singular values and corresponding singular vectors of a
  matrix}, ACM Transactions on Mathematical Software, 7 (1981), p.~149–169,
  \url{https://doi.org/10.1145/355945.355946}.

\bibitem{golub2000eigenvalue}
{\sc G.~H. Golub and H.~A. {van der Vorst}}, {\em Eigenvalue computation in the
  20th century}, Journal of Computational and Applied Mathematics, 123 (2000),
  pp.~35--65, \url{https://doi.org/10.1016/S0377-0427(00)00413-1}.

\bibitem{goto2008anatomy}
{\sc K.~Goto and R.~A. v.~d. Geijn}, {\em Anatomy of high-performance matrix
  multiplication}, ACM Transactions on Mathematical Software, 34 (2008),
  \url{https://doi.org/10.1145/1356052.1356053}.

\bibitem{gu2015subspace}
{\sc M.~Gu}, {\em Subspace iteration randomization and singular value
  problems}, SIAM Journal on Scientific Computing, 37 (2015), pp.~A1139--A1173,
  \url{https://doi.org/10.1137/130938700}.

\bibitem{GuEi96}
{\sc M.~Gu and S.~C. Eisenstat}, {\em Efficient algorithms for computing a
  strong rank-revealing qr factorization}, SIAM Journal on Scientific
  Computing, 17 (1996), pp.~848--869, \url{https://doi.org/10.1137/0917055}.

\bibitem{sun93}
{\sc J.~guang Sun}, {\em Backward perturbation analysis of certain
  characteristic subspaces}, Numerische Mathematik, 65 (1993), pp.~357--382,
  \url{https://doi.org/10.1007/bf01385757}.

\bibitem{halko2011algorithm}
{\sc N.~Halko, P.-G. Martinsson, Y.~Shkolnisky, and M.~Tygert}, {\em An
  algorithm for the principal component analysis of large data sets}, SIAM
  Journal on Scientific Computing, 33 (2011), pp.~2580--2594,
  \url{https://doi.org/10.1137/100804139}.

\bibitem{halko2011finding}
{\sc N.~Halko, P.~G. Martinsson, and J.~A. Tropp}, {\em Finding structure with
  randomness: {P}robabilistic algorithms for constructing approximate matrix
  decompositions}, SIAM Review, 53 (2011), pp.~217--288,
  \url{https://doi.org/10.1137/090771806}.

\bibitem{hie2019efficient}
{\sc B.~Hie, B.~Bryson, and B.~Berger}, {\em Efficient integration of
  heterogeneous single-cell transcriptomes using {S}canorama}, Nature
  Biotechnology, 37 (2019), pp.~685--691,
  \url{https://doi.org/10.1038/s41587-019-0113-3}.

\bibitem{horn2012matrix}
{\sc R.~A. Horn and C.~R. Johnson}, {\em Matrix Analysis}, Cambridge University
  Press, 2~ed., 2012, \url{https://doi.org/10.1017/CBO9781139020411}.

\bibitem{international2010integrating}
{\sc {International HapMap 3 Consortium}}, {\em Integrating common and rare
  genetic variation in diverse human populations}, Nature, 467 (2010),
  pp.~52--58, \url{https://doi.org/10.1038/nature09298}.

\bibitem{isaac2015scalable}
{\sc T.~Isaac, N.~Petra, G.~Stadler, and O.~Ghattas}, {\em Scalable and
  efficient algorithms for the propagation of uncertainty from data through
  inference to prediction for large-scale problems, with application to flow of
  the {A}ntarctic ice sheet}, Journal of Computational Physics, 296 (2015),
  pp.~348--368, \url{https://doi.org/10.1016/j.jcp.2015.04.047}.

\bibitem{kershaw1983some}
{\sc D.~Kershaw}, {\em Some extensions of {W.} {G}autschi's inequalities for
  the gamma function}, Mathematics of Computation, 41 (1983), pp.~607--611,
  \url{https://doi.org/10.2307/2007697}.

\bibitem{kuczynski1992estimating}
{\sc J.~Kuczy\'{n}ski and H.~Wo\'{z}niakowski}, {\em Estimating the largest
  eigenvalue by the power and {L}anczos algorithms with a random start}, SIAM
  Journal on Matrix Analysis and Applications, 13 (1992), pp.~1094--1122,
  \url{https://doi.org/10.1137/0613066}.

\bibitem{kumar2019target}
{\sc R.~Kumar, M.~Graff, I.~Vasconcelos, and F.~J. Herrmann}, {\em
  Target-oriented imaging using extended image volumes: a low-rank
  factorization approach}, Geophysical Prospecting, 67 (2019), pp.~1312--1328,
  \url{https://doi.org/10.1111/1365-2478.12779}.

\bibitem{lanczos1950iteration}
{\sc C.~Lanczos}, {\em An iteration method for the solution of the eigenvalue
  problem of linear differential and integral operators}, Journal of Research
  of the National Bureau of Standards, 45 (1950).

\bibitem{larsen1998lanczos}
{\sc R.~M. Larsen}, {\em Lanczos bidiagonalization with partial
  reorthogonalization}, DAIMI Report Series, 27 (1998),
  \url{https://doi.org/10.7146/dpb.v27i537.7070}.

\bibitem{lee2014large}
{\sc J.~Lee and P.~K. Kitanidis}, {\em Large-scale hydraulic tomography and
  joint inversion of head and tracer data using the principal component
  geostatistical approach ({PCGA})}, Water Resources Research, 50 (2014),
  pp.~5410--5427, \url{https://doi.org/10.1002/2014WR015483}.

\bibitem{lehoucq1998arpack}
{\sc R.~B. Lehoucq, D.~C. Sorensen, and C.~Yang}, {\em ARPACK Users' Guide},
  Society for Industrial and Applied Mathematics, 1998,
  \url{https://doi.org/10.1137/1.9780898719628}.

\bibitem{levis2021inference}
{\sc A.~Levis, D.~Lee, J.~A. Tropp, C.~F. Gammie, and K.~L. Bouman}, {\em
  Inference of black hole fluid-dynamics from sparse interferometric
  measurements}, in IEEE/CVF International Conference on Computer Vision, 2021,
  \url{https://doi.org/10.1109/ICCV48922.2021.00234}.

\bibitem{li2017algorithm}
{\sc H.~Li, G.~C. Linderman, A.~Szlam, K.~P. Stanton, Y.~Kluger, and
  M.~Tygert}, {\em Algorithm 971: {A}n implementation of a randomized algorithm
  for principal component analysis}, ACM Transactions on Mathematical Software,
  43 (2017), \url{https://doi.org/10.1145/3004053}.

\bibitem{li2015convergence}
{\sc R.-C. Li and L.-H. Zhang}, {\em Convergence of the block {L}anczos method
  for eigenvalue clusters}, Numerische Mathematik, 131 (2014), pp.~83--113,
  \url{https://doi.org/10.1007/s00211-014-0681-6}.

\bibitem{low2016analytical}
{\sc T.~M. Low, F.~D. Igual, T.~M. Smith, and E.~S. Quintana-Orti}, {\em
  Analytical modeling is enough for high-performance {BLIS}}, ACM Transactions
  on Mathematical Software, 43 (2016), \url{https://doi.org/10.1145/2925987}.

\bibitem{martinsson2020randomized}
{\sc P.-G. Martinsson and J.~A. Tropp}, {\em Randomized numerical linear
  algebra: {F}oundations and algorithms}, Acta Numerica, 29 (2020),
  p.~403–572, \url{https://doi.org/10.1017/S0962492920000021}.

\bibitem{mason2002chebyshev}
{\sc J.~Mason and D.~C. Handscomb}, {\em Chebyshev Polynomials}, Chapman and
  Hall/{CRC}, 2002, \url{https://doi.org/10.1201/9781420036114}.

\bibitem{matlab2020version}
{\sc {MathWorks Inc.}}, {\em {MATLAB} {V}ersion 9.9 ({R}2010a)}, 2022,
  \url{https://www.mathworks.com}.

\bibitem{meyer2023unreasonable}
{\sc R.~A. Meyer, C.~Musco, and C.~Musco}, {\em On the unreasonable
  effectiveness of single vector {K}rylov methods for low-rank approximation},
  arXiv:2305.02535,  (2023), \url{https://arxiv.org/abs/2305.02535}.

\bibitem{musco2015randomized}
{\sc C.~Musco and C.~Musco}, {\em Randomized block {K}rylov methods for
  stronger and faster approximate singular value decomposition}, in Proceedings
  of the 28th International Conference on Neural Information Processing
  Systems, 2015, \url{https://dl.acm.org/doi/10.5555/2969239.2969395}.

\bibitem{nakatsukasa2020fast}
{\sc Y.~Nakatsukasa}, {\em Fast and stable randomized low-rank matrix
  approximation}, arXiv:2009.11392,  (2020),
  \url{https://arxiv.org/abs/2009.11392}.

\bibitem{nakatsukasa2021fast}
{\sc Y.~Nakatsukasa and J.~A. Tropp}, {\em Fast \& accurate randomized
  algorithms for linear systems and eigenvalue problems}, arXiv:2111.00113,
  (2021), \url{https://arxiv.org/abs/2111.00113}.

\bibitem{ala2}
{\sc F.~N{\"u}ske, H.~Wu, J.-H. Prinz, C.~Wehmeyer, C.~Clementi, and
  F.~No{\'{e}}}, {\em Alanine dipeptide}.
\newblock \url{https://markovmodel.github.io/mdshare/ALA2/}.

\bibitem{NWPW+17}
{\sc F.~N{\"u}ske, H.~Wu, J.-H. Prinz, C.~Wehmeyer, C.~Clementi, and F.~Noé},
  {\em Markov state models from short non-equilibrium simulations—{A}nalysis
  and correction of estimation bias}, The Journal of Chemical Physics, 146
  (2017), \url{https://doi.org/10.1063/1.4976518}.

\bibitem{OSV79}
{\sc D.~P. O'Leary, G.~W. Stewart, and J.~S. Vandergraft}, {\em Estimating the
  largest eigenvalue of a positive definite matrix}, Mathematics of
  Computation, 33 (1979), pp.~1289--1292,
  \url{https://doi.org/10.2307/2006463}.

\bibitem{orourke2018random}
{\sc S.~O'Rourke, V.~Vu, and K.~Wang}, {\em Random perturbation of low rank
  matrices: Improving classical bounds}, Linear Algebra and its Applications,
  540 (2018), pp.~26--59, \url{https://doi.org/10.1016/j.laa.2017.11.014}.

\bibitem{papadimitriou2000latent}
{\sc C.~H. Papadimitriou, P.~Raghavan, H.~Tamaki, and S.~Vempala}, {\em Latent
  semantic indexing: {A} probabilistic analysis}, Journal of Computer and
  System Sciences, 61 (2000), pp.~217--235,
  \url{https://doi.org/10.1006/jcss.2000.1711}.

\bibitem{PTRV98}
{\sc C.~H. Papadimitriou, H.~Tamaki, P.~Raghavan, and S.~Vempala}, {\em Latent
  semantic indexing: A probabilistic analysis}, in Proceedings of the
  Seventeenth ACM SIGACT-SIGMOD-SIGART Symposium on Principles of Database
  Systems, 1998, \url{https://doi.org/10.1145/275487.275505}.

\bibitem{parlett1998symmetric}
{\sc B.~N. Parlett}, {\em The Symmetric Eigenvalue Problem}, Society for
  Industrial and Applied Mathematics, 1998,
  \url{https://doi.org/10.1137/1.9781611971163}.

\bibitem{scikit-learn}
{\sc F.~Pedregosa, G.~Varoquaux, A.~Gramfort, V.~Michel, B.~Thirion, O.~Grisel,
  M.~Blondel, P.~Prettenhofer, R.~Weiss, V.~Dubourg, J.~Vanderplas, A.~Passos,
  D.~Cournapeau, M.~Brucher, M.~Perrot, and {{\'E}}douard Duchesnay}, {\em
  Scikit-learn: {M}achine learning in {P}ython}, Journal of Machine Learning
  Research, 12 (2011), pp.~2825--2830,
  \url{http://jmlr.org/papers/v12/pedregosa11a.html}.

\bibitem{PiWa15}
{\sc M.~Pilanci and M.~J. Wainwright}, {\em Randomized sketches of convex
  programs with sharp guarantees}, IEEE Transactions on Information Theory, 61
  (2015), pp.~5096--5115, \url{https://doi.org/10.1109/TIT.2015.2450722}.

\bibitem{press2005applied}
{\sc S.~J. Press}, {\em Applied Multivariate Analysis: {U}sing {B}ayesian and
  Frequentist Methods of Inference}, Robert E. Krieger Publishing Co, 2~ed.,
  1982.

\bibitem{rokhlin2010randomized}
{\sc V.~Rokhlin, A.~Szlam, and M.~Tygert}, {\em A randomized algorithm for
  principal component analysis}, SIAM Journal on Matrix Analysis and
  Applications, 31 (2010), pp.~1100--1124,
  \url{https://doi.org/10.1137/080736417}.

\bibitem{Row97}
{\sc S.~Roweis}, {\em Em algorithms for pca and spca}, in Proceedings of the
  10th International Conference on Neural Information Processing Systems, 1997,
  \url{https://dl.acm.org/doi/10.5555/3008904.3008993}.

\bibitem{rudelson2007sampling}
{\sc M.~Rudelson and R.~Vershynin}, {\em Sampling from large matrices: {A}n
  approach through geometric functional analysis}, Journal of the ACM, 54
  (2007), p.~21–es, \url{https://doi.org/10.1145/1255443.1255449}.

\bibitem{saad2011numerical}
{\sc Y.~Saad}, {\em Numerical Methods for Large Eigenvalue Problems}, Society
  for Industrial and Applied Mathematics, 2011,
  \url{https://doi.org/10.1137/1.9781611970739}.

\bibitem{Sar06}
{\sc T.~S{\'a}rlos}, {\em Improved approximation algorithms for large matrices
  via random projections}, in 2006 47th Annual IEEE Symposium on Foundations of
  Computer Science, 2006, \url{https://doi.org/10.1109/FOCS.2006.37}.

\bibitem{schmid2022dynamic}
{\sc P.~J. Schmid}, {\em Dynamic mode decomposition and its variants}, Annual
  Review of Fluid Mechanics, 54 (2022), pp.~225--254,
  \url{https://doi.org/10.1146/annurev-fluid-030121-015835}.

\bibitem{smith2014anatomy}
{\sc T.~M. Smith, R.~v.~d. Geijn, M.~Smelyanskiy, J.~R. Hammond, and F.~G.~V.
  Zee}, {\em Anatomy of high-performance many-threaded matrix multiplication},
  in IEEE 28th International Parallel and Distributed Processing Symposium,
  2014, \url{https://doi.org/10.1109/IPDPS.2014.110}.

\bibitem{stewart1976simultaneous}
{\sc G.~W. Stewart}, {\em Simultaneous iteration for computing invariant
  subspaces of non-{H}ermitian matrices}, Numerische Mathematik, 25 (1976),
  pp.~123--136, \url{https://doi.org/10.1007/bf01462265}.

\bibitem{stewart1991perturbation}
{\sc G.~W. Stewart}, {\em Perturbation theory for the singular value
  decomposition}, in SVD and Signal Processing {II}: {A}lgorithms, Analysis and
  Applications, R.~J. Vaccaro, ed., Elsevier, 1991, pp.~99--109,
  \url{http://hdl.handle.net/1903/552}.

\bibitem{stewart2001matrix}
{\sc G.~W. Stewart}, {\em Matrix Algorithms}, vol.~II: Eigensystems, Society
  for Industrial and Applied Mathematics, 2001,
  \url{https://doi.org/10.1137/1.9780898718058}.

\bibitem{tropp2011improved}
{\sc J.~A. Tropp}, {\em Improved analysis of the subsampled randomized
  {H}adamard transform}, Advances in Adaptive Data Analysis, 03 (2011),
  pp.~115--126, \url{https://doi.org/10.1142/S1793536911000787}.

\bibitem{tropp2018analysis}
{\sc J.~A. Tropp}, {\em Analysis of randomized block {K}rylov methods}, ACM
  Report 2018-02, Caltech, Pasadena, 2018,
  \url{https://resolver.caltech.edu/CaltechAUTHORS:20210624-180721369}.

\bibitem{tropp2022randomized}
{\sc J.~A. Tropp}, {\em Randomized block {K}rylov methods for approximating
  extreme eigenvalues}, Numerische Mathematik, 150 (2021), pp.~217--255,
  \url{https://doi.org/10.1007/s00211-021-01250-3}.

\bibitem{tropp2017fixed}
{\sc J.~A. Tropp, A.~Yurtsever, M.~Udell, and V.~Cevher}, {\em Fixed-rank
  approximation of a positive-semidefinite matrix from streaming data}, in
  Proceedings of the 31st International Conference on Neural Information
  Processing Systems, 2017,
  \url{https://dl.acm.org/doi/10.5555/3294771.3294888}.

\bibitem{tropp2017practical}
{\sc J.~A. Tropp, A.~Yurtsever, M.~Udell, and V.~Cevher}, {\em Practical
  sketching algorithms for low-rank matrix approximation}, SIAM Journal on
  Matrix Analysis and Applications, 38 (2017), pp.~1454--1485,
  \url{https://doi.org/10.1137/17M1111590}.

\bibitem{TYUC19}
{\sc J.~A. Tropp, A.~Yurtsever, M.~Udell, and V.~Cevher}, {\em Streaming
  low-rank matrix approximation with an application to scientific simulation},
  SIAM Journal on Scientific Computing, 41 (2019), pp.~A2430--A2463,
  \url{https://doi.org/10.1137/18M1201068}.

\bibitem{tsuyuzaki2020benchmarking}
{\sc K.~Tsuyuzaki, H.~Sato, K.~Sato, and I.~Nikaido}, {\em Benchmarking
  principal component analysis for large-scale single-cell {RNA}-sequencing},
  Genome Biology, 21 (2020), \url{https://doi.org/10.1186/s13059-019-1900-3}.

\bibitem{VSM11}
{\sc V.~Vanhoucke, A.~Senior, and M.~Z. Mao}, {\em Improving the speed of
  neural networks on {CPU}s}, in Deep Learning and Unsupervised Feature
  Learning Workshop, NIPS 2011, 2011,
  \url{https://research.google/pubs/pub37631/}.

\bibitem{vershynin2018high}
{\sc R.~Vershynin}, {\em High-Dimensional Probability: An Introduction with
  Applications in Data Science}, Cambridge University Press, 2018,
  \url{https://doi.org/10.1017/9781108231596}.

\bibitem{wallace1992singular}
{\sc J.~M. Wallace, C.~Smith, and C.~S. Bretherton}, {\em Singular value
  decomposition of wintertime sea surface temperature and 500-mb height
  anomalies}, Journal of Climate, 5 (1992), pp.~561 -- 576,
  \url{https://doi.org/10.1175/1520-0442(1992)005<0561:SVDOWS>2.0.CO;2}.

\bibitem{wang2015improved}
{\sc S.~Wang, Z.~Zhang, and T.~Zhang}, {\em Improved analyses of the randomized
  power method and block {L}anczos method}, arXiv:1508.06429,  (2015),
  \url{https://arxiv.org/abs/1508.06429}.

\bibitem{Wedin_1972}
{\sc P.-{\AA}. Wedin}, {\em Perturbation bounds in connection with singular
  value decomposition}, {BIT}, 12 (1972), pp.~99--111,
  \url{https://doi.org/10.1007/bf01932678}.

\bibitem{williams2001using}
{\sc C.~K.~I. Williams and M.~Seeger}, {\em Using the {N}ystr\"{o}m method to
  speed up kernel machines}, in Proceedings of the 13th International
  Conference on Neural Information Processing Systems, 2000,
  \url{https://dl.acm.org/doi/10.5555/3008751.3008847}.

\bibitem{woolfe2008fast}
{\sc F.~Woolfe, E.~Liberty, V.~Rokhlin, and M.~Tygert}, {\em A fast randomized
  algorithm for the approximation of matrices}, Applied and Computational
  Harmonic Analysis, 25 (2008), pp.~335--366,
  \url{https://doi.org/10.1016/j.acha.2007.12.002}.

\bibitem{yang2021low}
{\sc M.~Yang, M.~Graff, R.~Kumar, and F.~J. Herrmann}, {\em Low-rank
  representation of omnidirectional subsurface extended image volumes},
  Geophysics, 86 (2021), pp.~S165--S183,
  \url{https://doi.org/10.1190/geo2020-0152.1}.

\bibitem{yu2018efficient}
{\sc W.~Yu, Y.~Gu, and Y.~Li}, {\em Efficient randomized algorithms for the
  fixed-precision low-rank matrix approximation}, SIAM Journal on Matrix
  Analysis and Applications, 39 (2018), pp.~1339--1359,
  \url{https://doi.org/10.1137/17M1141977}.

\bibitem{yuan2018superlinear}
{\sc Q.~Yuan, M.~Gu, and B.~Li}, {\em Superlinear convergence of randomized
  block {L}anczos algorithm}, in IEEE International Conference on Data Mining,
  2018, \url{https://doi.org/10.1109/ICDM.2018.00193}.

\end{thebibliography}
\end{document}